\documentclass[12pt]{article}
\usepackage{fancyhdr, array,calc,graphicx,url,tabularx}
\usepackage{geometry}
\usepackage{latexsym}
\usepackage{amssymb}
\usepackage{amsmath}
\usepackage{makeidx}
\usepackage{enumerate}
\usepackage{verbatim}
\usepackage{tikz}
\usepackage[colorlinks=true,linkcolor=blue]{hyperref}
\usepackage{changepage}
\usepackage{fancybox}

\usepackage{amsthm}

\usepackage[utf8]{inputenc}
\usepackage{tikz}
\tikzset{
     block/.style={rectangle, draw, fill=red!40, text width=6em,
                   text centered, rounded corners, minimum height=3em},
     arrow/.style={-{Stealth[]}}
     }

\makeindex

{
   \end{minipage}
   \vspace*{\stretch{3}}
   \clearpage
}

\renewcommand{\gg}{\gamma}

\newcommand{\bR}{{\mathbb{R}}}

\newcommand{\rest}{\restriction}

\newcommand{\card}[1]{{\vert #1 \vert} }

\renewcommand{\models}{\vDash}
\newcommand{\powerset}{{\wp}}
\newcommand{\dom}{{\rm dom}}
\newcommand{\rge}{{\rm rge}}
\newcommand{\cp}{{\rm crit }}

\newcommand{\cf}{{\rm cf}}
\newcommand{\lh}{{\rm lh}}

\newtheorem{theorem}{Theorem}[section]
\newtheorem{proposition}[theorem]{Proposition}
\newtheorem{definition}[theorem]{Definition}

\newtheorem{lemma}[theorem]{Lemma}
\newtheorem{corollary}[theorem]{Corollary}
\newtheorem{claim}[theorem]{Claim}

\newtheorem{conjecture}[theorem]{Conjecture}
\newtheorem{open}[theorem]{Open Problem}
\newtheorem{question}[theorem]{Question}

\newtheorem{remark}[theorem]{Remark}

\numberwithin{figure}{section}

\newcommand{\rcon}[1]{Conjecture~\ref{#1}}

\newcommand{\rprop}[1]{Proposition~\ref{#1}}
\newcommand{\rthm}[1]{Theorem~\ref{#1}}
\newcommand{\rlem}[1]{Lemma~\ref{#1}}

\newcommand{\rcor}[1]{Corollary~\ref{#1}}
\newcommand{\rdef}[1]{Definition~\ref{#1}}

\newcommand{\rsec}[1]{Section~\ref{#1}}

\newcommand{\rsubsec}[1]{Section~\ref{#1}}
\newcommand{\rrem}[1]{Remark~\ref{#1}}

\def\inseg{\trianglelefteq}

\def\k{\kappa}
\def\a{\alpha}
\def\b{\beta}
\def\d{\delta}

\def\l{\lambda}
\def\m{{\rm m}}

\def\P{{\mathcal{P} }}
\def\W{{\mathcal{W} }}
\def\Q{{\mathcal{ Q}}}
\def\mH{{\mathcal{ H}}}
\def\K{{\mathcal{ K}}}
\def\J{{\mathcal{ J}}}

\def\R{{\mathcal R}}

\def\H{{\rm{HOD}}}
\def\M{{\mathcal{M}}}
\def\N{{\mathcal{N}}}
\def\F{{\mathcal{F}}}
\def\T {{\mathcal{T}}}
\def\U{{\mathcal{U}}}
\def\S{{\mathcal{S}}}
\def\V{{\mathcal{V}}}
\def\X{{\mathcal{X}}}
\def\Y{{\mathcal{Y}}}

\def\F{{\mathcal{F}}}

\def\VT{{\vec{\mathcal{T}}}}
\def\VU{{\vec{\mathcal{U}}}}

\def\card#1{\left|#1\right|}

\def\iff{\mathrel{\leftrightarrow}}

\def\and{\mathrel{\kern1pt\&\kern1pt}}

\def\inseg{\triangleleft}
\def\insegeq{\trianglelefteq}

\def\<#1>{\langle\,#1\,\rangle}

 \input xy
 \xyoption{all}
       
\newcounter{nameOfYourChoice}

\title{The exact consistency strength of the generic absoluteness for the universally Baire sets}

\author{Grigor Sargsyan\footnote{Institute of Mathematics of Polish Academy of Sciences, Poland. Email: gsargsyan@impan.pl}
\\
Nam Trang \footnote{Department of Mathematics, University of North Texas, Denton, TX, USA. Email: Nam.Trang@unt.edu}}

\date{\today}

\begin{document}
\maketitle

\begin{abstract}

A set of reals is \textit{universally Baire} if all of its continuous preimages in topological spaces have the Baire property.  $\sf{Sealing}$ is a type of generic absoluteness condition introduced by Woodin that asserts in strong terms that the theory of the   universally Baire sets cannot be changed by forcing. 
  
 The $\sf{Largest\ Suslin\ Axiom}$ ($\sf{LSA}$) is a determinacy axiom isolated by Woodin. It asserts that the largest Suslin cardinal is inaccessible for ordinal definable bijections. Let $\sf{LSA-over-uB}$ be the statement that  in all (set) generic extensions there is a model of $\sf{LSA}$ whose Suslin, co-Suslin sets are the universally Baire sets. 
 
 We show that over some mild large cardinal theory, $\sf{Sealing}$ is equiconsistent with $\sf{LSA-over-uB}$. In fact, we isolate an exact large cardinal theory that is equiconsistent with both (see \rdef{dfn:hod_pm}). As a consequence, we obtain that $\sf{Sealing}$ is weaker than  the theory $``\sf{ZFC} + $there is a Woodin cardinal which is a limit of Woodin cardinals". 
 
A variation of $\sf{Sealing}$, called $\sf{Tower \ Sealing}$, is also shown to be equiconsistent with $\sf{Sealing}$ over the same large cardinal theory.
 
 The result is proven via Woodin's $\sf{Core\ Model\ Induction}$ technique, and is essentially the ultimate equiconsistency that can be proven via the current interpretation of $\sf{CMI}$ as explained in the paper. 
%
\end{abstract}

\tableofcontents
\newpage

\section{Introduction}

Soon after Cohen discovered forcing and established the consistency of the failure of the $\sf{Continuum\ Hypothesis \ (CH)}$ with $\sf{ZFC}$, thus establishing the independence of $\sf{CH}$ from $\sf{ZFC}$\footnote{Cohen proved that $\sf{ZFC+\neg CH}$ is consistent. Earlier, G\"odel showed that $\sf{ZFC+CH}$ is consistent by showing that $\sf{CH}$ holds in the \textit{constructible universe} $L$. Forcing can also be used to show that $\sf{ZFC+CH}$ is consistent.}, many natural and useful set theoretic principles have been discovered to remove independence from set theory. Perhaps the two best known ones are $\sf{Shoenfield's\ Absoluteness\ Theorem}$ and $\sf{Martin's\ Axiom}$. 

As is well known, $\sf{Shoenfield's\ Absoluteness\ Theorem}$, proved in \cite{Sho}, asserts that there cannot be any independence result expressible as a $\Sigma^1_2$ fact.  In the language of real analysis, $\Sigma^1_2$ sets of reals are projections of co-analytic sets\footnote{A set of reals is analytic if it is a projection of a closed set. A co-analytic set is the complement of an analytic set.}. Shoenfield's theorem says that a co-analytic set is empty if and only if its natural interpretations\footnote{As open sets are unions of open intervals, it must be clear that they can be easily interpreted in any extension of the reals.} in all generic extensions are empty. What is so wonderful about $\sf{Shoenfield's\ Absoluteness}$ $\sf{Theorem}$ is that it is a theorem of $\sf{ZFC}$. We will discuss $\sf{Martin's\ Axiom}$ and its generalization later on. 

 The goal of this paper is to establish an equiconsistency result between one Shoenfield-type generic absoluteness principle known as $\sf{Sealing}$ and a determinacy axiom that we abbreviated as $\sf{LSA-over-uB}$. $\sf{LSA}$ stands for the \textit{Largest-Suslin-Axiom}.  To state the main theorem, we need a few definitions.

A set of reals is \textit{universally Baire} if all of its continuous preimages in topological spaces have the property of Baire. Let $\Gamma^\infty$ be the collection of universally Baire sets\footnote{The superscript $\infty$ in this notation, which is due to Woodin, makes sense as one can define $\k$-universally Baire sets as those sets whose continuous preimages in all $<\k$ size topological spaces have the property of Baire. We then set $\Gamma^\k$ be the collection of all of these sets. Clearly $\Gamma^\infty=\cap_\k\Gamma^\k$.}.  Given a generic $g$, we let $\Gamma^\infty_g=_{def}(\Gamma^\infty)^{V[g]}$ and $\bR_g=_{def}\bR^{V[g]}$. $\powerset(X)$ is the powerset of $X$. $\sf{AD}$ stands for the $\sf{Axiom\ of\ Determinacy}$ and $\sf{AD}^+$ is a strengthening of $\sf{AD}$ due to Woodin. The reader can ignore the $+$ or can consult \cite[Definition 9.6]{Woodin}.

Motivated by Woodin's $\sf{Sealing\ Theorem}$ (\cite[Theorem 3.4.17]{StationaryTower} and \cite[Sealing Theorem]{SealingTheorem}), we define $\sf{Sealing}$, a key notion in this paper. We say $V[g], V[h]$ are two successive generic extensions (of $V$) if $g, h$ are $V$-generic and $V[g]\subseteq V[h]$.

\begin{definition}\label{dfn:ub_sealing} $\sf{Sealing}$ is the conjunction of the following statements.
\begin{enumerate}
\item For every set generic $g$, $L(\Gamma^\infty_g, \mathbb{R}_g)\models \sf{AD}^+$ and $\powerset(\bR_g)\cap L(\Gamma^\infty_g, \mathbb{R}_g)=\Gamma^\infty_g$.
\item  For every two successive set generic extensions $V[g]\subseteq V[h]$, there is an elementary embedding 
\begin{center}
$j: L(\Gamma^\infty_g, \mathbb{R}_g)\rightarrow L(\Gamma^\infty_h, \mathbb{R}_h)$.
\end{center}
\end{enumerate}
 such that for every $A\in \Gamma^\infty_g$, $j(A)=A_h$\footnote{The meaning of $A_h$ is explained below. It is the canonical extension of $A$ to $V[h]$.}.  
 \end{definition}
 
To introduce $\sf{LSA-over-uB}$, we first need to introduce the $\sf{Largest\ Suslin\ Axiom}$ $(\sf{LSA})$. A cardinal $\k$ is $\sf{OD}$-inaccessible if for every $\a<\k$ there is no surjection $f: \powerset(\a)\rightarrow \k$ that is definable from ordinal parameters. A set of reals $A\subseteq \bR$ is $\k$-\textit{Suslin} if for some tree $T$ on $\k$, $A=p[T]$\footnote{Given a cardinal $\k$, we say $T\subseteq \bigcup_{n<\omega} \omega^n \times \k^n$ is a \textit{tree} on $\k$ if $T$ is closed under initial segments. Given a tree $T$ on $\k$, we let $[T]$ be the set of its branches, i.e., $b\in [T]$ if $b\in \omega^\omega\times \k^\omega$ and letting $b=(b_0, b_1)$, for each $n\in \omega$, $(b_0\rest n, b_1\rest n)\in T$. We then let $p[T]=\{ x\in \bR: \exists f((x, f)\in [T])\}$.}. A set $A$ is \textit{Suslin} if it is $\kappa$-Suslin for some $\kappa$; $A$ is \textit{co-Suslin} if its complement $\mathbb{R}\backslash A$ is Suslin.  A set $A$ is \textit{Suslin, co-Suslin} if both $A$ and its complement are Suslin. A cardinal $\k$ is a \textit{Suslin cardinal} if there is a set of reals $A$ such that $A$ is $\k$-Suslin but $A$ is not $\l$-Suslin for any $\l<\k$. Suslin cardinals play an important role in the study of models of determinacy as can be seen by just flipping through the Cabal Seminar Volumes (\cite{Cabal3}, \cite{Cabal4}, \cite{CabalReprint1}, \cite{CabalReprint2}, \cite{CabalReprint3}, \cite{Cabal2}, \cite{Cabal1}). 

The $\sf{Largest\ Suslin\ Axiom}$ was introduced by Woodin in \cite[Remark 9.28]{Woodin}. The terminology is due to the first author. Here is the definition. 
\begin{definition}\label{dfn:lsa} 
The $\sf{Largest\ Suslin\ Axiom}$, abbreviated as $\sf{LSA}$, is the conjunction of the following statements:
\begin{enumerate}
\item $\sf{AD}^+$.
\item There is a largest Suslin cardinal.
\item The largest Suslin cardinal is $\sf{OD}$-inaccessible.
\end{enumerate}
 \end{definition}
 
 In the hierarchy of determinacy axioms, which one may appropriately call the $\sf{Solovay\ Hierarchy}$\footnote{Solovay defined what is now called the $\sf{Solovay\ Sequence}$ (see \cite[Definition 9.23]{Woodin}). It is a closed sequence of ordinals with the largest element $\Theta$, where $\Theta$ is the least ordinal that is not a surjective image of the reals. One then obtains a hierarchy of axioms by requiring that the $\sf{Solovay\ Sequence}$ has complex patterns. $\sf{LSA}$ is an axiom in this hierarchy. The reader may consult \cite{BSL} or \cite[Remark 9.28]{Woodin}.}, $\sf{LSA}$ is an anomaly as it belongs to the successor stage of the $\sf{Solovay\ Hierarchy}$ but does not conform to the general norms of the successor stages of the $\sf{Solovay\ Hierarchy}$. Prior to \cite{hod_mice_LSA}, $\sf{LSA}$ was not known to be consistent. In \cite{hod_mice_LSA}, the first author showed that it is consistent relative to a Woodin cardinal that is a limit of Woodin cardinals. Nowadays, the axiom plays a key role in many aspects of inner model theory, and features prominently  in Woodin's $\sf{Ultimate\ L}$ framework (see \cite[Definition 7.14]{UltimateL} and Axiom I and Axiom II on page 97 of \cite{UltimateL}\footnote{The requirement in these axioms that there is a strong cardinal which is a limit of Woodin cardinals is only possible if $L(A, \bR)\models \sf{LSA}$.}). 
  
  \begin{definition}\label{dfn:generic_LSA}
Let $\sf{LSA-over-uB}$ be the statement: For all $V$-generic $g$, in $V[g]$, there is $A\subseteq \bR_g$ such that $L(A, \bR_g)\models \sf{LSA}$ and $\Gamma^\infty_g$ is the Suslin co-Suslin sets of $L(A, \bR_g)$.  
\end{definition}

 The following is our main theorem. We say that $\phi$ and $\psi$ are equicosnsitent over theory $T$ if there is a model of $T\cup\{\phi\}$ if and only if there is a model of $T\cup\{\psi\}$.  
 
\begin{theorem}\label{thm:main_theorem}  $\sf{Sealing}$ and $\sf{LSA-over-uB}$ are equiconsistent over  the theory ``there exists a proper class of Woodin cardinals and the class of measurable cardinals is stationary".
 
\end{theorem}

In \rthm{thm:main_theorem}, ``$T_1$ and $T_2$ are equiconsistent" is used in the following stronger sense: there is a well-founded model of $T_1$ if and only if there is a well-founded model of $T_2$. 

 \begin{remark}\label{rmk:weak_sealing}  
It is our intention to consider $\sf{Sealing}$ under large cardinals. The reason for doing this is that universally Baire sets do not in general behave nicely when there are no large cardinals in the universe. One may choose to drop clause 1 from the definition of $\sf{Sealing}$. Call the resulting principle $\sf{Weak\ Sealing}$. If there is an inaccessible cardinal $\k$ which is a limit of Woodin cardinals  and strong cardinals then $\sf{Weak\ Sealing}$ implies $\sf{Sealing}$. This is because one may arrange so that $\Gamma^\infty$ is \textit{the derived model} after Levy collapsing $\k$ to be $\omega_1$ (see \rthm{dmt}). We do not know the consistency strength of $\sf{Weak\ Sealing}$ or $\sf{Sealing}$ in the absence of large cardinals. But one  gets that $\sf{Weak \ Sealing}$ and $\sf{Sealing}$ are equiconsistent over the large cardinal hypothesis in \rthm{thm:main_theorem}. 

 \end{remark}

%
%
%

Based on the above theorems, it is very tempting to conjecture that: $\sf{Sealing}$ and $\sf{LSA-over-uB}$ are equivalent over ``there exists a proper class of Woodin cardinals and the class of measurable cardinals is stationary". However, \cite{Sealing_iter} shows that this conjecture is false. The following variation of $\sf{Sealing}$, called $\sf{Tower \ Sealing}$, is also isolated by Woodin.  

\begin{definition}\label{dfn:tower_sealing}
$\sf{Tower \ Sealing}$ is the conjunction of:
\begin{enumerate}
\item For any set generic $g$, $L(\Gamma^\infty_g) \models \sf{AD}^+$, and $\Gamma^\infty_g = \powerset(\mathbb{R})\cap L(\Gamma^\infty_g,\mathbb{R}_g)$.
\item For any set generic $g$, in $V[g]$, suppose $\delta$ is Woodin. Whenever $G$ is $V[g]$-generic for either the $\mathbb{P}_{<\delta}$-stationary tower or the $\mathbb{Q}_{<\delta}$-stationary tower at $\delta$, then 
\begin{center}
$j(\Gamma^\infty_g) = \Gamma^\infty_{g*G}$,
\end{center}
where $j: V[g] \rightarrow M \subset V[g*G]$ is the generic elementary embedding given by $G$.
\end{enumerate}
\end{definition}

\begin{theorem}\label{thm:tower_sealing}
$\sf{Tower \ Sealing}$ and $\sf{Sealing}$ are equiconsistent over ``there exists a proper class of Woodin cardinals and the class of measurable cardinals is stationary".
\end{theorem}

\begin{remark}
\begin{enumerate}
\item The proof of Theorems \ref{thm:main_theorem} and \ref{thm:tower_sealing} shows that over the large cardinal assumption stated in \rthm{thm:main_theorem}, $\sf{LSA-over-uB}$ and $\sf{Sealing}$ are equiconsistent relative to the following consequence of $\sf{Sealing}$ and of $\sf{Tower \ Sealing}$ (cf. \rprop{adr in gamma ub}): \\ 
 \noindent $\sf{Sealing^-}$: \ \  ``for any set generic $g$, $\Gamma^\infty_g = \powerset(\mathbb{R})\cap L(\Gamma^\infty_g,\mathbb{R}_g)$ and there is no $\omega_1$-sequence of distinct reals in $L(\Gamma^\infty_g,\mathbb{R}_g)$."
 
\item  As mentioned above, \cite{Sealing_iter} shows that $\sf{LSA-over-UB}$ is not equivalent to $\sf{Sealing}$ (over some large cardinal theory). However, the equivalence of $\sf{Sealing},$ $\sf{Tower \ Sealing}$, other weak forms of these theories may still hold (over the large cardinal theory of Theorem \ref{thm:main_theorem}). See Conjecture \ref{equivalence conjecture}.
\item Woodin has observed that assuming a proper class of Woodin cardinals which are limits of strong cardinals, $\sf{Tower Sealing}$ implies $\sf{Sealing}$.
\end{enumerate}
\end{remark}


Before giving the proof, in the next few sections we will explain the context of \rthm{thm:main_theorem}.\\\\
\textbf{Generic Absoluteness}\\
 
 As was mentioned in the opening paragraph, the discovery of forcing almost immediately initiated the study of removing independence phenomenon from set theory. Large cardinals were used to establish a plethora of results that generalize $\sf{Shoenfield's\ Absoluteness\ Theorem}$ to more complex formulas than $\Sigma^1_2$. In another direction, new axioms were discovered that imply what is forced is already true. These axioms are called \textit{forcing axioms}, and $\sf{Martin's\ Axiom}$ is the first one.  
 
Forcing axioms assert that analogues of the $\sf{Baire\ Category\ Theorem}$ hold for any collection of $\aleph_1$-dense sets. A consequence of this is that the $\aleph_1$-fragment of the generic object added by the relevant forcing notion exists as a set in the ground model, implying that what is forced by the $\aleph_1$-fragment of the generic is already true in the ground model.  $\sf{Martin's\ Axiom}$ and its generalizations do not follow from $\sf{ZFC}$. Many axioms of this type have been introduced and extensively studied. Perhaps the best known ones are  $\sf{Martin's\ Axiom}$ (\cite{MA}), the $\sf{Proper\ Forcing\ Axiom}$ ($\sf{PFA}$, see \cite{Baumgartner}) and $\sf{Martin's\ Maximum}$ (see \cite{FMSI}).  

 The general set theoretic theme described above is known as generic absoluteness. The interested reader can consult \cite{VialeII}, \cite{FMSI}, \cite{FMSII}, \cite{Hauser}, \cite{StationaryTower},  \cite{GodelSteel}, \cite{Todorcevic}, \cite{Viale}, \cite{WilsonAbs}, \cite{Woodin} and the references appearing in those papers. We will not be dealing with forcing axioms in this paper, but $\sf{PFA}$ will be used for illustrative purposes.   

The largest class of sets of reals for which a Shoenfield-type generic absoluteness can hold is the collection of the universally Baire sets. We will explain this claim below. The story begins with the fact that the universally Baire sets have canonical interpretations in all generic extensions, and in a sense, they are the only ones that have this property. The next paragraph describes exactly how this happens. The proofs appear in \cite{UB}, \cite{StationaryTower} and \cite{DMT}.

 In \cite{UB}, it was shown by Feng, Magidor and Woodin  that a set of reals $A$ is universally Baire if and only if for each uncountable cardinal $\k$ there are trees $T$ and $S$ on $\k$ such that $p[T]=A$ and in all set generic extensions  $V[g]$ of $V$ obtained by a poset of size $<\k$, $V[g]\models p[T]=\bR-p[S]$. The canonical interpretation of $A$ in $V[g]$ is just $A_g=_{def}(p[T])^{V[g]}$ where $T$ is chosen on a $\kappa$ that is bigger than the size of the poset that adds $g$. It is not hard to show, using the absoluteness of well-foundedness, that $A_g$ is independent of the choice of $(T, S)$. 

Woodin showed that if $A$ is a universally Baire set of reals and the universe has a class of Woodin cardinals then the theory of $L(A, \bR)$ cannot be changed. He achieved this by showing that if there is a class of Woodin cardinals then for any universally Baire set $A$ and any two successive set generic extensions $V[g]\subseteq V[h]$, there is an elementary embedding $j:L(A_g, \bR_g)\rightarrow L(A_h, \bR_h)$\footnote{Unfortunately, the authors do not know a reference for this theorem of Woodin. But it can be proven via the methods of \cite{StationaryTower} and \cite{DMT}.}.  

Moreover, if sufficient generic absoluteness is true about a set of reals then that set is universally Baire. More precisely, suppose $\phi$ is a property of reals. Let $A_{\phi}$ be the set of reals defined by $\phi$. If sufficiently many statements about $A_\phi$ are generically absolute  then it is because $A_\phi$ is universally Baire (see the $\sf{Tree\ Production\ Lemma}$ in \cite{StationaryTower} or in \cite{DMT})\footnote{The exact condition is that club of countable Skolem hulls are generically correct.}. Thus, the next place to look for absoluteness is the set of all universally Baire sets. 

Is it possible that there is no independence result about the set of universally Baire sets? $\sf{Sealing}$, introduced in the preamble of this paper, is the formal version of the English sentence asserting that much like individual universally Baire sets, much like integers, the theory of universally Baire sets is immune to forcing. It is stated in the spirit of Woodin's aforementioned theorem for the individual universally Baire sets.  

While the definition of $\sf{Sealing}$ is very natural and its statement is seemingly benign, $\sf{Sealing}$ has drastic consequences on the $\sf{Inner\ Model\ Program}$, which is one of the oldest set theoretic projects, and is also the next set theoretical theme that we introduce.\\\\
\textbf{The Inner Model Program and The Inner Model Problem}\\

The goal of the Inner Model Program ($\sf{IMP}$) is to build canonical $L$-like inner models with large cardinals. The problem of building a canonical inner model for a large cardinal axiom $\phi$ is known as the $\sf{Inner\ Model\ Problem}$ ($\sf{IMPr}$) for $\phi$. There are several expository articles written about $\sf{IMP}$ and $\sf{IMPr}$. The reader who wants to learn more can consult \cite{JensenBSL}, \cite{BSL}, \cite{ABCMice}.
 
In \cite{Neeman}, Neeman, assuming the existence of a Woodin cardinal that is a limit of Woodin cardinals, solved the $\sf{IMPr}$ for  a Woodin cardinal that is a limit of Woodin cardinals and for large cardinals somewhat beyond. Neeman's result is the best current result on $\sf{IMPr}$.  However, this is only a tiny fragment of the large cardinal paradise, and also \textit{its solution is specific to the hypothesis} (we will discuss this point more).
 
 Dramatically, $\sf{Sealing}$ implies that $\sf{IMP}$, as is known today, cannot succeed as if $\M$ is a model that conforms to the norms of modern inner model theory and has some very basic closure properties then $\M\models ``$there is a well-ordering of reals in $L(\Gamma^\infty, \bR)$". As $\sf{AD}$ implies the reals cannot be well-ordered, $\M$ cannot satisfy $\sf{Sealing}$.  Thus, we must have the following\footnote{$\sf{Sealing\ Dichotomy}$ is well-known among inner model theorists, we do not mean that we were the first to notice it.}.\\\\
 $\sf{Sealing\ Dichotomy}$\\
Either no large cardinal theory implies $\sf{Sealing}$ or the $\sf{Inner\ Model\ Problem}$ for some large cardinal cannot have a solution conforming to the modern norms.\\

Intriguingly,  Woodin, assuming the existence of a supercompact cardinal and a class of Woodin cardinals, has shown that $\sf{Sealing}$ holds after collapsing the powerset of the powerset of a supercompact cardinal to be countable (for a proof, see \cite[Theorem 3.4.17]{StationaryTower}). Because we are collapsing the supercompact to be countable, it seems that Woodin's result does not imply that $\sf{Sealing}$ has dramatic effect on $\sf{IMP}$, or at least this impact cannot be seen in the large cardinal region below supercompact cardinals, which is known as the \textit{short extender region}. 

As part of proving \rthm{thm:main_theorem} and \rthm{thm:tower_sealing}, we will establish that
\begin{theorem}\label{upper bound} $\sf{Sealing}$ is consistent relative to a Woodin cardinal that is a limit of Woodin cardinals. So is $\sf{Tower \ Sealing}$.
\end{theorem}
One consequence of \rthm{upper bound} is that $\sf{Sealing}$ is within the short extender region. While \rthm{upper bound} doesn't illustrate the impact of $\sf{Sealing}$, its exact impact on $\sf{IMP}$ in the short extender region can also be precisely stated. But to do this we will need extenders.\\\\
\textbf{Extender Detour}\\

Before we go on, let us take a minute to introduce extenders, which are natural generalization of ultrafilters. In fact, extenders are just a coherent sequence of ultrafilters. As was mentioned above, the goal of $\sf{IMP}$ is to build canonical $L$-like inner models for large cardinals. The current methodology is that such models should be constructed in G\"odel's sense from extenders, the very objects whose existence large cardinal axioms assert. Perhaps the best way to introduce extenders is via the elementary embeddings that they induce. 

Suppose $M$ and $N$ are two transitive models of set theory and $j:M\rightarrow N$ is a non-trivial elementary embedding. Let $\kappa=\cp(j)$ and let $\l\in [\k, j(\k))$ be any ordinal. Set 
\begin{center}
$E_j=\{(a, A)\in [\l]^{<\omega}\times \powerset([\k]^{\card{a}})^M: a\in j(A)\}$.
\end{center}
$E_j$ is called the $(\k, \l)$-extender derived from $j$. $E_j$ is really an $M$-extender as it measures the sets in $M$. As with more familiar ultrafilters, one can define extenders abstractly without using the parent embedding $j$, and then show that each extender, via an ultrapower construction, gives rise to an embedding. Given a $(\k, \l)$-extender $E$ over $M$, we let $\pi_E: M\rightarrow Ult(M, E)$ be the ultrapower embedding.  A computation that involves chasing the definitions shows that $E$ is the extender derived from $\pi_E$. Similar computations also show that $\k=\cp(\pi_E)$ and $\pi_E(\k)\geq \l$. It is customary to write $\cp(E)$ for $\k$ and $lh(E)=\l$\footnote{``$lh(E)$ is the length of $E$". We note that when discussing Mitchell-Steel extender models, $lh(E)$ is the cardinal successor of the natural length of $E$. The natural length of $E$ is the supremum of generators of $E$. For more details, see \cite{steel2010outline}. }. It is also not hard to see that for each $a\in [\l]^{<\omega}$, $E_a$ is an ultrafilter concentrating on $[\k]^{\card{a}}$, and that if $a\subseteq b$ then $E_b$ \textit{naturally} projects to $E_a$. 

 The motivation behind extenders is the fact that extenders capture more of the universe in the ultrapower than one can achieve via the usual ultrapower construction. In particular, under large cardinal assumptions, one can have $(\k, \l)$-extender $E$ such that $V_\l\subseteq Ult(V, E)$. Because of this all large cardinal notions below superstrong cardinals can be captured by extenders.
 
 The extenders as we defined them above are called \textit{short} extenders, where shortness refers to the fact that all of its ultrafilters concentrate on its critical point. Large cardinal notions such as supercompactness, hugeness and etc cannot be captured by such short extenders as embeddings witnessing supercompactness give rise to measures that do not concentrate on the critical point of the embedding. However, one can capture these large cardinal notions by using the so-called \textit{long} extenders. We do not need them in this paper, and so we will not dwell on them. 
 
The large cardinal region that can be captured by  short extenders is the region of \textit{superstrong cardinals}. A cardinal $\k$ is called superstrong if there is an embedding $j:V\rightarrow M$ with $\cp(j)=\k$ and $V_{j(\k)}\subseteq M$. Superstrong cardinals are close to the optimal cardinal notions that can be expressed via short extenders.
 
Currently, to solve $\sf{IMPr}$ for a large cardinal, one attempts to build a model of the form $L[\vec{E}]$ where $\vec{E}$ is a carefully chosen sequence of extenders. The reader interested in learning more about what $L[\vec{E}]$ should be can consult \cite{steel2010outline}. This ends our detour.\\\\
\textbf{The Core Model Induction}\\

What is a solution to the $\sf{IMPr}$ for a given large cardinal? In the short extender region, $\sf{IMPr}$ for a large cardinal notion such as superstrong cardinals has a somewhat precise meaning. One is essentially asked to build a model of the form $L[\vec{E}]$ which has a superstrong cardinal and $\vec{E}$ is a \textit{fine extender sequence} as defined in \cite[Definition 2.4]{steel2010outline}. However, one may do this construction under many different hypotheses. 

As was mentioned above, Neeman solved the $\sf{IMPr}$ for a Woodin cardinal that is a limit of Woodin cardinals assuming the existence of such a cardinal. One very plausible precise interpretation of $\sf{IMPr}$ is exactly in this sense. Namely given a large cardinal axiom $\phi$, assuming large cardinals that are possibly stronger than $\phi$, build an $\M=L[\vec{E}]$ such that $\M\models \exists \k \phi(\k)$.  

Our interpretation of $\sf{IMPr}$ is influenced by John Steel's view on G\"odel's\ Program (see \cite{GodelSteel}). In a nutshell, the idea is to develop a theory that connects various foundational frameworks such as Forcing Axioms, Large Cardinals, Determinacy Axioms etc with one another\footnote{Our goal here is to avoid philosophical discussions, but if we were to go in this direction we would call this view approach to $\sf{IMPr}$ Steel's Program.}. In this view, $\sf{IMPr}$ is the bridge between all of these natural frameworks and $\sf{IMPr}$ needs to be solved under variety of hypotheses, such as $\sf{PFA}$ or failure of Jensen's $\square$ principles. Our primary tool for solving $\sf{IMPr}$ in large-cardinal-free contexts is the $\sf{Core\ Model\ Induction}$ ($\sf{CMI}$), which is a technique invented by Woodin and developed by many set theorists during the past 20-25 years\footnote{In some contexts, $K^c$ theory can also be used. See \cite{jensen2009stacking}. But solving $\sf{IMPr}$ via a $K^c$ theory won't in general provide such bridges between frameworks. The $K^c$ approach will not in general connect say $\sf{PFA}$ with the $\sf{Solovay\ Hierarchy}$. See \rcon{pfa conjecture}.}.

In the earlier days, $\sf{CMI}$ was perceived as an inductive method for proving determinacy in models such as $L(\mathbb{R})$. The goal was to prove that $L_\a(\bR)\models \sf{AD}$ by induction on $\a$. In those earlier days, which is approximately the period 1995-2010, the method worked by establishing intricate connections between large cardinals, universally Baire sets and determinacy\footnote{For example, the reader may try to understand the meaning of $W^*_\a$ in \cite{PFA}.}. The fundamental work done by Jensen, Neeman, Martin, Mitchell, Steel and Woodin were, and still are, at the heart of current developments of $\sf{CMI}$. The following is a non-exhaustive list of influential papers: most papers in the Cabal Seminar Volumes that discuss scales or playful universes (\cite{Cabal3}, \cite{Cabal4}, \cite{CabalReprint1}, \cite{CabalReprint2}, \cite{CabalReprint3}, \cite{Cabal2}, \cite{Cabal1}), \cite{Jensen}, \cite{IT}, \cite{PD}, \cite{FSIT}, \cite{neeman1995optimal}, \cite{CMIP}. Several fundamental papers were written implicitly developing this view of $\sf{CMI}$. For example, the reader can consult \cite{Ketchersid}, \cite{PFA} and \cite{cmi}. As $\sf{CMI}$ evolved, it became more of a tool for \textit{deriving maximal determinacy models} from non-large cardinal hypotheses. 

In a seminal work, Woodin has developed a technique for \textit{deriving} determinacy models from large cardinals. The theorem is known as the $\sf{Derived\ Model\ Theorem}$. A typical situation works as follows. Suppose $\l$ is a limit of Woodin cardinals and $g\subseteq Coll(\omega, <\l)$ is generic. Let $\bR^*=\cup_{\a<\l}\bR^{V[g\cap Coll(\omega, \a)]}$. Working in $V(\mathbb{R}^*)$\footnote{This is the minimal transitive model $W$ of $\sf{ZF}$ such that $V\subseteq W$ and $\bR^*\in W$. It can be shown that $\bR^W=\bR^*$.} let $\Gamma=\{A\subseteq \bR: L(A, \bR)\models \sf{AD}\}$. Then 
\begin{theorem}[Woodin, \cite{DMT}]\label{dmt} $L(\Gamma, \bR)\models \sf{AD}$.
\end{theorem}
In Woodin's theorem, $\Gamma$ is \textit{maximal} as there are no more (strongly) determined sets in the universe that are not in $\Gamma$. If one assumes that $\l$ is a limit of strong cardinals then $\Gamma$ above is just $(\Gamma^\infty)^{V(\bR^*)}$.

The aim of $\sf{CMI}$ is to do the same for other natural set theoretic frameworks, such as forcing axioms, combinatorial statements etc. Suppose $T$ is a natural set theoretic framework and $V\models T$. Let $\kappa$ be an uncountable cardinal. One way to perceive $\sf{CMI}$ is the following. \\\\
($\sf{CMI}$ at $\k$) Saying that one is doing $\sf{Core\ Model\ Induction}$ at $\k$ means that for some $g\subseteq Coll(\omega, \k)$\footnote{This is the poset that collapses $\k$ to be countable.}, in $V[g]$, one is proving that $L(\Gamma^\infty, \bR)\models \sf{AD}^+$.\\
($\sf{CMI}$ below $\k$) Saying that one is doing $\sf{Core\ Model\ Induction}$ below $\k$ means that for some $g\subseteq Coll(\omega, <\k)$\footnote{This is the poset that collapses everything $<\k$ to be countable.}, in $V[g]$, one is proving that $L(\Gamma^\infty, \bR)\models \sf{AD}^+$.\\\\
In both cases, the aim might be less ambitious. It might be that one's goal is to just produce $\Gamma\subseteq \Gamma^\infty$ such that $L(\Gamma, \bR)$ is a determinacy model with desired properties. 

$\sf{CMI}$ can even help proving versions of the $\sf{Derived\ Model\ Theorem}$. Here is an example. We use the set up introduced for stating the  $\sf{Derived\ Model\ Theorem}$ (\rthm{dmt}). Suppose $A\in V(\bR^*)$ is a set of reals such for some $\a<\l$ there is a $<\l$-universally Baire set $B\in V[g\cap Coll(\omega, \a)]$ such that $A=B_g$. The tools developed during the earlier period of $\sf{CMI}$ can be used to show that $L(A, \bR^*)\models \sf{AD}$. The point here is just that $\sf{CMI}$ is the most general method for proving derived model type of results. In the  $\sf{Derived\ Model\ Theorem}$, the presence of large cardinals makes $\sf{CMI}$ unnecessary, but in other cases it is the only method we currently have. One can also attempt to prove the full $\sf{Derived\ Model\ Theorem}$ via $\sf{CMI}$, but this seems harder and some of the main technical difficulties associated with other non-large cardinal frameworks resurface. 

The goal, however, is not to just derive a determinacy model from natural set theoretic frameworks, but to establish that the determinacy model has the same set theoretic \textit{complexity} as $V$ has. 

 Let $M$ be the maximal model of determinacy derived from $V$. One natural\footnote{That this way of stating the desired closeness is \textit{natural} is a consequence of several decades of research carried out on $\H$ of models of determinacy. See the fragment of the introduction titled \textbf{HOD Analysis}.}  way of saying that $M$ has the same complexity as $V$ is by saying that the large cardinal complexity of $V$ is reflected into $M$, and one particularly elegant way of saying this is to say that $\H^M$, the universe of the hereditarily ordinal definable sets of $M$, acquires these large cardinals. A typical conjecture that we can now state in this language is as follows. 

\begin{conjecture}\label{pfa conjecture} Assume the $\sf{Proper\ Forcing\ Axiom}$ and suppose $\k\geq \omega_2$. Let $g\subseteq Coll(\omega, \k)$. Then $\H^{L(\Gamma^\infty_g, \bR_g)}\models ``$there is a superstrong cardinal".
\end{conjecture}

A less ambitious conjecture would be that $\sf{PFA}$ implies that whenever $g\subseteq Coll(\omega, \k)$ is $V$-generic, there is a set of reals $A\in \Gamma^\infty_g$ such that $\H^{L(A, \bR_g)}\models ``$there is a superstrong cardinal". However, we believe that the stronger conjecture is also true. One can change $\sf{PFA}$ to any other natural framework that is expected to be stronger than superstrong cardinals.\footnote{In some cases, we work in $V[g]$ for $g\subseteq Coll(\omega,\kappa)$ for some $\kappa$. In other cases, we may work in $V[g]$ for $g\subseteq Coll(\omega, <\kappa)$. Whether one does $\sf{CMI}$ at $\kappa$ or below $\kappa$ is hypothesis dependent.} As we brought up $\H$, it is perhaps important to discuss its use in $\sf{CMI}$.\\\\
\textbf{HOD Analysis and Covering}\\

\rcon{pfa conjecture} is a product of many decades of work that goes back to the UCLA's Cabal Seminar, where the study of \textit{playful universes} originates (see, for example,  \cite{YannisII} and \cite{YannisI}). Our attempt is to avoid a historical introduction to the subject, and so we will avoid the long history of studying $\H$ and its \textit{playful} inner models assuming determinacy. 

Nowadays, we know that $\H$ of many models satisfying $\sf{AD}$ is an $L$-like model carrying many large cardinals,\footnote{See for example \cite{hod_as_core_model},\cite{hod_mice},\cite{hod_mice_LSA}.} and the problem of showing that $\H$ of every model of $\sf{AD}$ is an $L$-like model is one of the central open problems of \textit{descriptive inner model theory} (see \cite{BSL} and \cite{normalization_comparison}). 

The current methodology for proving that $\H^{L(\Gamma^\infty, \bR)}$ has the desired large cardinals is via a failure of certain covering principle involving $\H^{L(\Gamma^\infty, \bR)}$. Recall that under determinacy, $\Theta$ is defined to be the least ordinal that is not a surjective image of the reals. Set $\mH^-=\H^{L(\Gamma^\infty, \bR)}|\Theta$. 

To define the aforementioned covering principle, we first need to extend $\mH^-$ to a model $\mH$ in which $\Theta$ is the largest cardinal. This is a standard construction in inner model theory. We simply let $\mH$ be the union of all \textit{hod mice} extending $\mH$ whose countable submodels have iteration strategies in $L(\Gamma^\infty, \bR)$. This sentence perhaps means little to a general reader. It turns out, however, that in many situations it is possible to describe $\mH$ without any reference to inner model theoretic objects.

Here is one such example. Suppose $\k$ is a measurable cardinal which is a limit of strong cardinals and suppose we are doing $\sf{CMI}$ below $\k$. Let $g\subseteq Coll(\omega, <\k)$ be our generic. We also make the assumption that all sets of reals produced by the $\sf{CMI}$ below $\kappa$ are universally Baire\footnote{This is why we assume that $\kappa$ is a limit of strong cardinals, as this hypothesis implies what we stated.}. Let $j:V\rightarrow M$ be any embedding with $\cp(j)=\k$. We furthermore assume that we have succeeded in showing that $\sup(j[\Theta])<j(\Theta)$\footnote{This condition happens quite often}. Setting $\nu=\sup(j[\Theta])$,  let $\mathcal{C}(\mH^-)$ be the set of all $A\subseteq \Theta$ such that $j(A)\cap \nu\in j(\mH^-)$. Then $\mH$ is the transitive model extending $\mH^-$ that is coded by the elements of $\mathcal{C}(\mH^-)$\footnote{Fix a pairing function $\pi:\Theta^2\rightarrow \Theta$. Given $A\subseteq \mathcal{C}(\mH^-)$ we say $A$ is a code if $M_A=(\Theta, E_A)$ is a well-founded model where $E_A\subseteq \Theta^2$ is given by $(\a, \b)\in E_A\iff \pi(\a, \b)\in A$. If $A\in \mathcal{C}(\mH^-)$ is a code then let $\M_A$ be the transitive collapse of $M_A$. Then $\mH$ is the union of models of the form $\M_A$.}. 

At any rate, one can simply regard $\mH$ as a canonical one-cardinal extension of $\mH^-$. In fact, that $\mH$ is a canonical extension of $\mH^-$ is the central point. The next paragraph explains this.

Continuing with the above scenario, let now $h\subseteq Coll(\omega, \kappa)$ be $V[g]$-generic\footnote{Recall that above we were doing CMI below $\k$ and $g\subseteq Coll(\omega, <\k)$. Also, $\kappa=\omega_1^{V[g]}$.}. Because $\card{V_\k}=\k$, we have that $\card{\mH^-}^{V[g*h]}=\aleph_0$ and $\card{\mH}^{V[g*h]}\leq \aleph_1$. Letting $\eta=Ord\cap \mH$, 
\begin{center}
$L(\Gamma^\infty_{g*h}, \bR_{g*h})\models ``$there is an $\eta$-sequence of distinct reals".
\end{center}
Assuming $\sf{Sealing}$, we get that $\eta<\omega_1$ as under $\sf{Sealing}$, $L(\Gamma^\infty_{g*h}, \bR_{g*h})\models \sf{AD}$, and under $\sf{AD}$ there is no $\omega_1$-sequence of reals. Therefore, in $V$, $\eta<\kappa^+$ as we have that $(\kappa^+)^V=\omega_1^{V[g*h]}$. Letting now\\\\
$\sf{UB-Covering:}$ $\cf^V(Ord\cap \mH)\geq \k$,\\\\
$\sf{Sealing}$ implies that $\sf{UB-Covering}$ fails at measurable cardinals that are limits of strong cardinals. A similar argument can be carried out by only assuming that $\k$ is a singular strong limit cardinal. \footnote{In this case, $\mathcal{H}$ is defined in $V(\mathbb{R}^*)$, where $\mathbb{R}^* = \bigcup_{\alpha<\kappa} \mathbb{R}^{V[h\cap Col(\omega,\alpha)]}$ and $h\subseteq Coll(\omega,<\kappa)$ is $V$-generic.}

All other sufficiently strong frameworks also imply that the $\sf{UB-Covering}$ fails but for different reasons. One particular reason is that $\sf{UB-Covering}$ implies that Jensen's $\square_\k$ holds at singular cardinal $\kappa$, while a celebrated theorem of Todorcevic says that under $\sf{PFA}$, $\square_\k$ has to fail for all $\k\geq \omega_2$. 

The argument that has been used to show that $\mH$ has large cardinals proceeds as follows. Pick a target large cardinal $\phi$, which for technical reasons we assume is a $\Sigma_2$-formula. Assume $\mH\models \forall \gamma \neg \phi(\gamma)$. Thus far, in all applications of the $\sf{CMI}$, the facts that\\\\
$\sf{\phi-Minimality}:$ $\mH\models \forall \gamma \neg \phi(\gamma)$\\ 
and\\
$\sf{\neg\ UB-Covering}$: cf$^V({\mH\cap Ord})<\k$\\\\
hold have been used to prove that there is a universally Baire set not in $\Gamma^\infty_g$ where $g\subseteq Coll(\omega, \k)$ or $g\subseteq Coll(\omega, <\k)$ (depending where we do $\sf{CMI}$), which is obviously a contradiction. 

Because of the work done in the first 15 years of the 2000s, it seemed as though this is a general pattern that will persist through the short extender region. That is, for any $\phi$ that is in the short extender region, either $\sf{\phi-Minimality}$ must fail or $\sf{UB-Covering}$ must hold. \textit{The main way \rthm{upper bound} affects $\sf{IMPr}$ in the short extender region is by implying that this prevalent view is false.} \footnote{Example of $\phi(\gamma)$ is: ``$\gamma$ is a Woodin cardinal which is a limit of Woodin cardinals".}

Almost all existing literature on $\sf{CMI}$ uses the argument outlined above in one way or another. The interested reader can consult \cite{Dominik},  \cite{hod_mice},  \cite{hod_mice_LSA}, \cite{UBH},\cite{PFA}, \cite{Trang2015PFA}, \cite{Zhu}.  \\\\
\textbf{The future of} $\sf{CMI}$\\

In the authors' view, $\sf{CMI}$ should be viewed as a technique for proving that certain type of covering holds rather than a technique for showing that $\H$ has large cardinals. The latter should be the corollary, not the goal. The type of covering that we have in mind is the following. We state it in the short extender region, and we use $\sf{NLE}$\footnote{``No Long Extender"} of \cite{normalization_comparison} to state that we are in the short extender region.

\begin{conjecture}\label{main conjecture} Assume $\sf{NLE}$ and suppose there are unboundedly many Woodin cardinals and strong cardinals. Let $\kappa$ be a limit of Woodin cardinals and strong cardinals such that either $\k$ is a measurable cardinal or has cofinality $\omega$. Then there is a transitive model $M$ of $\sf{ZFC-Powerset}$ such that
\begin{enumerate}
\item $Ord\cap M= \k^+$,
\item $M$ has a largest cardinal $\nu$,
\item for any $V$-generic $g\subseteq Coll(\omega, <\k)$, letting $\mathbb{R}^* = \bigcup_{\alpha<\k} \mathbb{R}^{V[g\cap Coll(\omega,\alpha)]}$, in $V(\bR^*)$, 
\begin{center}
$L(M, \bigcup_{\a<\nu}(M|\a)^\omega, \Gamma^\infty, \bR^*)\models \sf{AD}$.
\end{center}
\item If in addition there is no inner model with a subcompact cardinal then $M\models \square_\nu$.
\end{enumerate}
\end{conjecture}

In \cite{CCMDM}, the first author showed that \rcon{main conjecture} holds in hod mice. Assuming $V$ is a hod mouse and keeping the notation of \rcon{main conjecture}, $M$ is simply the direct limit of all iterates of $V|\kappa^+$ that are below $\kappa$ and have a countable length in $V[g]$.

With more work, the conjecture can also be stated without assuming the large cardinals. We do not believe that the conjecture is true in the long extender region because of the following general argument. Assume $\k$ is an indestructible  supercompact cardinal and suppose the conclusion of the conjecture holds at $\k$. Let $g\subseteq Coll(\k, \k^+)$ be $V$-generic. Then \textit{presumably} if $M$ satisfies the conclusion of \rcon{main conjecture}, then $M^V=M^{V[g]}$. The confidence that this is true comes from the fact that we expect that any $M$ satisfying clause 3 must have an absolute definition. Because $\k$ is still a supercompact in $V[g]$, clause 1 has to fail.

We believe that proving \rcon{main conjecture} should become the goal of $\sf{CMI}$. To prove it, one has to develop techniques for building \textit{third order} canonical objects, objects that are canonical subsets of $\Gamma^\infty$. 

 One possible source of such objects is described in forthcoming \cite{PFAKc}. There, the authors introduced the notion of $Z$-hod pairs and developed their basic theory. We should also note that even in this paper, to prove \rthm{thm:main_theorem} we build objects that resemble objects that are of third order. We build our third order objects more or less according to the current conventions following \cite{hod_mice_LSA}. What we meant above is that we believe that to get to superstrongs entirely new kinds of canonical objects need to be constructed. The reader can read more about such speculations in \cite{CCMDM}.
 
The abstract claimed that \rthm{thm:main_theorem} is the ultimate equiconsistency proved via $\sf{CMI}$. This does not mean that there are no other equiconsistencies in the region of $\sf{LSA}$. All it means is that to go beyond, one has to start thinking of $\sf{CMI}$ as a method of building third order objects. 

The authors view  \rthm{thm:main_theorem} as a natural accumulation point in the development of their understanding of $\sf{CMI}$ and the way it is used to translate set theoretic strength between natural set theoretic frameworks, namely between forcing axioms, large cardinals, determinacy and other frameworks. It has been proven by arriving at it via a 15 year long process of trying to understand $\sf{CMI}$. Because of this, we feel that it is a theorem proven by the entire community rather than by the authors. We especially thank Hugh Woodin and John Steel for their influential ideas throughout the first 25 years of the $\sf{Core\ Model\ Induction}$.

The history of \rthm{thm:main_theorem} is as follows. The first author, in \cite{sargsyan2013covering}, stated a conjecture that in his view captured the ideas of the first 15 years of the 2000s, namely that $\sf{\phi-Minimality}$ and $\sf{\neg\ UB-Covering}$ cannot co-exist in the short extender region. Unfortunately, very soon after finishing that paper he realized that the covering conjecture of that paper has to fail in the region of $\sf{LSA}$\footnote{The exact theorem was that if $\P$ is an lsa type hod premouse, $\d$ is the largest Woodin cardinal of $\P$, $\k<\d$ is the least $<\d$-strong cardinal that reflects the set of $<\d$-strong cardinals and $\mu$ is a $<\d$-strong cardinal larger than $\k$ then in $\P$, $\sf{UB-Covering}$ must fail at $\mu$. This theorem was presented at the Fourth European Set Theory Conference in Mon St Benet in 2013.}. However, no easily quotable theorem was proven by him. It was not until Fall of 2018 when the second author was visiting the first author, that they realized that \rthm{thm:main_theorem} says exactly that $\sf{\phi-Minimality}$ and $\sf{\neg\ UB-Covering}$ can coexist in the short extender region.  


\textbf{Acknowledgement.} We would like to thank Hugh Woodin for many comments made about the earlier drafts of this paper, and especially for pointing out an easier argument demonstrating that $\sf{LSA-over-uB}$ is not equivalent to $\sf{Sealing}$ and for bringing the $\sf{Tower\ Sealing}$ to our attention. We also thank Ralf Schindler for useful suggestion regarding our use of second order set theory. The authors are grateful for the referees' enormous efforts to make the manuscript better. 
The authors would like to thank the NSF for its generous support. The first author was supported by NSF Career Award DMS-1352034 and by the National Science Center, Poland under the Weave-UNISONO
call in the Weave programme, registration number UMO-2021/03/Y/ST1/00281. The second author is supported by NSF Grant No DMS-1565808 and DMS-1849295.

\section{An overview of the fine structure of the minimal $\sf{LSA}$-hod mouse and excellent hybrid mice}\label{sec:overview}

As was mentioned above, the proof of \rthm{thm:main_theorem} is an accumulation of many ideas developed in the last 20 years. We will try to develop enough of the required background in general terms so that a reader  familiar with the terminology of descriptive inner model theory can follow the arguments. The main technical machinery used in the proof is developed more carefully in \cite{hod_mice_LSA}. In the next few sections we will write an introduction to this technical machinery intended for set theorists who are familiar with \cite{BSL}. 

We say that $M$ is a minimal model of $\sf{LSA}$ if
 \begin{enumerate}
 \item $M\models \sf{LSA}$, 
 \item $M=L(A, \mathbb{R})$ for some $A\subseteq \mathbb{R}$, and
 \item  for any $B\in \powerset(\mathbb{R})\cap M$ such that $w(B)<w(A)$, $L(B,\mathbb{R})\models \neg \sf{LSA}$. 
 \end{enumerate}
 It makes sense to talk about ``the" minimal model of $\sf{LSA}$. When we say $M$ is \textit{the minimal} model of $\sf{LSA}$ we mean that $M$ is a minimal model of $\sf{LSA}$ and $Ord, \bR\subseteq M$. Clearly from the prospective of a minimal model of $\sf{LSA}$, the universe is the minimal model of $\sf{LSA}$. The proof of \cite[Theorem 10.3.1]{hod_mice_LSA} implies that there is a unique minimal model of $\sf{LSA}$ such that $Ord, \bR\subseteq M$\footnote{This proof of \cite[Theorem 10.3.1]{hod_mice_LSA} shows that the common part of a divergent models of $AD$ contains a minimal model of $\sf{LSA}$.}. This unique minimal model of $\sf{LSA}$ is \textit{the} minimal model of $\sf{LSA}$.
  
One of the main contributions of \cite{hod_mice_LSA} is the detailed description of $V_\Theta^\H$ assuming that the universe is the minimal model of $\sf{LSA}$. The early chapters of \cite{hod_mice_LSA} deal with what is commonly referred to as the $\H$ analysis. These early chapters introduce the notion of a \textit{short-tree-strategy} mouse, which is the most important technical notion studied by \cite{hod_mice_LSA}. To motivate the need for this concept, we first recall some of the other aspects of the $\H$ analysis.

Recall the $\sf{Solovay\ Sequence}$ (for example, see \cite[Definition 0.9]{hod_mice} or \cite[Definition 9.23]{Woodin}). Recall that $\Theta$ is the least ordinal that is not a surjective image of the reals. The $\sf{Solovay\ Sequence}$ is a way of measuring the complexity of the surjections that can be used to map the reals onto the ordinals below $\Theta$. Assuming $\sf{AD}$, let $(\theta_\a: \a\leq \Omega)$ be a closed in $\Theta$ sequence of ordinals such that
\begin{enumerate}
\item $\theta_0$ is the least ordinal $\eta$ such that $\mathbb{R}$ cannot be mapped surjectively onto $\eta$ via an ordinal definable function,
\item for $\a+1\leq \Omega$, fixing a set of reals $A$ such that $A$ has Wadge rank $\theta_\a$, $\theta_{\a+1}$ is the least ordinal $\eta$ such that $\mathbb{R}$  cannot be mapped surjectively onto $\eta$ via a function that is ordinal definable from $A$,
\item for limit ordinal $\l\leq \Omega$, $\theta_\l=\sup_{\a<\l}\theta_\a$, and
\item $\Omega$ is least such that $\theta_{\Omega}=\Theta$.
\end{enumerate}

It follows from the definition of $\sf{LSA}$ (\rdef{dfn:lsa})  that if $\k$ is the largest Suslin cardinal then it is a member of the $\sf{Solovay\ Sequence}$. It is not hard to show that $\sf{LSA}$ is a much stronger axiom than $\textsf{AD}_{\mathbb{R}}+``\Theta$ is regular". Under $\sf{LSA}$, letting $\k$ be the largest Suslin cardinal, there is an $\omega$-club $C\subseteq \k$ such that for every $\l\in C$, $L(\Gamma_\l, \mathbb{R})\models ``\textsf{AD}_\mathbb{R} + \l=\Theta+\Theta$ is regular", where $\Gamma_\l=\{A\subseteq \mathbb{R}: w(A)<\l\}$.\footnote{This theorem is probably due to Woodin. The outline of the proof is as follows. By an unpublished theorem of Woodin (but see \cite[Theorem 1.9]{ADRUB}), $\kappa$ is a measurable cardinal, as it is a regular cardinal. It follows that there is an $\omega$-club $C$ consisting of members of the Solovay sequence such that for all $\l\in C$, $\H\models ``\l$ is regular". Hence, $L(\Gamma_\l, \mathbb{R})\models ``\textsf{AD}_\mathbb{R} + \l=\Theta+\Theta$ is regular". For the proof of the last inference see \cite[Theorem 2.3]{SquarePaper}.}

 Assume now that $V$ is the minimal model of $\sf{LSA}$. It follows from the work done in \cite{hod_mice_LSA} that for every $\kappa$ that is a member of the $\sf{Solovay\ Sequence}$ but is not the largest Suslin cardinal there is a hod pair $(\P, \Sigma)$ such that
 \begin{enumerate}
 \item the Wadge rank of $\Sigma$ (or rather the set of reals coding $\Sigma$) is $\geq \k$ and
 \item for some $\eta\in \P$, letting $\M_\infty(\P, \Sigma)$ be the direct limit of all countable $\Sigma$-iterates $\Q$ of $\P$ such that the iteration embedding $\pi^\Sigma_{\P, \Q}$ is defined and letting $\pi^\Sigma_{\P, \infty}:\P\rightarrow \M_\infty(\P, \Sigma)$ be the iteration map, then $V_\k^\H$ is the universe of $\M_\infty(\P, \Sigma)|\pi^\Sigma_{\P, \infty}(\eta)$\footnote{Thus, $\pi^\Sigma_{\P, \infty}(\eta)=\k$.}.
 \end{enumerate}
 A technical reformulation of the above fact appears as \cite[Theorem 7.2.2]{hod_mice_LSA}.
 
 The situation, however, is drastically different for the largest Suslin cardinal. Let $\k$ be the largest Suslin cardinal. The inner model theoretic object that has Wadge rank $ \k$ cannot be an iteration strategy. This is because if $\Sigma$ is an iteration strategy with nice properties like \textit{hull condensation}\footnote{$\Sigma$ must also satisfy some form of generic interpretability, i.e., there must be a way to interpret $\Sigma$ on the the generic extensions of $\M_1^{\#, \Sigma}$.} then assuming $\sf{AD}$ holds in $L(\Sigma, \bR)$, $L(\Sigma, \bR)\models ``\M_1^{\#, \Sigma}$ exists and is $\omega_1$-iterable"\footnote{This can be proved by a $\Sigma^2_1$-reflection argument.}. This then easily implies that $\Sigma$ is both Suslin and co-Suslin. It then follows that no nice iteration strategy can have Wadge rank $\geq \k$, as any such strategy is both Suslin and co-Suslin\footnote{It follows from the theory of Suslin cardinals under $AD$ that $\k$ cannot be the largest Suslin cardinal, see \cite[Chapter 3]{Jackson}.}. 

The inner model theoretic object that has Wadge rank $ \k$ is a \textit{short tree strategy}, which is a partial iteration strategy. 
 Suppose $\P$ is any iterable structure and $\Sigma$ is its iteration strategy. Suppose $\d$ is a Woodin cardinal of $\P$. Given $\T\in dom(\Sigma)$ that is based on $\P|\d$, we say that $\T$ is $\Sigma$-\textit{short} if letting $\Sigma(\T)=b$, either the iteration map $\pi^\T_b$ is undefined or $\pi^\T_b(\d)>\d(\T)$. If $\T$ is not $\Sigma$-short then we say that it is $\Sigma$-\textit{maximal}. We then set $\Sigma^{stc}$ be the fragment of $\Sigma$ that acts on short trees. 
 
  Following \cite[Definition 3.1.4]{hod_mice_LSA} we make the following definition.
 \begin{definition}\label{mtsharp} Suppose $\T$ is a normal iteration tree of limit length. We then let 
 \begin{center}$\m(\T)=\cup_{\a < \lh(\T)}\M_\a^\T|\lh(E_\a^\T)$ and $\m^+(\T)=(\m(\T))^\#$.\end{center}
 \end{definition}
 
In the language of the above definition, the convention used in \cite{hod_mice_LSA} is the following: $\Sigma^{stc}(\T)=b$ if and only if 
 \begin{enumerate}
 \item $\T$ is $\Sigma$-short and $\Sigma(\T)=b$, or
 \item $\T$ is $\Sigma$-maximal and $b=\m^+(\T)$.
 \end{enumerate}
 Thus, $\Sigma^{stc}$ tells us the branch of a $\Sigma$-short tree or the last model of a $\Sigma$-maximal tree.

The reader can perhaps imagine many ways of defining the notion of \textit{short tree strategy} without a reference to an actual strategy. The convention that we adopt in this paper is the following. If $\Lambda$ is a short tree strategy for $\P$ then we will require that
 \begin{enumerate}
 \item for some $\P$-cardinal $\d$, $\P=(\P|\d)^\#$ and $\P\models ``\d$ is a Woodin cardinal",
 \item if $\d$ is as above and $\nu$ is the least $<\d$-strong cardinal of $\P$ then $\P\models ``\nu$ is a limit of Woodin cardinals",
 \item given an iteration tree $\T\in dom(\Lambda)$, $\Lambda(\T)$ is either a cofinal well-founded branch of $\T$ or is equal to $\m^+(\T)$,
 \item for all iteration trees $\T\in dom(\Lambda)$, if $\Lambda(\T)$ is a branch $b$ then $\pi^\T_b(\d)>\d(\T)$,
 \item for all iteration trees $\T\in dom(\Lambda)$, if $\Lambda(\T)$ is a model then $\m^+(\T)\models ``\d(\T)$ is a Woodin cardinal".
 \end{enumerate}
 If a hod mouse $\P$ has properties 1 and 2 above then we say that $\P$ is of $\#$-\textit{lsa type}.  \cite[Definition 2.7.3]{hod_mice_LSA} introduces other types of LSA hod premice.
 
 The set of reals that has Wadge rank $\k$ is some short tree strategy $\Lambda$. The hod mouse $\P$ that $\Lambda$ iterates has a unique Woodin cardinal $\d$ such that
if $\nu<\d$ is the least cardinal that is $<\d$-strong in $\P$, then $\P\models ``\nu$ is a limit of Woodin cardinals". The aforementioned Woodin cardinal $\d$ is also the largest Woodin cardinal of $\P$. This fact is proven in \cite{hod_mice_LSA} (for example, see \cite[Theorem 7.2.2]{hod_mice_LSA} and \cite[Chapter 8]{hod_mice_LSA}). There is yet another way that the $\sf{LSA}$ stages of the $\sf{Solovay\ Sequence}$ are different from other points. 
 
 We continue assuming that $V$ is the minimal model of $\sf{LSA}$. If $\Sigma$ is a strategy of a hod mouse with nice properties then ordinal definability with respect to $\Sigma$ is captured by $\Sigma$-mice. More precisely, \cite[Theorem 10.2.1]{hod_mice_LSA} implies that if $x$ and $y$ are reals then $x$ is ordinal definable from $y$ using $\Sigma$ as a parameter if and only if there is a $\Sigma$-mouse $\M$ over $y$\footnote{The difference between a mouse and a mouse over $y$ is the same as the difference between $L$ and $L[x]$.} such that $x\in \M$. 
 
 \cite[Theorem 10.2.1]{hod_mice_LSA} also implies that the same conclusion is true for short tree strategies. Namely, if $\Lambda$ is a short tree strategy then for $x$ and $y$ reals, $x$ is ordinal definable from $y$ using $\Lambda$ as a parameter if and only if there is a $\Lambda$-mouse $\M$ over $y$ such that $x\in \M$. Theorems of this sort are known as $\sf{Mouse\ Capturing}$ theorems. Such theorems are very important when analyzing models of determinacy using inner model theoretic tools. 
 
For a strategy $\Sigma$ the concept of a $\Sigma$-mouse has appeared in many places. The reader can consult \cite[Definition 1.20]{hod_mice} but the notion probably was first mentioned in \cite{PFA} and was finally fully developed in \cite{operators}. 

A $\Sigma$-mouse $\M$, besides having an extender sequence also has a predicate that indexes the strategy. The idea, which is due to Woodin, is that the strategy predicate should index the branch of the least tree that has not yet been indexed. 

Unfortunately this idea doesn't quite work for $\Lambda$-mice where $\Lambda$ is a short tree strategy. In the next subsection, we will explain the solution presented in \cite{hod_mice_LSA}.

\subsection{Short tree strategy mice}\label{subsec: short tree strategy mice}
 
  We are assuming that $V$ is the minimal model of $\sf{LSA}$. Suppose $\Lambda$ is a short tree strategy for a hod mouse $\P$. We let $\d$ be the largest Woodin cardinal of $\P$. Thus, $\P=(\P|\d)^\#$. In this subsection, we would like to convince the reader that the concept of $\Lambda$-mouse, while much more involved, behaves very similarly to the concept of a $\Sigma$-mouse where $\Sigma$ is an iteration strategy. 
  
In general, when introducing any notion of a mouse one has to keep in mind the procedures that allow us to build such mice. Formally speaking, many notions of $\Lambda$-mice might make perfect sense, but when we factor into it the constructions that are supposed to produce such mice we run into a key issue.

 In any construction that produces some sort of mouse (e.g. $K^c$-constructions, fully backgrounded constructions, etc) there are stages where one has to consider certain kinds of Skolem hulls, or as inner model theorists call them fine structural \textit{cores}. The reader can view these cores as some carefuly defined Skolem hulls. To illustrate the aformentioned problem, imagine we do have some notion of $\Lambda$-mice and let us try to run a construction that will produce such mice. Suppose  $\T$ is a tree according to $\Lambda$ that appears in this construction. Having a notion of a $\Lambda$-mouse means that we have a prescription for deciding whether $\Lambda(\T)$ should be indexed in the strategy predicate or not. 
 
Suppose $\T$ is a $\Lambda$-maximal tree. It is hard to see exactly what one can index so that the strategy predicate remembers that $\T$ is maximal. And this ``remembering" is the issue. Imagine that at a later stage we have a Skolem hull $\pi: \M\rightarrow \N$ of our current stage such that $\T\in rng(\pi)$. It is possible that $\U=_{def}\pi^{-1}(\T)$ is $\Lambda$-short. If we have indexed $X$ in our strategy that proves $\Lambda$-maximality of $\T$ then $\pi^{-1}(X)$ now can no longer prove that $\U$ is $\Lambda$-maximal. Thus, the notion of $\Lambda$-mouse cannot be first order. 

The solution is simply not to index anything for $\Lambda$-maximal trees. This doesn't quite solve the problem as the above situation implies that nothing should be indexed for many $\Lambda$-short trees as well. To solve this problem, we will only index the branches of some $\Lambda$-short trees, those that we can locally prove are $\Lambda$-short. We explain this below in  more details.
  
Fix an lsa type hod premouse $\P$ and let $\Lambda$ be its short tree strategy. Let $\d$ be the largest Woodin cardinal of $\P$ and $\nu$ be the least $<\d$-strong of $\P$. To explain the exact prescription that we use to index $\Lambda$, we explain some properties of the models that have already been constructed according to this indexing scheme.  Suppose $\M$ is a $\Lambda$-premouse. 


Call $\T\in \M$  \textit{universally short} ($\sf{uvs}$) if $\T$ is obviously \textit{short} (see \cite[Definition 3.3.2]{hod_mice_LSA}).  For instance, it can be that the $\#$-operator provides a $\Q$-structure and determines a branch $c$ of $\T$ such that $\Q(c, \T)$\footnote{$\M^\T_c$ is a direct limit along the models of $c$. $\Q(c, \T)$ is the largest initial segment of $\M^\T_c$ such that $Q(c, \T)\models ``\d(\T)$ is a Woodin cardinal". It is only defined provided that $\d(\T)$ is not a Woodin cardinal for some function definable over $\M^\T_c$.} exists and $\Q(c, \T)\unlhd\m^+(\T)$. Another way that a tree can be obviously short is that there could be a model $\Q$ in $\T$ such that $\pi^\T_{\P, \Q}:\P\rightarrow \Q$ is defined and the portion of $\T$ that comes after $\Q$ is based on $\Q^b$. Here $\Q^b$ is defined as $\Q|(\kappa^+)^\Q$, where $\kappa$ is the supremum of the Woodin cardinals below the largest Woodin of $\Q$. The reader should keep in mind that there is a formula $\zeta$ in the language of $\Lambda$-premice such that for any $\Lambda$-premouse $\M$ and for any iteration tree $\T\in \M$, $\T$ is $\sf{uvs}$ if and only if $\M\models \zeta[\T]$.

Unfortunately there can be trees that are not universally short ($\sf{nuvs}$). Suppose then $\T$ is $\sf{nuvs}$. In this case whether we index $\Lambda(\T)$ or not depends on whether we can find a $\Q$-structure that can be authenticated  to be the correct one. There can be many ways to certify a $\Q$-structure, and \cite{hod_mice_LSA}
provides one such method. An interested reader can consult \cite[Chapter 3.7]{hod_mice_LSA}. Notice that because $\P$ has only one Woodin cardinal, not being able to find a $\Q$-structure is equivalent to the tree being maximal. Thus, in a nutshell the solution proposed by \cite{hod_mice_LSA} is that we index only branches that are given by \textit{internally} authenticated  $\Q$-structures. 

Suppose now that we have the above Skolem hull situation, namely that we have $\pi:\M\rightarrow \N$ and $\T$ in $\N$ that is $\Lambda$-maximal but $\pi^{-1}(\T)$ is short. There is no more indexing problem. The reason is that in order to index $\Lambda(\pi^{-1}(\T))$ in $\M$ we need to find an authenticated  $\Q$-structure for $\pi^{-1}(\T)$. The authentication process is first order, and so if $\N$ does not have such an authenticated  $\Q$-structure for $\T$ then $\M$ cannot have such an authenticated  $\Q$-structure for $\pi^{-1}(\T)$. 

The reader of this paper does not need to know the exact way the authentication procedure works. However, the reader should keep in mind that the authentication procedure is internal to the mouse. More precisely, the following holds:\\

\textbf{Internal Definability of Authentication:} there is a formula $\phi$ in the appropriate language such that whenever $(\P, \Lambda)$ is as above and $\M$ is a $\Lambda$-mouse over some set $X$ such that $\P\in X$, for any iteration tree $\T\in \M$, $\M\models \phi[\T]$ if and only if $\T\in dom(\Lambda)$, $\T$ is short and $\Lambda(\T)\in \M$. \\\\
We again note that the Internal Definability of Authentication (IDA) is only shown to be true for the minimal model of $\sf{LSA}$. In general, IDA cannot be true as there can be short trees without $\Q$-structures. The authors have recently discovered another short tree indexing scheme that can work in all cases, but has some weaknesses compared to the one introduced in \cite{hod_mice_LSA}. 

Using the notation in \cite{hod_mice_LSA}, recall that $\P^b$ is the ``bottom part" of $\P$, i.e $\P^b = \P|(\nu^+)^{\P}$, where $\nu$ is the supremum of the Woodin cardinals below the top Woodin of $\P$.

We now describe another key feature of the indexing scheme of \cite{hod_mice_LSA} that is of importance here. 
We say $\Sigma$ is a \textit{low level component} of $\Lambda$ if there is a tree $\T$ on $\P$ according to $\Lambda$ such that $\pi^{\T, b}$ exists\footnote{$\pi^{\T,b}$ is the restriction of the iteration embedding to $\P^b$. See \cite{hod_mice_LSA}, just after Definition 2.7.21, for a more detailed definition. Note that in some cases, $\pi^{\T,b}$ may exist but $\pi^\T$ may not.} ($\T$ may be $\emptyset$) 
and for some $\R\unlhd \pi^{\T, b}(\P^b)$, $\Sigma=\Lambda_\R$. Let $LLC(\Lambda)$ be the set of $\Sigma$ that are a low level components of $\Lambda$. What is shown in \cite{hod_mice_LSA} is that $\Lambda$ is determined by $LLC(\Lambda)$ in a strong sense. 

Given a transitive model $M$ of a fragment of $\sf{ZFC}$ such that $\P\in M$ we say $M$ is closed under $LLC(\Lambda)$ if whenever $\T\in M$ is a tree according to $\Lambda$ such that $\pi^{\T, b}$ exists, $\Lambda_{\pi^{\T, b}(\P^b)}$ has a universally Baire representation over $M$. More precisely, whenever $g\subseteq Coll(\omega, \pi^{\T, b}(\P^b))$ is $M$-generic, for every $M$-cardinal $\l$ there are trees $T, S\in M[g]$ on $\l$ such that $M[g]\models ``(T, S)$ are $<\l$-complementing" and for all $<\l$-generics $h$, $(p[T])^{M[g*h]}=Code(\Lambda_{\pi^{\T, b}(\P^b)})\cap M[g*h]$. Here $Code(\Phi)$ is the set of reals coding $\Phi$ (with respect to a fixed coding of elements of $HC$ by reals).

It is shown in \cite{hod_mice_LSA} that if, assuming $\sf{AD}^+$, $(M, \Sigma)$ is such that
\begin{enumerate}
\item $M$ is a countable model of a fragment of $\sf{ZFC}$,
\item $M$ has a class of Woodin cardinals,
\item $\Sigma$ is an $\omega_1$-iteration strategy for $M$ and
\item whenever $i:M\rightarrow N$ is an iteration via $\Sigma$, $N$ is closed under $LLC(\Lambda)$,
\end{enumerate}
then there is a formula $\psi$ such that whenever $g$ is $M$-generic, for any $\T\in M[g]$,
\begin{center}
$\T$ is according to $\Lambda$ if and only if $M[g]\models \psi[\T]$.  \ \ \ \ \ \ \ \ \ \  ($\star$)
\end{center}  
The interested reader can consult Chapters 5, 6 and 8 of \cite{hod_mice_LSA}.

The reason we explained the above is to give the reader some confidence that defining a short tree strategy $\Lambda$ for a hod premose $\P$ is equivalent to describing the set $LLC(\Lambda)$. This fact is the reason that the indexing schema of \cite{hod_mice_LSA} works in the following sense.

Being able to define short-tree-strategy mice is one thing, proving that they are useful is another. Usually what needs to be shown are the following two key statements. We let $\phi_{sts}$ be the formula that is mentioned in the Internal Definability of Authentication.\\

\textbf{The Eventual Authentication.} Suppose $(\P, \Lambda)$ is as above and $\M$ is a sound $\Lambda$-mouse over some set $X$ such that $\P\in X$ and $\M$ projects to $X$. Suppose $\T\in \M$ is according to $\Lambda$ and is $\Lambda$-short. Suppose further that $\M\models \neg \phi_{sts}[\T]$. Then there is a sound $\Lambda$-mouse $\N$ over $X$ such that $\M\unlhd \N$ and $\N\models \phi_{sts}[\T]$.\footnote{One can then prove that there is such an $\N$ that projects to $X$.} \\

\textbf{Mouse Capturing for $\Lambda$:}  Suppose $(\P, \Lambda)$ is as above. Then for any $x\in \mathbb{R}$ that codes $\P$ and any $y\in \mathbb{R}$, $y$ is ordinal definable from $x$ and $\Lambda$ if and only if there is a $\Lambda$-mouse $\M$ over $x$ such that $y\in \M$. \\

Both The Eventual Authentication and Mouse Capturing for $\Lambda$ are proven in \cite{hod_mice_LSA} (see \cite[Chapter 8, Lemma 8.1.3, Lemma 8.1.5]{hod_mice_LSA} and \cite[Theorem 10.2.1]{hod_mice_LSA}). 

The next subsection discusses the $\Q$-structure authentication process mentioned above. 

\subsection{The authentication method}\label{authentication method subsec}

Suppose $\P$ is a $\#$-lsa type hod premouse. Recall from the previous subsections that this means that $\P$ has a largest Woodin cardinal $\d$ such that $\P=(\P|\d)^\#$ and the least $<\d$-strong cardinal of $\P$ is a limit of Woodin cardinals. We let $\d^\P$ be the largest Woodin cardinal of $\P$ and $\k^\P$ be the least $<\d^\P$-strong cardinal of $\P$. We shall also require that $\P$ is \textit{tame}, meaning that for any $\nu<\d^\P$, if $(\P|\nu)^\#$ is of lsa type and $\M\insegeq \P$ is the largest such that $\M\models ``\nu$ is a Woodin cardinal" then $\nu$ is not overlapped in $\M$\footnote{This means that if $E\in \vec{E}^\M$ then $\nu\not \in(\cp(E), index(E))$.}.

Our goal here is to explain the $\Q$-structure authentication procedure employed by \cite{hod_mice_LSA}. Recall our discussion of $\sf{uvs}$ and $\sf{nuvs}$ trees.  The $\Q$-structure authentication procedure applies to only $\sf{nuvs}$ trees, trees that are not obviously short.

 \cite[Chapters 3.6-3.9]{hod_mice_LSA} develop the aforementioned authentication procedure. \cite[Definition 3.8.9, 3.8.16, 3.8.17]{hod_mice_LSA} introduce the sts indexing scheme. For illustrative purposes, it is better to think of the indexing scheme introduced there as a hierarchy of indexing schemes indexed by ordinals. Naturally, this hierarchy is defined by induction. For illustrative purposes we call $\gg$th level of the  hierarchy $sts_\gg$. Thus, $sts_\gg(\P)$ is the set of all sts premice that are based on $\P$ (i.e., their short tree strategy predicate describes a short tree strategy for $\P$) and have rank $\leq \gg$.

To begin the induction, we let $sts_0(\P)$ be the set of all sts premice that do not index a branch for any $\sf{nuvs}$ tree. More precisely, if $\M\in sts_0(\P)$ and $\T\in dom(S^\M)$ then if $S^\M(\T)$ is defined then $\T$ is $\sf{uvs}$. 

Below and elsewhere, $S^\M$ is the strategy predicate of $\M$. Given $sts_{\a}(\P)$ we let $sts_{\a+1}(\P)$ be the set of all sts premice that index branches of those $\sf{nuvs}$ trees that have a $\Q$-structure in $sts_\a(\P)$. More precisely, suppose $\M\in sts_{\a+1}(\P)$ and $\T\in dom(S^\M)$ and $S^\M(\T)$ is defined. Then either
\begin{enumerate}
\item $\T$ is $\sf{uvs}$ or
\item $\T$ is $\sf{nuvs}$ and there is $\Q\in\M$ such that $\M\models ``\Q\in sts_\a(\P)"$, $\m^+(\T)\insegeq\Q$, $\Q\models ``\d(\T)$ is a Woodin cardinal" but $\d(\T)$ is not a Woodin cardinal with respect to some function definable over $\Q$\footnote{This can be written as $\mathcal{J}_1(\Q)\models ``\d(\T)$ is not a Woodin cardinal".}  and there is a cofinal branch $b$ of $\T$ such that $\Q\insegeq \M^\T_b$.
\end{enumerate}
When $\Q$ exhibits the properties listed in clause 2 we say that $\Q$ is a $\Q$-structure for $\T$. It follows from the zipper argument of \cite[Theorem 2.2]{IT} that for each $\Q$-structure $\Q$ there is at most one branch $b$ with properties described in clause 2 above. However, there is nothing that we have said so far that guarantees the uniqueness of the $\Q$-structure itself. The uniqueness is usually a consequence of iterability and comparison (see \cite[Theorem 3.11]{steel2010outline})\footnote{In general, the theory of $\Q$-structures doesn't have much to do with sts mice. It will help if the reader develops some understanding of \cite[Chapter 6.2 and Definition 6.11]{steel2010outline}.}. Thus, to make the definition of $sts_{\a+1}$ complete, we need to impose an iterability condition on $\Q$.

The exact iterability condition that one needs is stated as clause 5 of \cite[Definition 3.8.9]{hod_mice_LSA}. This clause may seem technical, but there are good reasons for it. For the purposes of identifying a unique branch $b$ saying that $\Q$ in clause 2 is sufficiently iterable in $\M$ would have sufficed. However, recall the statement of the Internal Definability of Authentication. The problem is that when we require that an $\M$ as above is a $\Lambda$-premouse we in addition must say that the branch $b$ that the $\Q$-structure $\Q$ defines is the exact same branch that $\Lambda$ picks. To guarantee this, we need to impose a condition on $\Q$ such that $\Q$ will be iterable not just in $\M$ but in $V$. The easiest way of doing this is to say that $\Q$ has an iteration strategy in some derived model as then, using genericity iterations (see \cite[Chapter 7.2]{steel2010outline}), we can extend such a strategy for $\Q$ to a strategy that acts on iterations in $V$. 

For limit $\a$, $sts_\a(\P)$ is essentially $\bigcup_{\b<\a}sts_\b(\P)$. What has been left unexplained is the kind of strategy that the $\Q$-structure $\Q$ must have in some derived model. Let $\Sigma$ be this strategy. If $\M\in sts_\a(\P)$ is a $\Lambda$-mouse then $\Q$ must be a $\Lambda_{\m^+(\T)}$-mouse over $\m^+(\T)$. Thus, our next challenge is to find a first order way of guaranteeing that $\Sigma$-iterates of $\Q$ are  $\Lambda_{\m^+(\T)}$-mice, even those iterates that we will obtain after blowing up $\Sigma$ via genericity iterations. 

The solution that is employed in \cite{hod_mice_LSA} is that if $\R$ is a $\Sigma$-iterate of $\Q$ and $\U\in dom(S^\R)$ then $\U$ itself is authenticated  by the extenders of $\M$. Below we refer to this certification as \textit{tree certification}. This is again a rather technical notion, but the following essentially illustrates the situation.

 Let us suppose $\R=\Q$ and $\U\in dom(S^\Q)$. The indexing scheme of \cite{hod_mice_LSA} does not index all trees in $\P$. In other words, $S^\M$ is never total. $dom(S^\M)$ consists of trees that are built via a comparison procedure that iterates $\P$ to a background construction of $\M$. Set $\N=\m^+(\U)$. One requirement is that $\N$ also iterates to one such background construction to which $\P$ also iterates. Let $\S$ be this common background construction and suppose $\a+1<lh(\U)$ is such that $\a$ is a limit ordinal. First assume $\U\rest \a$ is $\sf{uvs}$. What is shown in \cite{hod_mice_LSA} is that knowing the branch of $\P$-to-$\S$ tree there is a first order procedure that identifies the branch of $\U\rest \a$, and that procedure is the tree certification procedure applied to $\U\rest \a$.

 Suppose next that $\U\rest \a$ is $\sf{nuvs}$. Then because $\a+1<lh(\U)$, $\U\rest \a$ must be short and the branch chosen for it in $\Q$ must have a $\Q$-structure $\Q_1$ which is itself an sts mouse. We have that $\Q_1\in \Q$ and $\Q_1$ must have the same certification in $\Q$ that $\Q$ has in $\M$. Again, the $\sf{nuvs}$ trees in $\Q_1$ have a tree certification in $\Q$ according to the above procedure. The $\sf{uvs}$ ones produce another $\Q_2\in \Q_1$. Because we cannot have an infinite descent, the definition of tree certification is meaningful. 
 
 \begin{remark}\label{important remark}  It is sometimes convenient to think of a short tree strategy as one having two components, the \textit{branch component} and the \textit{model component}. Given a short tree strategy $\Lambda$, we let $b(\Lambda)$ be the set of those trees $\T\in dom(\Lambda)$ such that $\Lambda(\T)$ is a branch of $\T$, and we let $m(\Lambda)$ be the set of those trees $\T\in dom(\Lambda)$ such that $\Lambda(\T)$ is a model.

The convention adopted in this paper is that if $\T\in m(\Lambda)$ then $\Lambda(\T)=\m^+(\T)$\footnote{It is not up to us to decide whether $\Lambda(\T)\in m(\Lambda)$ or $\Lambda(\T)\in b(\Lambda)$. The short-tree strategy itself decides this.}. Thus, if $\M$ is an sts premouse then $S^\M$ is a short tree strategy in the above sense, i.e., for $\T\not \in b(S^\M)$, $S^\M(\T)$ is simply left undefined. 
 \end{remark}

This ends our discussion of sts premice. Of course, a lot has been left out and the mathematical details are unfortunately excruciating, but we hope that the reader has gained some level of intuition to proceed with the paper.  

In the next subsection we will deal with one of the most important aspects of hod mice, namely the generic interpretability of iteration strategies. 

\subsection{Generic interpretability}

There are several situations when one has to be careful when discussing sts premice and $\Lambda$-premice in general. First, for an iteration strategy $\Sigma$, $\M_1^{\#, \Sigma}$ makes complete sense. It is the minimal active $\Sigma$-mouse with a Woodin cardinal. For short tree strategy $\Lambda$ the situation is somewhat different. The expression``$\M_1^{\#, \Lambda}$ is the minimal active $\Lambda$-mouse with a Woodin cardinal" doesn't say much as we do not say how closed $\M_1^{\#, \Lambda}$ must be. One must also add statements of the form ``in which all $\Lambda$-short trees are indexed". This is because it could be that $\Lambda$-premouse $\M$ is active and has a Woodin cardinal but there is a $\Lambda$-short tree $\T\in \M$ that has not yet been indexed in $\M$ (see The Eventual Certification above). In particular, without extra assumptions, it may be the case that given a $\Lambda$-sts mouse $\M\models \sf{ZFC}$, $\Lambda\rest \M$ is not definable over $\M$. Clearly such definability holds for many strategy mice.

The above issue becomes somewhat of a problem when dealing with  \textit{generic interpretability}, which is the statement that the internal strategy predicate can be uniquely extended onto generic extensions. For ordinary strategy mice, generic interpretability is in general easier to prove. For short tree strategy mice the situation is somewhat parallel to the above anomaly. Suppose $\M$ is a $\Lambda$-mouse where $\Lambda$ is a short tree strategy and suppose $g$ is $\M$-generic. In general, we cannot hope to prove that $\Lambda\rest \M[g]$ is definable over $\M[g]$. 
In this subsection, we introduce some properties of short tree strategies that allow us to prove generic interpretability, albeit in a  somewhat weaker sense.
  
The most important concept that is behind most arguments of \cite{hod_mice_LSA} is the concept of \textit{branch condensation} (see \cite[Chapter 4.9]{hod_mice_LSA}). It is very possible that the concept of \textit{full normalization} introduced in \cite{normalization_comparison} can be used instead of \textit{branch condensation} to obtain a greater generality. In fact, the authors have recently discovered a new notion of a short-tree-strategy mouse utilizing full normalization. 

Branch condensation implies generic interpretability. The following is our generic interpretability theorem which is essentially  \cite[Theorem 6.1.5]{hod_mice_LSA}. The aforementioned theorem is stated for strategies with branch condensation that are associated with a pointclass $\Gamma$. Here, we need strategies whose association with pointclasses is a consequence of some abstract properties that it has, not something explicitly assumed about them.  Such strategies can be obtained working inside a model of determinacy. The specific properties that we need are the following properties:
\begin{enumerate}
\item hull condensation,
\item strong branch condensation,
\item branch condensation for pull-backs.
\end{enumerate}
The meaning of clause 3 above is as follows. Suppose $(\P, \Lambda)$ is an sts hod pair.  $\Lambda$ has branch condensation for pullbacks if whenever $\xi\in \P$ is a limit of Woodin cardinals of $\P$ such that $\P\models ``\cf(\xi)=\omega"$ and  $\pi:\Q\rightarrow \P|\xi$ is elementary, the $\pi$-pullback of $\Lambda_{\P|\xi}$ has branch condensation. For more on branch condensation, the reader may consult \cite[Chapter 4.9]{hod_mice_LSA}.
\begin{definition}\label{splendid strategies} We say that a short tree strategy is splendid if it satisfies the above 3 properties.
\end{definition}
Clause 3 above implies that the pullback of splendid strategies are splendid\footnote{Simply because ``being a pullback" is a transitive property.}. It might help to consult \rrem{important remark} before reading the next theorem.

\begin{theorem}\label{generic interpretability} Suppose $\P$ is an lsa type hod premouse and $\Lambda$ is a splendid short tree strategy for $\P$. Suppose $\N$ is a $\Lambda$-premouse satisfying $\sf{ZFC}$ and that $\N$ has unboundedly many Woodin cardinals. Then for any $\N$-generic $g$, $\Psi=_{def}S^\N$ has a unique extension $\Psi^g\subseteq \N[g]$ that is definable from $\Psi$ over $\N[g]$ and $b(\Psi^g)\subseteq b(\Lambda)\rest \N[g]$. 
\end{theorem}

Our \rthm{generic interpretability} is weaker than  \cite[Theorem 6.1.4]{hod_mice_LSA}. The conclusion of the aforementioned theorem is that $\Psi^g=\Lambda\rest \N[g]$. However, in  \cite[Theorem 6.1.4]{hod_mice_LSA} $\N$ satisfies a strong iterability hypothesis. Without this iterability hypothesis, $b(\Psi^g)\subseteq b(\Lambda)\rest \N[g]$ is all the proof of \cite[Theorem 6.1.4]{hod_mice_LSA} gives.

In the next subsection we will introduce a type of short tree strategy mouse that we will use to establish \rthm{thm:main_theorem}.

\subsection{Excellent hod premice}

Our proof of \rthm{thm:main_theorem} is an example of how one can translate set theoretic strength from one set of principles to another set of principles by using inner model theoretic objects as intermediaries. Below we introduce the notion of an \textit{excellent hybrid premouse}. We will then use this notion to show that both $\sf{Sealing}$ and $\sf{LSA-over-uB}$ hold in a generic extension of an excellent hybrid premouse. Conversely, we will show that in any model of either $\sf{Sealing}$ or $\sf{LSA-over-uB}$ there is an excellent hybrid premouse. We start by introducing some terminology and then introduce the excellent hybrid premice.

\begin{remark} Below and elsewhere, when discussing iterability we usually mean with respect to the extender sequence of the structures in consideration. Sometimes our definitions will be stated with no reference to such an extender sequence, but these definitions will always be applied in contexts where there is  a distinguished extender sequence. 
\end{remark}

To state our generic interpretability results we need to introduce a form of self-iterability, namely \textit{window-based self-iterability}. We say that $[\nu, \d]$ is a \textit{window} if there are no Woodin cardinals in the interval $(\nu, \d)$. Given a window $w$, we let $\nu^w$ and $\d^w$ be such that $w=[\nu^w, \d^w]$. We say that a window $w$ is above $\k$ if $\nu^w\geq\k$. We say that a window $w$ is not \textit{overlapped} if there is no $\nu^w$-strong cardinal. We say $w$ is maximal if $\nu^w=\sup\{\a+1: \a<\nu^w$ is a Woodin cardinal$\}$ and $\d^w$ is a Woodin cardinal.\\\\
\textbf{Window-Based Self-Iterability.} Suppose $\k$ is a cardinal. We say $\sf{WBSI}$ holds at $\k$ if for any window $w$ that is above $\k$ and for any successor cardinal $\eta\in (\nu^w, \d^w)$, setting $Q= H_{\eta^+}$, $Q$ has an $Ord$-iteration strategy $\Sigma$ which acts on iterations that only use extenders with critical points $>\nu^w$. \\\\
One usually says that $Q$ is $Ord$-iterable above $\nu^w$ to mean exactly what is written above. 

\begin{definition}\label{theory t0} We let $T_0$ be the conjunction of the following statements.
\begin{enumerate}
\item $\sf{ZFC}$,
\item There are unboundedly many Woodin cardinals.
\item The class of measurable cardinals is stationary.
\item No measurable cardinal that is a limit of Woodin cardinals carries a normal ultrafilter concentrating on the set of measurable cardinals.
\end{enumerate}
\end{definition}

When we write $M\models T_0$ and $M$ has a distinguished extender sequence then we make the tacit assumption that the large cardinals and specific ultrafilters mentioned in \rdef{theory t0} are witnessed by extenders on the sequence of $M$.  

\begin{definition}\label{dfn:hod_pm}
Suppose $\P$ is hybrid premouse. We say that $\P$ is \textbf{almost excellent} if
\begin{enumerate}
\item $\P\models T_0$.
\item There is a Woodin cardinal $\d$ of $\P$ such that $\P\models ``\P_0=_{def}(\P|\d)^{\#}$ is a hod premouse of $\#$-lsa type", $\P$ is an sts premouse based on $\P_0$ and $\P\models ``S^\P$, which is a short tree strategy for $\P_0$,  is splendid".
\item Given any $\tau<\d^{\P_0}$ such that $(\P_0|\tau)^{\#}$ is of $\#$-lsa type, there is $\M\insegeq \P$ such that $\tau$ is a cutpoint of $\M$ and $\M\models ``\tau$ is  not a Woodin cardinal".

We say that $\P$ is \textbf{excellent} if in addition to the above clauses,  $\P$ satisfies:
\item Letting $\d$ be as in clause 2, $\P\models ``\sf{WBSI}$ holds at $\d$".
\end{enumerate}
If $\P$ is excellent then we let $\d^\P$ be the $\d$ of clause 2 above and $\P_0=((\P|\d^\P)^\#)^\P$.
\end{definition}

\begin{remark} In the previous subsection, we were mainly concerned with the structure of hod mice associated with the minimal model of $\sf{LSA}$. An excellent hybrid premouse is beyond the minimal model of $\sf{LSA}$. Indeed, the arguments used in the proofs of \cite[Lemma 8.1.10 and Theorem 8.2.6]{hod_mice_LSA} apply to show that if $\P$ is excellent and $\l>\d^\P$ is a limit of Woodin cardinals of $\P$ then the (new) derived model at $\l$ is a model of $\sf{LSA}$. It follows from a standard Skolem hull argument and the derived model theorem that there is $A\in \powerset(\bR^\P)\cap \P$ such that $L(A, \bR^\P)\models \sf{LSA}$. 

Nevertheless, everything that we have said in the previous subsection about short tree strategies and sts mice carries over to the level of excellent hybrid mice. The methods of \cite{hod_mice_LSA} work through the \textit{tame}\footnote{A non-tame hod premouse is one that has an extender overlapping a Woodin cardinal.} hod mice. The authors recently have discovered a new sts indexing scheme that works for arbitrary hod mice.  This work is not relevant to the current work as the indexing of \cite{hod_mice_LSA} just carries over verbatim. 
\end{remark}

For the rest of this paper we assume the following minimality hypothesis $\neg (\dagger)$, where \\\\
$(\dagger):$ \ \  In some generic extension there is a (possibly class-sized) excellent hybrid premouse.\\\\
We will periodically remind the reader of this. One consequence of this assumption is the following fact, which roughly says that all local non-Woodin cardinals of a hod premouse (or hybrid premouse) are witnessed by $\Q$-structures which are initial segments of the model and are tame. It also shows that if $\P$ is a hod mouse such that there is an lsa initial segment $\P_0$ of $\P$ and there is a Woodin cardinal $\delta > o(\P_0)$ inside $\P$, then we can construct an excellent hybrid premouse in $\P$ by essentially performing a fully backgrounded sts construction in $\P|\delta$ above $\P_0$ (with respect to the short-tree component of $\P_0$).

\begin{proposition}[$\neg (\dagger)$]\label{simple structure} Suppose $\P$ is a hod premouse. Let $\kappa$ be a measurable limit of Woodin cardinals of $\P$ and let $\xi\leq o^\P(\kappa)$. Suppose $(\P|\xi)^\#\models ``\xi$ is a Woodin cardinal" but either $\xi$ is not the largest Woodin cardinal of $\P$ or $\xi<o^\P(\kappa)$. Then there is $\M\insegeq \P$ such that $\xi$ is a cutpoint in $\M$, $\rho(\M)\leq \xi$ and $\M\models ``\xi$ is not a Woodin cardinal".  
\end{proposition} 
\begin{proof} Towards a contradiction assume that there is no such $\M$. Suppose first $\xi$ is a Woodin cardinal of $\P$. It must then be a cutpoint cardinal as otherwise we easily get an excellent hybrid premouse by performing a fully backgrounded construction inside $\P|\kappa$ with respect to the short-tree component of $\P_0$, where $\P_0$ is an lsa hod initial segment of $\P|\kappa$. The existence of $\P_0$ follows from the fact that $(\P|\xi)^\sharp$ is an lsa initial segment of $\P$ and $\xi < o^\P(\kappa)$. 

It then follows that there is a Woodin cardinal $\zeta$ of $\P$ above $\xi$. Now we can use \cite[Lemma 8.1.4]{hod_mice_LSA} to build an excellent hybrid premouse via a backgrounded construction of $\P|\zeta$ as above (with respect to the short-tree component of $(\P|\xi)^\sharp$).

Suppose next that $\xi$ is not a Woodin cardinal. Because no $\M$ as above exists, it follows that $\xi<o^\P(\kappa)$. We can now repeat the above steps in $Ult(\P, E)$ where $E$ is the least extender overlapping $\xi$. 
\end{proof}

There are a few important facts that we will need about excellent hybrid premice that one can prove by using more or less standard ideas, and that in one form or another have appeared in \cite{hod_mice_LSA}. We will use the next subsection recording some of these facts.

\subsection{More on self-iterability}

Here we prove that window based strategy acts on the entire model. The main theorem that we would like to prove is the following.

\begin{theorem}\label{thm:full iterability in excellent hybrid mice generic} Suppose $\P$ is an excellent hybrid premouse, $w$ is a maximal window of $\P$ above $\d^\P$ and $\eta\in [\nu^w, \d^w)$ is a regular cardinal. Let $\Sigma$ be the $Ord$-strategy of $\Q=\P|\eta$ that acts on iterations that are above $\nu^w$. Let $g$ be $\P$-generic for some poset $\mathbb{P}\in \P$. Then $\Sigma$ has a unique extension $\Sigma^g$ definable over $\P[g]$ such that in $\P[g]$, $\Sigma^g$ is an $Ord$-iteration strategy for $\P$ that acts on iterations that are based on $\Q$ and are above $\nu^w$.
\end{theorem}
The proof will be presented as a sequence of lemmas. First we make a few  observations. Suppose $\P$ is excellent and for some $\P$-cardinal $\xi>\d^\P$ is a limit of Woodin cardinals of $\P$,  $\pi:\N\rightarrow \P|\xi$ is an elementary embedding in $\P$ such that $\N$ is countable. It follows that $\eta=_{def}\sup(\pi[\d^\N])<\d^\P$, and therefore, letting $\Lambda=(\pi$-pullback of $S^\P$),\\\\
(O1)  $\P\models ``\Lambda^{stc}$ is a splendid $Ord$-strategy for $\N_0\footnote{The definition of $\N_0$ appears in \rdef{dfn:hod_pm}. The fact that $\Lambda$ has branch condensation follows from generic interpretability. Because $\P|\xi\models ``$the generic interpretation of $S^\P$ has branch condensation", we have the same holds over $\N$.}$",\\
(O2) $\P\models ``\N$ is a $\Lambda^{stc}$-premouse\footnote{\cite[Definition 3.1.8]{hod_mice_LSA} introduces the short-tree-component of an iteration strategy. Roughly speaking $\Lambda^{stc}(\T)=b$ if and only if letting $\Lambda(\T)=c$, either (i) $b=c$, $\pi^\T_c$ is undefined or $\pi^\T_c(\d^\N)>\d(\T)$ or (ii) $b=\m^+(\T)$ and $\pi^\T_c$ is defined and $\pi^\T_c(\d^\N)=\d(\T)$.}.\\
(O3) in $\P$, \rthm{generic interpretability} applies to $\N$ and $\Lambda$.\\
(O4) if $i:\N\rightarrow \N_1$ is such that $\cp(i)>\d^\N$, and for some $\sigma:\N_1\rightarrow \P|\xi$, $\pi=\sigma\circ i$, then $\N_1$ is a $\Lambda^{stc}$-premouse. \\\\
(O4) will be key in many arguments in this paper, but often we will ignore stating it for the sake of succinctness. In each case, however, the reader can easily find the realizable embeddings. The reason (O4) is important is that without it we cannot really prove any self-iterability results, as if iterating $\N$ above destroyed the fact that the resulting premouse is a $\Lambda^{stc}$-premouse then we couldn't find the relevant $\Q$-structures using $\Lambda$ or comparison techniques. 

\begin{lemma}\label{full iterability in excellent hybrid mice} Suppose $\P$ is an excellent hybrid premouse, $w$ is a maximal window of $\P$ above $\d^\P$ and $\eta\in [\nu^w, \d^w)$ is a regular cardinal. Let $\Sigma$ be the $Ord$-strategy of $\Q=\P|\eta$ that acts on iterations that are above $\nu^w$. Then $\Sigma$ is an $Ord$-iteration strategy for $\P$ that acts on iterations that are based on $\Q$ and are above $\nu^w$. 
\end{lemma}
\begin{proof}
We set $V=\P$. Suppose $\T$ is an iteration tree on $\Q$ according to $\Sigma$. We can then naturally regard $\T$ as a tree on $\P$. We claim that all the models of this tree are well-founded. Towards a contradiction assume not. Fix an inaccessible $\xi>\d^w$ such that when regarding $\T$ as a tree on $\P|\xi$, some model of it is illfounded. Let $\T^+$ be the result of applying $\T$ to $\P|\xi$, and let $\pi: \M\rightarrow \P|\xi$ be such that 
\begin{enumerate}
\item $w, \T^+ \in rng(\pi)$,
\item $\card{\M}=\eta$, and
\item $\cp(\pi)>\eta$.
\end{enumerate}  
Let $\U=\pi^{-1}(\T)$ and $\U^+=\pi^{-1}(\T^+)$. We thus have that some model of $\U^+$ is ill-founded.

 Let $\R=\P|\eta^+$ and let $\U^\R$ be the result of applying $\U$ to $\R$. Notice that $\M\in \R$. Because $\P$ has no Woodin cardinals in the interval $(\nu^w, \eta^+)$, we have that $\U^\R$ is according to any $Ord$-strategy of $\R$. Thus $\U^\R$ only has well-founded models. It is not hard to show, however, that for each $\a<lh(\U)$, if $[0, \a]_\U\cap \mathcal{D}^\U=\emptyset$ then there is an elementary embedding $\sigma_\a:\M_\a^{\U^+}\rightarrow \pi^{\U^\R}_{0, \a}(\M)$. In the case $[0, \a]_\U\cap \mathcal{D}^\U\not =\emptyset$, $\M_\a^{\U^+}=\M_\a^\U=\M_\a^{\U^\R}$. 
\end{proof}

In fact more is true.
\begin{lemma}\label{full iterability in excellent hybrid mice generic} Suppose $\P$ is an excellent hybrid premouse, $w$ is a maximal window of $\P$ above $\d^\P$ and $\eta\in [\nu^w, \d^w)$ is a regular cardinal. Let $\Sigma$ be the $Ord$-strategy of $\Q=\P|\eta$ that acts on iterations that are above $\nu^w$. Let $\mathbb{P}\in \P$ be a poset and $g\subseteq \mathbb{P}$ be $\P$-generic. Then $\Sigma$ has a unique extension $\Sigma^g$ definable over $\P[g]$ such that in $\P[g]$, $\Sigma^g$ is an $Ord$-iteration strategy for $\Q$ acting on iterations that are above $\nu^w$.

Moreover, in $\P[g]$, $\Sigma^g$ can be regarded as an $Ord$-iteration strategy for $\P$ that acts on iterations that are based on $\Q$ and are above $\nu^w$.\footnote{The proof of this clause is very similar to the proof of \rlem{full iterability in excellent hybrid mice}.} 
\end{lemma}
\begin{proof} The proof is by now a standard argument in descriptive inner model theory. It has appeared in several publications. For example, the reader can consult the proof of \cite[Lemma 3.9 and Theorem 3.10]{hod_mice} or \cite[Proposition 1.4-1.7]{sargsyan_2022}. We will only give an outline of the proof. 

Fix $\zeta$ such that $\mathbb{P}\in \P|\zeta$. Fix now a maximal window $v$ such that $\nu^v>\max(\zeta, \nu^w)$. Let $(\M_\xi, \N_\xi: \xi<\Omega)$ be the output of the fully background $S^\P$-construction done over $\P|\nu^w$ with critical point $>\nu^v$. Because $\Sigma$ is an $Ord$-strategy, we must have a $\xi<\Omega$ such that $\N_\xi$ is a normal iterate of $\Q$ via an iteration $\T$ that is according to $\Sigma$ and is such that the iteration embedding $\pi^\T:\Q\rightarrow \N_\xi$ is defined.\footnote{I.e., $[0, lh(\T)-1]_\T\cap D^\T=\emptyset$ the final model iteration doesn't drop.}

Assume now that we have determined that an iteration $\U\in \P|\nu^v[g]$ of $\Q$ is according to $\Sigma^g$ and has limit length. For simplicity, let us assume $\U$ has no drops. We want to describe $\Sigma^g(\U)$. Set $\Sigma^g(\U)=b$ if and only if there is $\sigma:\M^\U_b\rightarrow \N_\xi$ such that $\pi^\T=\sigma\circ \pi^\U_b$. To show that this works, we need to show that there is a unique branch $b$ with the desired property. Such a branch $b$ is called $\pi^\T$-realizable.

 Towards a contradiction, assume that either there is no such branch or there are two. Let $\l=((\d^v)^{+2})^\P$. Let now $\pi:\N\rightarrow \P|\l$ be a pointwise definable countable hull of $\P|\l$. It follows that we can find a maximal window $u$ of $\N$, an $\N$-regular cardinal $\zeta\in (\nu^u, \d^u)$, a partial ordering $\mathbb{Q}\in \N$ and a maximal window $z$ of $\N$ such that
 \begin{enumerate}
 \item $\mathbb{Q}\in \N|\nu^z$,
 \item for some $\W$ that is a model appearing in the fully backgrounded construction of $\N|\d^z$ done over $\N|\nu^u$ with respect to $S^\N$ using extenders with critical points $>\nu^z$, there is an iteration $\K\in \N$ on $\R=_{def}\N|\zeta$ with last model $\W$ such that $\pi^\K$ is defined,
 \item some condition $q\in \mathbb{Q}$ forces that whenever $h\subseteq \mathbb{Q}$ is $\N$-generic, there is an iteration $\X\in (\N|\nu^z)[h]$ of $\R$ with no drops such that either there is  no $\pi^\K$-realizable branch or there are at least two $\pi^\K$-realizable branches.
 \end{enumerate}
 Let $h\in \P$ be $\N$-generic for $\mathbb{Q}$. Let $\X\in (\N|\nu^z)[h]$ be as in clause 3 above. Because $\pi(\R)$ is fully iterable in $\P$ above $\pi(\nu^u)$, we have that $\R$ is fully iterable in $\P$ above $\nu^u$. Let $b$ be the branch of $\X$ according to the strategy of $\R$ that is obtained as the $\pi$-pullback of the strategy of $\pi(\R)$ (recall that $\pi(\R)$ is iterable as a $S^\P$-mouse). Because $\R$ has no Woodin cardinals above $\nu^u$, we have a largest $\S\insegeq \M^\X_b$ such that $\S\models ``\d(\X)$ is a Woodin cardinal" but $rud(\S)\models ``\d^\X$ is not a Woodin cardinal".
 
 We claim that\\
 
 \textit{Claim.} $\S\in \N[h]$.
 \begin{proof} To see this, as $\N$ is closed under $\#$, we can assume that $\m^+(\X)\models ``\d(\X)$ is a Woodin cardinal". Let $\V=\m^+(\X)$. We now compare $\V$ with the construction producing $\W$. As $\W$ has no Woodin cardinals above $\nu^u$, we get that there are models $\V^*\insegeq \V^{**}$ appearing on the construction producing $\W$, a tree $\Y$ on $\V$ and a branch $c$ of $\Y$ such that $\V^*=\M^\Y_c$ and $\V^{**}$ is the least model appearing on the construction producing $\W$ such that $\V^{**}\models ``\pi^\Y_c(\d(\X))$ is a Woodin cardinal" but $rud(\V^{**})\models ``\pi^\Y_c(\d(\X))$ is not a Woodin cardinal". It follows that $\S=Hull_n^{\V^{**}}(\{p\}\cup rng(\pi^\Y_c))$ where  $n$ is the fine structural level at which a counterexample to Woodiness of $\d(\X)$ can be defined over $\S$ and $p$ is the $n$-th standard parameter of $\V^{**}$. Because $\V^{**}, \Y, c\in \N[h]$, we have that $\S\in \N[h]$.
 \end{proof} 
 
 It now follows that $b\in \N[h]$, and as $\N$ is pointwise definable, we must have that $b$ is $\pi^\K$-realizable in $\N[h]$ (notice that the argument from the above paragraph implies that $\M^\X_b$ iterates to $\W$). If $d$ is another $\pi^\K$-realizable branch in $\N[h]$ (or in $\P$) then as both $\M^\X_b$ and $\M^\X_d$ are iterable as $S^\N$-mice, we have that $\M^\X_b=\M^\X_d$. This is a contradiction as $\d(\X)$ is not Woodin in either $\M^\X_b$ or $\M^\X_d$.
 
In the case $\U$ has drops, the argument is very similar to the proof of the claim above. In this case, we cannot hope to find a realizable branch, but we can find the appropriate $\Q$-structure using the proof of the claim. 
\end{proof}

Putting the proofs of \rlem{full iterability in excellent hybrid mice} and \rlem{full iterability in excellent hybrid mice generic} together we obtain the proof of \rthm{thm:full iterability in excellent hybrid mice generic}.  \\\\
\textit{The proof of} \rthm{thm:full iterability in excellent hybrid mice generic}. \\\\
We outline the proof. We use the notation introduced in \rthm{thm:full iterability in excellent hybrid mice generic}. Let $\xi$ be a $\P$-inaccessible limit of $\P$-Woodin cardinals and such that $\mathbb{P}\in \P|\xi$, and let $\pi:\N\rightarrow \P|\xi$ be such that $\card{\N}=\eta$, $\cp(\pi)>\eta$ and $\mathbb{P}\in rng(\pi)$. Let $\mathbb{Q}=\pi^{-1}(\mathbb{P})$. Let $h$ be $\P$-generic for $\mathbb{Q}$. Notice now that \rlem{full iterability in excellent hybrid mice generic} applies both in $\N[h]$ and $\P[h]$. Moreover, the proof of \rlem{full iterability in excellent hybrid mice generic} shows that\\\\
(1) $(\Sigma^h)^{\P[h]}\rest (\N[h])=(\Sigma^h)^{\N[h]}$.\\\\
 To see (1), notice that as $\Q$ has no Woodin cardinals, both $(\Sigma^h)^{\P[h]}$ and $(\Sigma^h)^{\N[h]}$ are guided by $\Q$-structures. To see that (1) holds we need to show that both $(\Sigma^h)^{\P[h]}$ and $(\Sigma^h)^{\N[h]}$ pick the same $\Q$-structures, and this would follow if we show that the $\Q$-structures picked by  $(\Sigma^h)^{\N[h]}$ are iterable in $\P[h]$. To see this, we have to recall our definition of $(\Sigma^h)^{\N[h]}$. The iterability of any $\Q$-structure picked by $(\Sigma^h)^{\N[h]}$ is reduced to iterabilty of $\N$ in some non-maximal window $u^\N$. The iterability of this window is reduced to the iterability of $\P$ in some non-maximal window $\pi(u^\N)$, and according to \rlem{full iterability in excellent hybrid mice generic} this last iterability holds.
 
 Finally notice that if we let $\Q^+=\P|(\eta^+)^\P$ and $\Lambda^h$ be the strategy of $\Q^+$ given by  \rlem{full iterability in excellent hybrid mice generic}, $\Lambda^h_\Q=\Sigma^h$ (again this is simply because they both are $\Q$-structure guided strategies). It now just remains to repeat the argument from  \rlem{full iterability in excellent hybrid mice}. Given any tree $\T\in \N[h]$ according to $\Sigma^h$ such that $\pi^\T$ exists, $\pi^\T$ can be applied to $\Q^+$ and hence, to $\N$. This finishes the proof of \rthm{thm:full iterability in excellent hybrid mice generic}. $\qedsymbol{}$
 \subsection{Iterability of countable hulls.}
 
 Here we would like to prove that countable hulls of an excellent hybrid premouse have iteration strategies. The reason for doing this is to show that if $\P$ is an excellent hybrid premouse and $g$ is $\P$-generic then any universally Baire set $A$ in $\P[g]$ is reducible to some iteration strategy which is Wadge below $S^\P$. We will use this to show that $\sf{Sealing}$ holds in a generic extension of an excellent hybrid premouse (see \rthm{thm:upper_bd}).
 
\begin{proposition}\label{citerability in excellent hybrid mice} Suppose $\P$ is an excellent hybrid premouse and $(w_i: i<\omega)$ are infinitely many consecutive  windows of $\P$. Set $\xi=\sup_{i<\omega}\d^{w_i}$. Suppose $\mathbb{P}\in \P|\nu^{w_0}$ is a poset and $g\subseteq \mathbb{P}$ is $\P$-generic.  Working in $\P[g]$, let $\pi: \N\rightarrow \P|(\xi^+)^\P[g]$ be a countable transitive hull. Then in $\P[g]$, $\N$ has a $\nu^{w_0}$-strategy $\Sigma$ that acts on non-dropping trees that are based on the interval $[\pi^{-1}(\nu^{w_0}), \pi^{-1}(\xi)]$. 
\end{proposition}
\begin{proof} Set $u_i=\pi^{-1}(w_i)$ and $\zeta=\pi^{-1}(\xi)$. Our intention is to lift trees from $\N$ to $\P$ and use $\P$'s strategy. However, as $\P$-moves, we lose  \rthm{thm:full iterability in excellent hybrid mice generic}, it is now only applicable inside the iterate of $\P$. To deal with this issue we will use Neeman's ``realizable maps are generic" theorem (see \cite[Theorem 4.9.1]{Neeman}). That it applies is a consequence of the fact that the strategy of $\P$ we have described in \rthm{thm:full iterability in excellent hybrid mice generic} is unique, thus the lift up trees from $\N$ to $\P$ pick unique branches (this is a consequence of Steel's result that UBH holds in mice, see \cite[Theorem 3.3]{LKC}, but can also be proved using methods of  \rthm{thm:full iterability in excellent hybrid mice generic}). One last wrinkle is to notice that when lifting trees from $\N$ to $\P$,  \rthm{thm:full iterability in excellent hybrid mice generic} applies. This is because for each $i$, $\sup(\pi[\d^{u_i}])<\d^{w_i}$. 

We now describe our intended strategy for $\N$. We call this strategy $\Lambda$. Notice that if $\T$ is a normal iteration of $\N$ based on the interval $[\nu^{u_0}, \zeta]$ then $\T$ can be re-organized as a stack of $\omega$-iterations $(\T_i, \N_i: i<\omega)$ where $\N_0=\N$, $\N_{i+1}$ is the last model of $\T_i$ and $\T_{i+1}$ is the largest initial segment of $\T_{\geq \N_i}$ that is based on the window $\pi^{\T_{\leq \N_i}}(u_{i+1})$. 

Suppose then $\T=(\T_i,  \N_i: i<\omega)$ is a normal non-dropping iteration of $\N$ based on $[\nu^{u_0}, \zeta]$. We say $\T$ is according to $\Lambda$ if and only if there is an iteration $\U=(\U_i, \P_i: i<\omega)$ of $\P$ and embeddings $\pi_i:\N_i\rightarrow \P_i$ such that
\begin{enumerate}
\item $\P_0=\P$,
\item $\U_i=\pi_i\T_i$ for each $i<\omega$, 
\item $\P_{i+1}$ is the last model of $\U_i$ for each $i<\omega$,
\item for $i<\omega$, letting $s_i=\pi^\T_{\N, \N_i}(u_i)$, $\l_i=\sup(\pi_i[\d^{s_i}])$ and $\Q_i=\P_i|(\l_i^+)^{\P_i}$, $\P_i[g][\pi_i]\models ``\U_i$ is according to the strategy (as described in \rlem{full iterability in excellent hybrid mice generic}) of $\Q_i$"\footnote{According to \cite[Theorem 4.9.1]{Neeman} the size of the poset that adds $\pi_i$ to $\P_i$ is less than the generators of $\U_{\leq \P_i}$, which is contained in $\pi_i(\nu^{s_i})$.}.
\end{enumerate}
The reader can now use \rthm{thm:full iterability in excellent hybrid mice generic} and Neeman's aforementioned result to show that $\N$ has a $\nu^{w_0}$-iteration strategy. The main point is that: for any $\P[g]$-generic $G \subseteq Coll(\omega,\T)$, in $\P[g][G]$, $\T$ is countable, so we can find generics $g_i$ for each $i$ such that $\pi_i\in \P_i[g][g_i]$. Furthermore, $\Q_i$'s strategy is unique and uniquely extends to all generic extensions of $\P_i$, so the procedure described above can be carried out in $\P[g]$ using the forcing relation of $Coll(\omega, \T)$.
\end{proof}
\subsection{A revised authentication method}\label{subsec:rev authentication}

Suppose $\P$ is an excellent hybrid premouse.  Let $g$ be $\P$-generic. We would like to know if $S^\P$ has a canonical interpretation in $\P[g]$. That this is possible follows from \rthm{generic interpretability}. Perhaps consulting \rrem{important remark} will be helpful. However, to make these notions more precise, we will need to dig deeper into the proof of   \rthm{generic interpretability} and understand how the definition of $\Psi^g$ works. For this we will need to understand expressions such as $``\Q$ is an \textit{authenticated} sts premouse" and etc. The intended meaning of ``authenticated" is the one used in the proof of \cite[Theorem 6.1.5]{hod_mice_LSA}. More specifically, the interested reader should consult \cite[Definition 3.7.3, 3.7.4, 6.2.1 and 6.2.2]{hod_mice_LSA}. Here we will briefly explain the meaning of the expression and state a useful consequence of it that equates this notion to the standard notion of being constructed by fully backgrounded constructions (see \rrem{revised authentication}). The new key concepts are $(\P, \Sigma, X)$-authenticated hybrid premouse and  $(\P, \Sigma, X)$-authenticated iteration. The essence of these two notions are as follows.

\begin{definition}[Authenticated hod premouse]\label{authenticated by a strategy} Suppose $(\P, \Sigma)$ is an sts pair, $X\subseteq \P^b$ and $\R$ is a hod premouse. We say $\R$ is $(\P, \Sigma, X)$-authenticated if there are \\
(e1) a $\Sigma$-iterate $\S$ of $\P$ such that the iteration embedding $\pi:\P\rightarrow \S$ exists and\\
(e2) an iteration $\U$ of $\R$ with last model some $\S^b||\xi$. \\
The iteration $\U$ is constructed using information given by $\pi[X]$. More precisely, for each maximal window $w$ of $\S^b$, consider 
\begin{center}
$s(\pi, X, w)=Hull^{\S^b}(\pi[X]\cup \nu^w)\cap \d^w$.
\end{center}
It is required that for each limit $\a<lh(\U)$, if $c=[0, \a]_\U$ then one of the following two conditions holds:\\\\
(C1) $\S\models ``\d(\U\rest \a)$ is not a Woodin cardinal", $\Q(c, \U\rest \a)$ exists and $\Q(c, \U\rest \a)\insegeq \S$.\\
(C2) $\S\models ``\d(\U\rest \a)$ is a Woodin cardinal" and letting $w$ be the maximal window of $\S$ such that $\d^w=\d(\U\rest \a)$, $s(\pi, X, w)\subseteq rng(\pi^{\U\rest \a}_c)$.\\\\
Usually, $X$ is chosen in a way that for each window $w$ of $\S$, $\sup(s(\pi, X, w))=\d^w$. For such $X$, conditions (C1) and (C2) completely determine $\U$.

Given a hod premouse $\P$ of lsa type, a set $X\subseteq \P^b$ and a set $\Gamma$ consisting of iterations of $\P$ we can similarly define $(\P, \Gamma, X)$-authenticated hod premice.
\end{definition}

\begin{definition}[Authenticated iteration]\label{authenticated iteration}
Suppose $\R$ is a $(\P, \Sigma, X)$ authenticated hybrid premouse and $\W$ is an iteration of $\R$. We say $\W$ is $(\P, \Sigma, X)$-authenticated if there is a triple $(\S, \U, \xi)$  $(\P, \Sigma, X)$-authenticating $\R$ such that $\pi^\U$-exists and $\W$ is according to $\pi^\U$-pullback of $\Sigma_{\S||\xi}$. 
\end{definition}

Suppose now that $\M$ is an sts premouse based on $\P$ and $g$ is $\M$-generic for a poset in $\M|\zeta$. Suppose $\R\in \M|\zeta[g]$ is an lsa type hod premouse such that $\R^b$ is $(\P, \S^\M, \P^b)$-authenticated and $\R=(\R|\d^\R)^\#$. In $\M[g]$, we can build an sts premouse $\W$ based on $\R$ using $(\P,  S^\M, \P^b)$-authenticated iterations. This means that whenever $\U$ is an iteration indexed in $\W$, $\a<lh(\U)$ is a limit ordinal such that $\pi^{\U\rest \a, b}$ exists and $\X$ is the longest initial segment of $\U_{\geq \a}$ that is based on $\V=_{def}(\M_\a^\U)^b$, then both $\V$ and $\X$ are $(\P, S^\M, \P^b)$-authenticated. In addition to the above, we also require that if $\Q$ is a $\Q$-structure for some $\sf{nuvs}$ tree in $\W$ that has been authenticated by $\W$ via the authentication procedure used in sts premice then any iteration indexed in $\Q$ is $(\P, \P^b, S^\M)$-authenticated. Moreover, the same holds for all iterates of $\Q$ via the strategy witnessing that $\Q$ is authenticated in $\W$. 

It is important to keep in mind that the above construction may fail simply because some non-$(\P, S^\M, \P^b)$-authenticated object has been constructed. Also, the same construction can be done using  $(\P, S^\M, X)$-authenticated objects where $X\subseteq \P^b$. 

\begin{remark}\label{revised authentication} Suppose now that $\M$ is an sts premouse based on $\P$ and $g$ is $\M$-generic for a poset in $\M|\zeta$. Suppose $\R\in \M|\zeta[g]$ is an lsa type hod premouse such that $\R^b$ is $(\P, \P^b, S^\M)$-authenticated and $\R=(\R|\d^\R)^\#$. Suppose $\M$ has a Woodin cardinal $\d$ above $\zeta$. To say that an sts premouse $\Q$ over $\R$ is $(\P, \P^b, S^\M)$-authenticated is equivalent to saying that $\Q\insegeq \W$ where $\W$ is a model in the $(\M,\P^b)$-authenticated fully backgrounded construction described in \cite[Definition 6.2.2]{hod_mice_LSA}.
\end{remark}
The reader maybe wondering why it is enough to only authenticate the lower level iterations. The reason is that in many situations the lower level strategies define the entire short tree strategy. This point was explained in \rsubsec{subsec: short tree strategy mice}.

The construction mentioned in \rrem{revised authentication} is called the $(\P, X, S^\M)$-authenticated hod pair construction over $\R$. The details of everything that we have said above appears in \cite[Chapter 6.2]{hod_mice_LSA}. The reader may choose to consult \cite[Definition 6.2.2]{hod_mice_LSA}.

\subsection{Generic Interpretability}

In this portion of the current section, we would like to outline the proof of generic interpretability. As was mentioned before, generic interpretability is somewhat tricky for short tree strategies. This is because given a $\Lambda$-sts premouse $\N$ and a tree $\T\in \N[h]$, $\T$ maybe short but $\N[h]$ may not be able to find the branch of $\T$ that is according to $\Lambda$, as this branch might have a $\Q$-structure that is more complex than $\N$.

Suppose $\P$ is an excellent hybrid premouse. For the purposes of this paper, we say that $\P$ satisfies \textit{weak generic interpretability} if for every poset $\mathbb{P}\in \P$ and for every $\P$-generic $g\subseteq \mathbb{P}$,  there is an sts strategy $\Lambda$ for $\P_0$\footnote{$\P_0$ was introduced in \rdef{dfn:hod_pm}.} that is definable (with parameters) over $\P[g]$ such that for every tree $\T\in dom(\Lambda)$,
\begin{enumerate}
\item if $\T$ is $\sf{uvs}$ then letting $\Lambda(\T)=c$, either 
\begin{enumerate}
\item for some node $\R$ of $\T$ such that $\pi^{\T_{\leq \R}, b}$ is defined, $\T_{\geq \R}$ is a tree on $\R^b$ and $\T^\frown\{c\}$ is $(\P_0, \P_0^b, S^\P)$-authenticated, or 
\item for some node $\R$ of $\T$ such that $\pi^{\T_{\leq \R}, b}$ is defined, $\T_{\geq \R}$ is a tree on $\R^b$ that is above $Ord\cap \R^b$, $\Q(c, \T)$ exists and $\Q(c, \T)\insegeq \m^+(\T)$,
\end{enumerate}
\item if $\T$ is $\sf{nuvs}$ then letting $c=\Lambda(\T)$, $c$ is a cofinal branch if and only if $\Q(c, \T)$ exists and $\T^\frown\{c\}$ is $(\P_0, \P_0^b, S^\P)$-authenticated.
\end{enumerate}

\begin{proposition}\label{generic interpretability in ehpm} Suppose $\P$ is an excellent  hybrid premouse. Then $\P$ satisfies  weak generic interpretability.  
\end{proposition}
\begin{proof} We outline the proof as the proof is very much like the proof of 
\cite[Theorem 6.1.5]{hod_mice_LSA}. Let $g$ be $\P$-generic. The definition of $\Lambda$ essentially repeats the above clauses.  We first consider trees that are $\sf{uvs}$. 

Suppose $\T$ is an $\sf{uvs}$ tree according to $\Lambda$, and suppose that for some node $\R$ on $\T$, $\pi^{\T_{\leq \R}, b}$ exists and $\T_{\geq \R}$ is a tree based on $\R^b$. Because $\T$ is according to $\Lambda$, we may assume that $\R$ is $(\P_0, \P_0^b, S^\P)$-authenticated. Thus, we can fix a window $w$ of $\P$ such that $g$ is $<\nu^w$-generic over $\P$ and letting $\W$ be the iterate of $\P_0$ constructed by the fully backgrounded hod pair construction of $\P|\d^w$ using extenders with critical point $>\nu^w$, we can find an embedding $\sigma:\R^b\rightarrow \W^b$ such that \\\\
(a) $\pi^{\U, b}=\sigma\circ \pi^{\T_{\leq \R}, b}$ where $\U$ is the $\P_0$-to-$\W$ tree according to $S^\P$ and\\
(b) $\T_{\geq \R}$ is according to the $\sigma$-pullback of $S^\P_{\W^b}$. \\\\
Letting $c$ be the branch according to the strategy as in (b), we have that $\T^\frown\{c\}$ is $(\P_0, \P_0^b, S^\P)$-authenticated.
Moreover, there is only one such branch $c$. To see this, we need to reflect. 

Let $\xi$ be large and let $\pi:\N\rightarrow \P|\xi$ be a countable hull. Fix an $\N$-generic $h\in \P$. Let $\U\in \N[h]$ be an  $\sf{uvs}$ tree on $\N_0$ such that for some node $\S$ on $\U$ with the property that $\pi^{\U, b}$ exists, $\U_{\geq \S}$ is a normal tree on $\S^b$, and moreover, $\U$ is $(\N_0, \N_0^b, S^\N)$-authenticated in $\N[h]$. Suppose now that there are two distinct branches $c$ and $d$ obtained in the above manner. We can then fix $\N$-windows $u_1, u_2$ that play the role of $w$ above and build $\K_1$ and $\K_2$, the equivalents of $\W$ above, inside $\N|\d^{u_1}$ and $\N|\d^{u_2}$. For $i\leq 2$, we  have maps $m_i:\S^b\rightarrow \K^b_i$  and  $S^\N$-iteration maps $\tau_i:\N_0^b\rightarrow \K^b_i$ such that $\tau_i=m_i\circ \pi^{\U_{\leq \S}, b}$. Let $\tau_i: \K^b_i \rightarrow \mathcal{Y}$ be the comparison map using strategies $(S^\N)_{\K^b_i}$. Then, for $i\in \{c,d\}$ there is an embedding $\l_i: \M^\U_i \rightarrow \mathcal{Y}$ that factors into the iteration map from $\N_0^b$ to $\mathcal{Y}$. It is then easy to see, using branch condensation of $S^\P$ and its $\pi$-pullback, that $c = d$.

The rest of the argument is very similar. For example, we outline the proof of clause 2 in the definition of weak generic interpretability. Suppose $\T\in \P[g]$ is $\sf{nuvs}$ and $\Q$ is a $(\P_0, \P_0^b, S^\P)$-authenticated $\Q$-structure for $\T$. We want to see that there is a cofinal well-founded branch $c\in \P[g]$ such that $\Q(c, \T)=\Q$. As above, instead of working with $\P$, we can work with a reflection. Thus, we assume that $\pi:\N\rightarrow \P|\xi$ is a countable elementary embedding, $h\in \P$ is $\N$-generic and $\T, \Q\in \N[h]$. Moreover, we can assume that $\N$ is pointwise definable. Let $\Psi$ be the $\pi$-pullback of $S^\P$. As $\Q$ is $(\N_0, \N_0^b, S^\N)$-authenticated, it follows from \rthm{citerability in excellent hybrid mice} that $\Q$ has an iteration strategy as a $\Psi_{\m^+(\T)}$-sts premouse. Let $c$ be the branch of $\T$ according to $\Psi$. As $\N$ is pointwise definable we have that $\Q(c, \T)$-exists, and hence $\Q(c, \T)=\Q$. 
 \end{proof}
 
 Next we show that the low level strategies are in fact universally Baire. However, \rprop{claim:not_uB} shows that $S^\P$ itself does not have a universally Baire representation. 
 
 \begin{proposition}\label{lls ub} Suppose $\P$ is excellent, $g$ is $\P$-generic and $\Sigma$ is the generic interpretation of $S^\P$ onto $\P[g]$. Let $\T\in \P[g]$ be an iteration tree on $\P$ of length $<\omega_1^{\P[g]}$ such that $\pi^{\T, b}$-exists. Set $\R=\pi^{\T, b}(\P^b)$. Then $(\Sigma_\R\rest HC^{\P[g]})\in \Gamma^\infty_g$. Moreover, for any $\P$-cardinal $\eta$ there are $\eta$-complementing trees $T, S\in \P[g]$ such that for all posets $\mathbb{Q}\in \P[g]$ such that $\card{\mathbb{Q}}^{\P[g]}<\eta$ and for all $\P[g]$-generic $h\subseteq \mathbb{Q}$, \begin{center}$(p[T])^{\P[g*h]}=\Sigma_\R^h\rest HC^{\P[g*h]}$.\end{center}  \end{proposition}
 \begin{proof} We again outline the proof as the proof uses standard ideas. Let $w$ be a maximal window of $\P$ such that $g$ is generic for a poset in $\P|\nu^w$. We now outline the construction producing $\nu^w$-complementing trees $(T, S)$ as in the statement of the proposition. 
 
 Let $\S$ be the model appearing on the hod pair construction of $\P|\d^w$ in which extenders used have critical points $>\nu^w$ and to which $\R$ normally iterates via $\Sigma_\R$. Let $i: \R\rightarrow \S$ be the iteration embedding. What we need to show is that club many hulls of $\P[g]$ are correct about $\Sigma_\R$, where we take $\Sigma_\R$ to be defined as $i$-pullback of the strategy of $\S$ that $\S$ inherits from $\P$ (see \rthm{full iterability in excellent hybrid mice generic}). That this works follows from the fact that the strategy of $\S$ has hull condensation. Let $\Psi$ be the strategy of $\S$. 
 
 More precisely, let $\phi(x, \R, \S, i)$ be the formula that says ``$x\in \bR$ codes an iteration of $\R$ that is according to the $i$-pullback of $\Psi$". Clearly, $\phi$ defines $\Sigma_\R\rest HC^{\P[g]}$. Let now $\xi$ be large and  $\pi:\N\rightarrow  \P|\xi[g]$ be countable such that $\R, (i, \S)\in rng(\pi)$. Let  $\Phi=\pi^{-1}(\Psi)$ and $j=\pi^{-1}(i)$. Let $h\in \P[g]$ be a $<\pi^{-1}(\nu^w)$-generic over $\N$, and let $\U\in \N[g]$ be a tree on $\R$. 
 
 Suppose first that $\N[h]\models ``j\U$ is according to the strategy of $\Phi^h$". Because $\Phi^h$ is the $\pi$-pullback of $\Psi$, we have that $i\U=\pi(j\U)$ is according to $\Psi$. Hence, $\U$ is according to $\Sigma_\R$. 
 
 Next suppose that $\U$ is according to $\Sigma_\R$. It then follows by the above reasoning that $\N[h]\models ``j\U$ is according to the strategy of $\Phi^h$". This finishes the proof that $\Sigma_\R$ has a uB representation. The rest of the proposition follows from the fact that the formula $\phi$ above defines $\Sigma_\R$ in all $<\nu^w$-generic extensions of $\P[g]$. Such calculations were carried out more carefully in \cite{hod_mice} and also in \cite{sargsyan_2022}. In particular, the reader may wish to consult \cite[Proposition 1.4]{sargsyan_2022}.
 \end{proof}
 
 \subsection{Fully backgrounded constructions inside excellent hybrid premice}
 
Given an excellent hybrid premouse $\P$, we would eventually like to show that collapsing $((\d^\P)^+)^\P$ to be countable forces both $\sf{Sealing}$ and $\sf{LSA-over-uB}$. Such an analysis of generic extensions of fine structural models usually requires some kind of re-constructibility property, which guarantees that the model can somehow see versions of itself inside it. In this subsection, we would like to establish some such facts about excellent hybrid premice. \rprop{fully backgrounded constructions} is a key proposition that we will need in this paper. Recall the definition of $\P_0$ from \rdef{dfn:hod_pm}.

\begin{proposition}\label{fully backgrounded constructions} Suppose $\P$ is excellent and $g$ is $\P$-generic. Let $\Sigma$ be the generic interpretation of $S^\P$ onto $\P[g]$, and suppose $\R$ is a $\Sigma$-maximal iterate of $\P_0$. Let $w$ be a maximal window of $\P$ such that $w$ is above $\d^\P$, $g$ is generic for a poset in $\P|\nu^w$ and $\R\in \P|\nu^w[g]$. Let $\xi<\d^{\R^b}$ be a Woodin cardinal of $\R$. Suppose $\N_0$ is the output of the fully backgrounded hod pair construction of $\P|\d^w[g]$ done relative to $\Sigma_{\R|\xi}$ and over $\R|\xi$ and using extenders with critical points $>\nu^w$. Then $Ord\cap \N_0=\d^w$.
\end{proposition}
\begin{proof}
Towards a contradiction, suppose that $\eta=_{def}Ord\cap \N_0<\d^w$. We will now work towards showing that $\eta$ is a Woodin cardinal of $\P$. As $\eta\in (\nu^w, \d^w)$, this is clearly a contradiction. Suppose then $\eta$ is not a Woodin cardinal of $\P$. As $\P$ has no Woodin cardinals in the interval $(\nu^w, \d^w)$, we must have that there is a $\Sigma$-mouse $\P|\eta\insegeq \Q\insegeq \P$ such that  $\eta$ is a cutpoint of $\Q$, $\Q\models ``\eta$ is a Woodin cardinal" but $\eta$ is not Woodin relative to functions definable over $\Q$. Unfortunately $\Q$ cannot be translated into a $\Sigma_{\N_0}$-mouse, but we can rebuild it in a sufficiently rich model extending $\N_0$.

Let $\N$ be the output of a fully backgrounded construction of $\P|\d^w[g]$ done with respect to $\Sigma_{\N_0}$ and over $\N_0$ using extenders with critical point $>\eta$. As $\N_0$ is a $\Sigma$-maximal iterate of $\P_0$\footnote{Suppose $\N_0$ is not $\Sigma$-maximal and let $\U$ be the $\P_0$-to-$\N_0$ tree. Let $b=\Sigma(\U)$. We then have that $\Q(b, \U)$ exists and so $\N_0$ couldn't be the final model of the fully backgrounded hod pair construction of $\P|\d^w[g]$. It follows from the universality of the fully backgrounded constructions, that continuing the construction further we will construct $\Q(b, \U)$. The reader may wish to consult \cite[Chapter 4]{hod_mice_LSA}.}, we have that $\N\models ``\eta$ is  a Woodin cardinal". We now want to rebuild $\Q$ inside $\N[\P|\eta]$. The idea here goes back to \cite[Theorem 8.1.13]{hod_mice_LSA} (for instance, the construction of $\N_2$ in the proof of the aforementioned theorem.). Notice that if $p$ is the $\P_0$-to-$\N_0$-iteration\footnote{$p$ is a stack of two normal trees.}  then $\pi^{p, b}$ exists and $\pi^{p, b}\in \N[\P|\eta]$. Let $X=\pi^{p, b}[\P_0^b]$. Working inside $\N[\P|\eta]$ we can build a $\Sigma$-premouse over $\P|\eta$ via a fully backgrounded $(\N_0, X, S^\N)$-authenticated construction. In this construction we only use extenders with critical point $>\eta$. Let $\W$ be the output of this construction. As $\W$ is universal, we have that $\Q\insegeq \W$. Thus, $\N[\P|\eta]\models ``\eta$ is not a Woodin cardinal".

However, standard arguments show that $\N[\P|\eta]\models ``\eta$ is a Woodin cardinal". Indeed, let $f:\eta\rightarrow \eta$ be a function in $\N[\P|\eta]$. Because $\P|\eta$ is added by an $\eta$-cc poset, we can find $g\in \N$ such that for every $\a<\eta$, $f(\a)<g(\a)$. Let $E\in \vec{E}^{\N_0}$  be any extender  witnessing Woodiness for $g$ and such that $\N\models ``\nu_E$ is a measurable cardinal". Thus, $\pi^{\N_0}_E(g)(\kappa)<\nu_E$, where $\kappa$ is the critical point of $E$. Let $F$ be the resurrection of $E$. We must have that $\pi^\P_F(f)(\kappa)<\nu_E$. Thus, $F\rest \nu_E\in \P|\eta$ is an extender witnessing Woodinness for $f$ in $\P|\eta$ and hence in $\N[\P|\eta]$. 
\end{proof}
 
Using \rprop{fully backgrounded constructions}, we can now prove that  $S^\P$ itself is not a universally Baire set. Its proof requires a few more facts from \cite{hod_mice_LSA}, which we now review. Given an lsa type pair $(\P, \Sigma)$, following \cite[Definition 3.3.9]{hod_mice_LSA}, we let $\Gamma^b(\P, \Sigma)$ be the set of all $A\subseteq \bR$ such that for some countable iteration $\T$ such that $\pi^{\T, b}$ exists, $A$ is Wadge reducible to $\Sigma_{\pi^{\T, b}(\P^b)}$. The following comparison theorem is essentially \cite[Theorem 4.13.1]{hod_mice_LSA}. 
 
 \begin{theorem}\label{comparison thm} Assume $\sf{AD}^+$ and suppose $(\P, \Sigma)$ and $(\Q, \Lambda)$ are two lsa type hod pairs such that $\Gamma^b(\P, \Sigma)=\Gamma^b(\Q, \Lambda)$ and both $\Sigma$ and $\Lambda$ are splendid. Then there is an lsa type hod pair $(\R, \Psi)$ such that $\R$ is a $\Sigma$-iterate of $\P$ and a $\Lambda$-iterate of $\Q$ and $\Sigma_\R=\Psi=\Lambda_\R$. 
 \end{theorem}
 
 \begin{proposition}\label{claim:not_uB} Suppose $\P$ is excellent and $g$ is $\P$-generic. Let $\Sigma$ be the generic interpretation of $S^\P$ onto $\P[g]$. 
Then $\P[g] \vDash \Sigma \notin \Gamma^\infty$.
\end{proposition}
\begin{proof}
Towards a contradiction, suppose that $\P[g]\models \Sigma^g \in \Gamma^\infty$. Let $(w_i: i<\omega)$ be a sequence of successive windows of $\P$ such that $g$ is generic over a poset in $\P|\nu^{w_0}$. Set $\d_i=\d^{w_i}$. 

For each $i$, let $\P_i$ be the $S^\P$-iterate\footnote{$S^\P$ is the internal strategy predicate of $\P$, which by itself is not a total iteration strategy but can be uniquely extended to a total iteration strategy. By $``S^\P$-iterate" we mean an iterate according to the total extension of $S^\P$. The reader may consult \cite[Chapter 5]{hod_mice_LSA}.} built via the fully backgrounded hod pair construction of $\P|\d_i$ using extenders with critical points $>\nu_i^{w_i}$. It follows from \rprop{fully backgrounded constructions} that \\\\
(1) $\d^{\P_i}=\d_i$, i.e., $\d_i$ is the largest Woodin cardinal of $\P_i$.\\\\
Fix now $k\subseteq Coll(\omega, <\d_\omega)$ generic over $\P[g]$ where $\d_{\omega}=\sup_{i<\omega}\d_i$. 
Recall next that Steel showed that if there are unboundedly many Woodin cardinals then every universally Baire set has a universally Baire scale (see \cite[Theorem 4.3]{DMT}\footnote{Recall that by a result of Martin, Steel and Woodin for a $\lambda$ a limit of Woodins, $Hom_{<\l}$ coincides with the $<\l$-universally Baire sets. See \cite[Theorem 2.1]{DMT} and \cite[Chapter 2]{DMT}.}). Let now $W$ be the derived model of $\P$ as computed in $\P[m]$ where $m=g*k$. It follows that\\\\
(2) the canonical set of reals coding $\Sigma^m\rest HC^W$\footnote{We will identify Code$(\Sigma)$ with $\Sigma$ itself in this paper.} has a scale in $W$. \\\\
This paragraph will be using Theorem \ref{comparison thm} and the notation introduced there. Working inside $W$, let $\Gamma=\Gamma^b(\P_0, \Sigma^m)$. Thus, $\Gamma$ is the set of reals that are generated by the low-level-components of $\Sigma^m$. More precisely, $A\in \Gamma$ if there is an iteration $\T$ on $\P_0$ according to $\Sigma^m$ such that $\pi^{\T, b}$ exists and $A$ is Wadge below $\Sigma^m_{\pi^{\T, b}(\P_0^b)}$.  As $\Sigma^m\rest HC^W$ is Suslin, co-Suslin in $W$, we must have a hod pair $(\S, \Lambda)\in W$ such that $\Gamma^b(\S, \Lambda)=\Gamma$ (this follows from the Generation of Mouse Full Pointclasses, see \cite[Theorem 10.1.2]{hod_mice_LSA}). We can further assume that $\S$ is a $\Sigma^m$-iterate of $\P_0$ and $\Lambda^{stc}=\Sigma^m_{\S}\rest HC^W$ (this extra possibility follows from Theorem \ref{comparison thm}). 

Fix now $i<\omega$ such that letting $n=k\cap Coll(\omega, \d_i)$, $\P[g*n]$ has a uB representation of $\Lambda$. It now follows that since \\\\
(3) $\P_{i+1}$ is a $\Lambda$-iterate of $\S$ and\\
(4) letting $l:\S\rightarrow \P_{i+1}$ be the iteration embedding, $l[\d^\S]$ is cofinal in $\d^{\P_{i+1}}$, \\\\
we have
\begin{center}
$\d^{\P_{i+1}}<\d_{i+1}$. 
\end{center}
This directly contradicts (1).
 \end{proof}
 
 \subsection{Constructing an iterate via fully backgrounded constructions}\label{backgrounded iterate}
 
 Suppose $\M$ is strategy-hybrid $\eta$-iterable mouse such that $\M\in V_\eta$, $\eta$ is an inaccessible cardinal and $\M$ has an $\eta$-strategy with hull condensation. Thus $\M$ has an extender sequence $\vec{E}$ and a strategy predicate $S^\M$, which can be a strategy of $\M$ itself (as in hod mice) or a strategy of some $\N\in \M$. We want to build an iterate $\X$ of $\M$ such that the extenders of $\X$ are all fully backgrounded. Here we describe this construction. 

We say $(\V_\xi, \W_\xi, \T_\xi, \X_\xi: \xi<\iota)$ are the models and iterations of the fully backgrounded $(\M, \Sigma)$-iterate-construction of $V_\eta$ if the following conditions are satisfied with $\a_\xi=Ord\cap\V_\xi$.
\begin{enumerate}
\item $\V_0=\W_0=J_0^\M$.
\item For every $\xi<\iota$, $\V_\xi=\W_\xi|\a_\xi$.
\item For $\xi<\iota$, $\T_\xi$ is an iteration of $\M$ according to $\Sigma$ with last model $\X_\xi$ such that $\V_\xi=\X_\xi|\a_\xi$ and the generators of $\T_\xi$ are contained in $\a_\xi$\footnote{There is only one such iteration $\T_\xi$.}.
\item  For $\xi<\iota$, if $\a_\xi\in dom(\vec{E}^{\X_\xi})\cup dom(S^{\X_\xi})$ then $\W_\xi= \X_\xi||\a_\xi$.
\item For $\xi<\iota$, if $\a_\xi\not \in (dom(\vec{E}^{\X_\xi})\cup dom(S^{\X_\xi}))$ and $\V_\xi\not =\X_\xi$ then $\W_\xi=J_1(\V_\xi)$.
\item For $\xi<\iota$, if there is a total extender $F\in V_\eta$ such that
\begin{center}
$\pi_F((\V_\zeta, \W_\zeta, \T_\zeta, \X_\zeta: \zeta<\xi))\rest \xi=(\V_\zeta, \W_\zeta, \T_\zeta, \X_\zeta: \zeta<\xi)$,
\end{center}
then $F\cap \V_\xi=E_{\a_\xi}^{\X_\xi}$. It follows that $\W_\xi=(\V_\xi, E_{\a_\xi}^{\X_\xi})$.
\item For $\xi+1<\iota$, $\V_{\xi+1}=\mathcal{C}(\W_\xi)$.
\item If $\xi<\iota$ is limit then $\W_\xi=liminf_{\zeta\rightarrow \xi}\W_\zeta$ and $\V_\xi=C(\W_\xi)$. More precisely, given $\W_\xi|\kappa$, $\W_\xi|(\kappa^+)^{\W_\xi}$ is the eventual value of $\W_\zeta|(\kappa^+)^{\W_\zeta}$. 
\end{enumerate}
 
 We then let $\sf{FBIC}(\M, \Sigma, \eta)$ be the models and iterations of the above construction. We can vary this construction in two ways. The first way is that fixing some $\l<\eta$ we can require that the extender $F$ in clause 6 has critical point $>\l$. This amounts to backgrounding extenders via total extenders that have critical points $>\l$. The second way is that we may choose to start the construction with any initial segment of $\M$. More precisely, given a cardinal cutpoint $\nu$ of $\M$, we can start by setting $\V_0=\M|\nu$. 
 
 Thus, by saying that $(\V_\xi, \W_\xi, \T_\xi, \X_\xi: \xi<\iota)$ are the models and iterations of the $\textsf{FBIC}(\M, \Sigma, \eta, \l, \nu)$ we mean that the sequence is built as above but starting with $\M|\nu$ and using backgrounded extenders that have critical points $>\l$. 
 
$\textsf{FBIC}(\M, \Sigma, \eta, \l, \nu)$ can break without reaching its eventual goal. We say 
\begin{center}
$\textsf{FBIC}(\M, \Sigma, \eta, \l, \nu)$ 
\end{center}
is \textit{successful} if one of the following conditions holds.
\begin{enumerate}
\item $\iota=\xi+1$, $\pi^{\T_\xi}$ exists and either $\V_\xi=\X_\xi$ or $\W_\xi=\X_\xi$,
\item $\iota$ is a limit ordinal and $liminf_{\xi\rightarrow \iota}\V_\xi$ is the last model of a normal $\Sigma$-iteration $\T$ of $\M$ such that $\pi^\T$ exists. 
\end{enumerate}
If $\textsf{FBIC}(\M, \Sigma, \eta, \l, \nu)$ is successful, then we say $\N$ is its output if it is the iterate of $\M$ described above.

The following is the main theorem that we will need from this section. We say $E$ is a \textit{strictly short extender} if its generators are bounded below $\pi_E(\cp(E))$. We say $\M$ is \textit{strictly short} if all of its extenders are strictly short. 
\begin{theorem}\label{success of fbic} Suppose $(\M, \Sigma)$ and $\eta$ are as above and in addition to the above data, $\eta$ is a Woodin cardinal and $\M$ is strictly short. Suppose $\l<\eta$ and $\nu$ is a cutpoint cardinal of $\M$. Then $\textsf{FBIC}(\M, \Sigma, \eta, \l, \nu)$ is successful. 
\end{theorem}
The proof is a standard combination of universality (see \cite[Lemma 11.1]{DMATM}) and stationarity (see \cite[Lemma 3.23]{scales_hybrid_mice}) of fully backgrounded constructions.

\section{An upper bound for $\sf{Sealing}$ and $\sf{LSA-over-uB}$}\label{sec:upperbound1}
The goal of this section is to prove \rthm{thm:upper_bd}. It reduces $\sf{Sealing}$, $\sf{Tower \ Sealing}$, and $\sf{LSA-over-uB}$ to a large cardinal theory. This essentially constitutes one half of \rthm{thm:main_theorem} and \rthm{thm:tower_sealing}. 
\begin{theorem}\label{thm:upper_bd} Suppose $\P$ is excellent and $g\subseteq Coll(\omega, \d^\P)$ is $\P$-generic. Then  both $\sf{Sealing}$ and $\sf{LSA-over-uB}$ hold in $\P[g]$. 
\end{theorem} 

We start the proof of \rthm{thm:upper_bd}. Let $\P$ be excellent (see \rdef{dfn:hod_pm}). Set $\d_0=\d^\P$
and let $g\subseteq Coll(\omega,\delta_0)$ be $\P$-generic. We first show that $\sf{Sealing}$ holds in $\P[g]$.  Let $\P_0=_{def}(\P|\d_0)^{\#}$. We write $\P = (|\P|,\in, \mathbb{E}^\P, S^\P)$ where $\mathbb{E}^\P$ is the extender sequence of $\P$ and $S^\P$ is the predicate coding the short-tree strategy of $\P_0$ in $\P$. Thus, $\P$ above $\d_0$ is a short tree strategy premouse over $\P_0$.

 Let $\Sigma^-$ be this short tree strategy. It follows from \rprop{generic interpretability in ehpm} that for any $\P[g]$-generic $h$, $\Sigma^-$ has a canonical extension $\Sigma^h$ in $\P[g*h]$\footnote{It is not correct to say that $\Sigma\in \P$. The correct language is that $\Sigma$ is a definable class of $\P$ and $\Sigma^g$ is a definable class of $\P[g]$.}. Let then $\Sigma$ be the extension of $\Sigma^-$ in $\P[g]$. 
 

%

\subsection{An upper bound for $\sf{Sealing}$}\label{sec:SealingBounds}

Let $h$ be $\P[g]$-generic. Working in $\P[g][h]$, let $\Delta^h=\Gamma^b(\P_0, \Sigma)$. Equivalently, $\Delta^h$ is the set of reals $A\subseteq \bR^{\P[g*h]}$ such that for some countable tree $\T$ on $\P_0$ with last model $\Q$ such that $\pi^{\T, b}$ exists, $A\in L(\Sigma^h_{\Q^b}, \bR^{\P[g*h]})$. It follows from \rprop{lls ub} that  if $\Q$ is as above then $\Sigma^h_{\Q^b}\in \Gamma^\infty_{g*h}$.

\begin{lemma} \label{lem:Gamma}
$\Gamma^\infty_{g*h}=\Delta^h$.
\end{lemma} 
\begin{proof}  It follows from \rprop{lls ub} that $\Delta^h\subseteq \Gamma^\infty_{g*h}$. Fix $A\subseteq \bR^{\P[g*h]}$ that is a universally Baire set in $\P[g*h]$. Work in $\P[g*h]$, and suppose $A\not \in \Delta^h$. Because there is a proper class of Woodin cardinals, any two universally Baire sets are Wadge comparable. Since $\Delta^h\cup\{A\}\subseteq \Gamma^\infty_{g*h}$ and $A\not \in \Delta^h$, we have that $\Delta^h\subseteq L(A, \bR_{g*h})$. Recall that by a result of Steel (\cite[Theorem 4.3]{DMT}), $A$ is Suslin, co-Suslin in $\Gamma^\infty_{g*h}$. Hence, we can assume, without losing generality, that in $L(A, \bR_{g*h})$, there are Suslin, co-Suslin sets beyond $\Delta^h$. 

 It follows from \cite[Theorem 10.1.1]{hod_mice_LSA} that there is a lsa type hod pair $(\S, \Lambda)\in \Gamma^\infty_{g*h}$ such that $\Gamma^b(\S, \Lambda)=\Delta^h$. Just like in the proof of \rprop{claim:not_uB} we can assume that $\S$ is a $\Sigma^h$-iterate of $\P_0$ and $\Lambda^{stc}=\Sigma^h_\S$. It then follows that $\Sigma^h\in \Gamma^\infty_{g*h}$, contradicting \rprop{claim:not_uB}. 

\end{proof}

By the results of \cite[Section 8.1]{hod_mice_LSA} (specifically, \cite[Theorem 8.1.1 clause 4]{hod_mice_LSA}), we get that in $\P[g*h]$,
\begin{equation}\label{eqn:closure}
\Delta^h = \powerset(\mathbb{R}_{g*h})\cap L(\Delta^h, \bR_{g*h}).
\end{equation}
The lemma and (\ref{eqn:closure}) immediately give us clause (1) of $\sf{Sealing}$. For clause (2), let $h$ be $\P[g]$-generic and $k$ be $\P[g*h]$-generic. We want to show that there is an elementary embedding
\begin{center}
$j:L(\Delta^h, \bR_{g*h})\rightarrow L(\Delta^{h*k}, \bR_{g*h*k})$
\end{center}
such that for every $A\in \Delta^h$, $j(A)$ is the canonical extension of $A$ in $\P[g*h*k]$. This will be accomplished in \rlem{elem embedding}. The next lemma provides a key step in the construction of the desired elementary embedding. It does so by realizing $L(\Delta^h, \bR_{g*h})$ as a derived model of an iterate of $\P_0$. 

Recall from \cite[Definition 2.7.2]{hod_mice_LSA} that if $\S$ is a hod premouse of limit type (including lsa type) then $\d^{\S^b}$ is the supremum of the Woodin cardinals of $\S^b$. In general, the reader may wish to review some of the notation concerning hod premice: the relevant notation can be found in \cite[Chapter 2 and 3]{hod_mice_LSA}. The Key Phenomenon stated before \cite[Definition 2.7.8]{hod_mice_LSA} might also be useful. 

\begin{lemma}\label{technical lemma} Suppose $\mathbb{P}\in \P[g]$ is a poset and $m\subseteq \mathbb{P}$ is $\P[g]$-generic. Suppose further that in $\P[g*m]$, $\S$ is a countable $\Sigma^m$-iterate of $\P_0$ such that the $\P_0$-to-$\S$ iteration embedding exists. Suppose $\k<\d^{\S^b}$ is a Woodin cardinal of $\S$ and $A\in \Gamma^{\infty}_{g*m}$ (in $\P[g*m]$). Then in $\P[g*m]$, there is a countable $\Sigma^m_\S$-iterate $\W$ of $\S$ such that the $\S$-to-$\W$ iteration embedding exists, the $\S$-to-$\W$ iteration is above $\k$ and $A$ is Wadge below $\Sigma^m_{\W^b}$. 
\end{lemma}
\begin{proof} The lemma follows from \rprop{fully backgrounded constructions}. Indeed, let $w$ be a window of $\P$ such that $g*m$ is generic for a poset in $\P|\nu^w$. Let $\N_0$ be the output of $\textsf{FBIC}(\S, \Sigma^m, \d^w, \nu^w, \k)$ (see \rthm{success of fbic}).  Thus, $\N_0$ is a $\Sigma^m$-iterate of $\S$ above $\k$, and all of its extenders with critical point $>\k$ have, in $\P[g*h]$, full background certificates whose critical points are strictly greater than $\nu^w$. We also have that\\\\
(A) $Ord\cap \N_0=\d^w$ (see \rprop{fully backgrounded constructions}).\\\\ Working inside $\N_0$, let $\N$ be the output of the hod pair construction of $\N_0$ done using extenders with critical point $>\d^{\N_0^b}$.

It follows from \rlem{lem:Gamma} that there is a countable iteration $p$ of $\P_0$ according to $\Sigma^m$ such that $\pi^{p, b}$ exists and letting $\R=\pi^{p, b}(\P_0^b)$, $A$ is Wadge below $\Sigma^m_{\R}$. Fix such a $(p, \R)$. We now claim that\\

\textit{Claim.} for some $\xi<\d^{\N^b}$, $\N|\xi$ is a $\Sigma^m_\R$-iterate of $\R$.
\begin{proof} To see this, we compare $\R$ with the construction producing $\N$. We need to see that $\R$ can be compared with $\N$. There are two ways such a comparison could go wrong.
\begin{enumerate}
\item $\N$ and $\R$ are not full with respect to the same $Lp$-operator. More precisely, for some normal $\Sigma^m_\R$-iteration $\T$ of limit length letting $b=\Sigma^m_\R(\T)$, either 
\begin{enumerate}
\item $\M^\T_b\models ``\d(\T)$ is a Woodin cardinal" and $\N\models ``\d(\T)$ is not a Woodin cardinal" or
\item  $\M^\T_b\models ``\d(\T)$ is not a Woodin cardinal" and $\N\models ``\d(\T)$ is a Woodin cardinal".
\end{enumerate}
\item A strategy disagreement is reached. More precisely, for some normal $\Sigma^m_\R$-iteration $\T$ with last model $\R^*$ and some $\xi$ which is a  Woodin cardinal of $\R^*$, $\R^*|\xi=\N|\xi$ yet $S^\N_{\R^*|\xi}\not= \Sigma^m_{\R^*|\xi}$.
\end{enumerate} 
It is easier to argue that case 2 cannot happen. This essentially follows from \cite[Theorem 4.13.2]{hod_mice_LSA}. Because $\N$ is backgrounded via extenders whose critical points are $>\d^{\N_0^b}$, the fragment of $\Sigma^m_{\N_0}$ we need to compute the strategy of $\N|\xi$ is the fragment that acts on non-dropping trees that are above $\d^{\N_0^b}$ and are based on $\N_0|\zeta$. Then \cite[Theorem 4.13.2]{hod_mice_LSA} implies that this fragment of $\Sigma^m_{\N_0}$ is induced by the unique strategy of $\P|\zeta$. The same strategy of $\P|\zeta$ also induces $\Sigma^m_{\R^*|\xi}$. Therefore, clause 2 cannot happen.

We now show that clause 1 also cannot happen. Suppose $\xi<\d^{\N^b}$ is a limit of Woodin cardinals or is a Woodin cardinal. Let $\zeta=o^\N(\xi)$, the Mitchell order of $\xi$, and let  $\T$ be a normal tree on $\R$ with last model $\W$ such that $\W|\zeta=\N|\zeta$ and the generators of $\T$ are contained in $\zeta$. Furthermore assume that $\zeta$ is a cutpoint in $\W$. Let $\nu$ be the least Woodin cardinal of $\W$ above $\zeta$, and let $\tau$ be the least Woodin cardinal of $\N$ above $\zeta$. It is enough to show that whenever $(\T, \W, \xi, \zeta, \nu, \tau)$ are as above then $\W|\nu$ normally iterates via $\Sigma^m_{\W|\nu}$ to $\N|\tau$. 

To see this it is enough to show that if $\U$ is a normal tree on $\W|\nu$ of limit length and $\m(\U)\insegeq \N|\tau$ then setting $b=\Sigma^m_{\W}(\U)$, either
\begin{enumerate}
\item $\d(\U)<\tau$ and $\Q(b, \U)$ exists and $\Q(b, \U)\insegeq \N|\tau$ or
\item $\d(\U)=\tau$ and $\pi^\U_b(\nu)=\tau$.
\end{enumerate}
To see the above, fix $\U$ and $b$ as above. Suppose first that $\d(\U)<\tau$. Let $\Q\insegeq \N|\tau$ be largest such that $\Q\models ``\d(\U)$ is a Woodin cardinal". Then, as $\Sigma^m$ is fullness preserving, $\Q\insegeq \M^\U_b$.

 Suppose then $\d(\U)=\tau$. If $\pi^\U_b(\nu)>\tau$ then $\Q(b, \U)$-exists and is $Ord$-iterable inside $\P[g*m]$. Working inside $\N$, let $\K$ be the output of the fully backgrounded construction of $\N$ done with respect to $S^\N_{\N|\tau}$ over $\N|\tau$ and using extenders with critical point $>\d^{\N^b}$. Because $\K$ is universal we must have that $\Q(\b, \U)\insegeq \K$. Thus, $\K\models ``\tau$ is not a Woodin cardinal", which implies that $\N\models ``\tau$ is not a Woodin cardinal". 
 \end{proof}
 
Let now $\Y^*$ be a normal tree on $\S$ according to $\Sigma^m_\S$ whose last model is $\N_0$. Let $\eta\in (\d^{\N_0^b}, \d^w)$ be such that $\R$ iterates to the hod pair construction of $\N_0|\eta$. Let $E\in \vec{E}^{\N_0}$ be such that $\cp(E)=\d^{\N_0^b}$ and $lh(E)>\eta$ (the existence of such an $E$ follows from (A) above). Let $\a<lh(\Y^*)$ be the least such that $E\in \M_\a^{\Y^*}$ and set $\Y^{**}=\Y^*\rest \a+1$. Finally, set $\Y=\Y^{**\frown} \{E\}$. Notice that if $\V$ is the last model of $\Y$ then $\pi^\Y$-exists.

 To finish the proof of the lemma, we need to take a countable Skolem hull of $\P|\l[g*m]$ where $\l=((\d^w)^+)^\P$.  Let $\pi:\M\rightarrow \P|\l[g*m]$ be a countable Skolem hull of $\P|\l[g*m]$ such that $\R, \N, \Y \in rng(\pi)$. Let $\X=\pi^{-1}(\Y)$ and let $\W$ be the last model of $\X$. By elementarity, $\X=\X^{*\frown}\{ F\}$ and $\R$ normally iterates via $\Sigma^m_\R$ to a hod pair construction of $\W|lh(F)$. It follows now that $\Sigma^m_{\W^b}$ is Wadge above $\Sigma^m_\R$, and hence, Wadge above $A$. Therefore, $\X$ is as desired. 
\end{proof}

\begin{lemma}\label{elem embedding}
There is an elementary embedding 
\begin{center}
$j: L(\Delta^h,\mathbb{R}^{\P[g*h]})\rightarrow L(\Delta^{h*k},\mathbb{R}^{\P[g*h*k]})$
\end{center}
 such that for each $A\in \Delta^h$, $j(A)=A^k$, the interpretation of $A$ in $\P[g*h*k]$. 
 \end{lemma}
\begin{proof}
Let $W_1=L(\Delta^h,\mathbb{R}^{\P[g*h]})$ and $W_2=L(\Delta^{h*k},\mathbb{R}^{\P[g*h*k]})$. Let $C$ be the set of inaccessible cardinals of $\P[g*h*k]$. Because we have a class of Woodin cardinals, it follows that $(\Delta^h)^\#$ exists. Moreover, for $\Gamma\subseteq \powerset(\bR)$, assuming $\Gamma^\#$ exists, any set in $L(\Gamma,\bR)$ is definable from a set in $\Gamma$, a real and a finite sequence of indiscernibles. It is then enough to show that\\\\
(*) if $s=(\a_0, ..., \a_n)\in C^{<\omega}$, $A\in \Delta^h$, $x\in \bR_{g*h}$ and $\phi$ is a formula then \begin{center}
$L(\Delta^h, \bR_{g*h})\models \phi[A, x, s]$ if and only if $L(\Delta^{h*k}, \bR_{g*h*k})\models \phi[A^k, x, s]$.
\end{center}
Indeed, we first show that (*) induces an elementary $j:W_1\rightarrow W_2$ as desired.  Let $Y$ be the set of $a$ that are definable over $L(\Delta^{h*k}, \bR_{g*h*k})$ from a member of $C^{<\omega}$, a set of the form $A^k$ for some $A\in \Delta^h$ and a real $x\in \bR_{g*h}$. Notice that (*) implies that\\

\textit{Claim 1.} $Y$ is elementary in $L(\Delta^{h*k}, \bR_{g*h*k})$.
\begin{proof}
 We show that $Y$ is $\Sigma_1$-elementary. The general case follows from Tarski-Vaught criteria. To see this fix $a\in Y$ and let $\phi$ be a $\Sigma_1$ formula. Suppose that 
 \begin{center}
 $W_2\models \phi[a]$.
 \end{center}
Fix a term $t$, $s\in C^{<\omega}$, a set of the form $A^k$ where $A\in \Delta^h$ and $x\in \bR_{g*h}$ such that $a=t^{W_2}[s, A^k, x]$. It then follows from (*) that if $b=t^{W_1}[s, A, x]$ , $W_1\models \phi[b]$. Let $\phi=\exists u \psi(u, v)$. Fix a term $t_1$, $s_1\in C^{<\omega}$, $B\in \Delta^h$ and $y\in \bR_{g*h}$ such that setting $c=t^{W_1}_1[s_1, B, y]$, $W_1\models \psi[c, b]$. Therefore, (*) implies that if $d=t^{W_2}[s_1, B^k, y]$ then $W_2\models \psi[d, a]$. As $d\in Y$, we have $Y\models \phi[a]$. 
 \end{proof}
 
 Let now $N$ be the transitive collapse of $Y$. It is enough to show that $N=W_1$. This easily follows from (*) and the proof of the claim. For example let us show that $\bR^N=\bR^{W_1}$. Fix $x\in \bR^N$. Let $t$ be a term, $s\in C^{<\omega}$, $A\in \Delta^h$ and $a\in \bR^{W_1}$ such that $x=t^{W_2}[s, A^k, a]$. Letting $y=t^{W_1}[s, A, a]$, it is easy to see that $x=y$. We now let $j:W_1\rightarrow W_2$ be the inverse of the transitive collapse of $Y$. Clearly $j$ is elementary and $j(A)=A^k$ for $A\in \Delta^h$. 
 
 By a similar reduction, using the definition of $\Delta^h$ and $\Delta^{h*k}$ it is enough to show that (**) holds where\\\\
(**) if $s=(\a_0, ..., \a_n)\in C^{<\omega}$, $\T$ is a countable iteration of $\P_0$ according to $\Sigma^h$ such that $\pi^{\T, b}$ exists, $\R=\pi^{\T, b}(\P_0^b)$, $x\in \bR_{g*h}$ and $\phi$ is a formula then \begin{center}
$W_1\models \phi[\Sigma^h_\R, x, s]$ if and only if $W_2\models \phi[\Sigma^{h*k}_\R, x, s]$.
\end{center}

To show (**), let $s\in C^{<\omega}$, $\T$ be a countable iteration of $\P_0$ according to $\Sigma^h$ such that $\pi^{\T, b}$ exists, $\R=\pi^{\T, b}(\P_0^b)$, $x\in \bR_{g*h}$ and $\phi$ be a formula such that $W_1\models \phi[\Sigma^h_\R, x, s]$. Notice that without losing generality, we can assume that $\pi^\T$ exists, as otherwise we can work with a shorter initial segment of $\T$ that produces the same bottom part $\R$. Let $\S^*$ be the last model of $\T$, and let $\S^{**}$ be the ultrapower of $\S^*$ by the least extender on the sequence of $\S^*$ with the critical point $\d^{\R}$. Let $\iota$ be the least Woodin of $\S^{**}$ that is $>\d^\R$.  Let $\W$ be the $\Sigma^h$-iterate of $\S^{**}$ that is obtained via an $x$-genericity iteration done in the window $(\d^\R, \iota)$. 

We would now like to see that $W_2\models  \phi[\Sigma^{h*k}_\R, x, s]$. The idea is to realize $W_1$ and $W_2$ respectively  as a derived model of $\W$. Given a transitive model of set theory $M$ with $\l$ a limit of Woodin cardinals of $M$ we let $D(M, \l)$ be the derived model at $\l$ as computed by some symmetric collapse of $\l$. While $D(M, \l)$ depends on this generic, its theory does not. Thus, expressions like $D(M, \l)\models \psi$ have an $\sf{uvs}$ meaning. If $u\subseteq Coll(\omega, <\l)$ is the generic then $D(\M, \l, u)$ is the derived model computed using $u$.\footnote{This is sometimes called the ``old" derived model. $D(\M,\lambda,u)$ has the form $L(\mathbb{R}^*_u, Hom^*_u)$ where $\mathbb{R}^*_u = \bigcup_{\alpha<\lambda} \mathbb{R}^{\M[u\rest \alpha]}$ and $Hom^*_u$ is the collection of $A\subseteq \mathbb{R}^*_u$ in $\M(\mathbb{R}^*_u)$ such that there are $< \lambda$-complementing trees $T,U\in \M[u\rest \beta]$ for some $\beta<\lambda$ such that $p[T]^{\M(\mathbb{R}^*_u)} = A = \mathbb{R}^*_u - p[U]$.} 

To finish the proof we will need a way of realizing $W_1$ as a derived model of an iterate of $\W$ that is obtained by iterating above $\xi$, where $\xi$ is the least Woodin cardinal of $\W$ above $\delta^\R$. The same construction will also realize $W_2$ as a derived model of a $\W$'s iterate. Let $l\subseteq Coll(\omega, \Gamma^\infty_{g*h})$ be $\P[g*h]$-generic. Working in $\P[g*h*l]$,  let $(A_i: i<\omega)$ be a generic enumeration of $\Gamma^\infty_{g*h}$, and let $(x_i: i<\omega)$ be a generic enumeration of $\bR_{g*h}$.\\\\
(1) There is sequence $(\W_i, p_i^*, p_i: i<\omega) \in \P[g*h*l]$ such that
\begin{enumerate}
\item for each $n<\omega$, $(\W_i, p_i^*, p_i: i\leq n)\in HC^{\P[g*h]}$,
\item $\W_0=\W$,
\item letting $E_i\in \vec{E}^{\W_i}$ be the Mitchell order 0 measure on $\d^{\W_i^b}$ and $\M_i=Ult(\W_i, E_i)$\footnote{We take the ultrapower by $E_i$ to have more cutpoint Woodin cardinals.}, $p_i^*$ is an iteration of $\M_i$ according to $\Sigma^h_{\M_i}$ that is above $\d^{\W_i^b}$, has a last model $\N_i$ and is such that $\pi^{p_i^*}$ exists and for some $\nu_i<\d^{\N_i^b}$ a Woodin cardinal of $\N_{i}$, $A_i<_w \Sigma^h_{\N_i|\nu_i}$,
\item fixing some $\nu_i$ as above and letting $\xi$ be the least Woodin cardinal of $\N_i$ that is $>\nu_i$, $\W_{i+1}$ is the $\Sigma^h_{\N_i}$-iterate of $\N_i$ that is above $\nu_i$ and makes $x_i$ generic at the image of $\xi$; $p_i$ is the corresponding iteration.
\end{enumerate}
The proof of (1) is a straightforward application of \rlem{technical lemma}. Let $\pi_{i, j}:\W_i\rightarrow \W_j$ be the iteration embedding, and let $\W_\omega$ be the direct limit of $(\W_i, \pi_{i, j}: i<j<\omega)$. It follows that for some $u\subseteq Coll(\omega, <\d^{\W_\omega^b})$-generic, \\\\
(2) $\bR^{\W_\omega[u]}=\bR_{g*h}$ and $\Gamma^\infty_{g*h}=\powerset(\bR_{g*h})\cap D(\W_\omega, \d^{\W_\omega^b}, u)$ and hence,\\\\
(3) $W_1=D(\W_\omega, \d^{\W_\omega^b}, u)$\\\\
Letting $S$ stand for the strategy predicate and $t$ be the sequence of the first $n$ indiscernibles of $\W|\d^{\W^b}$, we thus get by our assumption $W_1\models \phi[\Sigma^h_\R, x, s]$ and by elementarity that\\\\
(4) $D(\W[x], \d^{\S^b})\models \phi[S^\W_\R, x, t]$.\\\\
The same construction that gives (3) also gives $\N_\omega$ and $v$ such that\\\\
(5) $\N_\omega$ is a $\Sigma^{h*k}$-iterate of $\W$ above $\xi$, $v\subseteq Coll(\omega, <\d^{\N_\omega^b})$ is generic and $D(\N_{\omega}, \d^{\N_\omega^b}, v)=W_2$.\\\\
Thus, $W_2\models \phi[\Sigma^{h*k}_\R, x, t_1]$ where $t_1$ is the image of $t$ in $\N_\omega$. By indiscernability we get that  $W_2\models \phi[\Sigma^{h*k}_\R, x, s]$. 
\end{proof}

\subsection{An upper bound for $\sf{LSA-over-uB}$}\label{sec:upperbound2}

Let $(\P_0,\Sigma^-), \P, \Sigma, g$ be as before (see right after \rthm{thm:upper_bd}). Now we show $\sf{LSA-over-uB}$ is satisfied in $\P[g]$. Fix a poset $\mathbb{P}\in P[g]$ and let $h\subseteq \mathbb{P}$ be $\P[g]$-generic. We will show that
\begin{enumerate}
\item $L(\Sigma^{g*h}, \bR_{g*h})\models \sf{LSA}$ and
\item $\Gamma^{\infty}_{g*h}$ is the Suslin co-Suslin sets of $L(\Sigma^{g*h}, \bR_{g*h})$.
\end{enumerate}
Clause 2 above is an immediate consequence of Clause 1, and the results of the previous section.

We now show clause 1. Let $(\gamma_i : i<\omega)$ be the first $\omega$ Woodin cardinals of $\P[g*h]$ and $\gamma = sup_{i<\omega} \gamma_i$. Let $w_i$ be the corresponding consecutive windows determined by the $\gamma_i$'s. Write $\Lambda$ for $\Sigma^{h}$, the canonical interpretation of $\Sigma$ in $P[g*h]$. In $\P[g*h]$, let 
\begin{center}
$\pi: \M \rightarrow (\P[g*h]|\gamma^+)^\#$,
\end{center}
be elementary and such that $\M$ is countable and $\cp(\pi) > \d_0$. For each $i$, let $\delta_i = \pi^{-1}(\gamma_i)$, and $\lambda = \sup_{i<\omega} \delta_i$. Note that, because $\cp(\pi) > \d_0$, $\M|\l$ is closed under $\Lambda$, and $\lambda$ is the supremum of the Woodin cardinals of $\M$.  It follows from \rprop{citerability in excellent hybrid mice} that $\M$ has a $\nu^{w_0}$-strategy acting on non-dropping trees based on the interval $[\pi^{-1}(\nu^{w_0}), \lambda)$ in $\P[g*h]$; call this strategy $\Psi$.

Let $k\subseteq Coll(\omega, <\lambda)$ be $\M$-generic. Let $\mathbb{R}^*_k = \bigcup_{\xi<\lambda} \mathbb{R}^{\M[k|\xi]}$ and recall the "new" derived model of $\M$ at $\lambda$
\begin{center}
$D^+(\M,\lambda,k) = L(\{A\in \powerset(\mathbb{R}^*_k)\cap \M(\mathbb{R}^*_k): L(A,\mathbb{R}^*_k) \models \sf{AD}^+  \})$. 
\end{center}
By Woodin's derived model theorem, cf. \cite{DMT}, $D^+(\M,\l,k)\models \sf{AD}^+$. Again, the theory of $D^+(\M,\lambda,k)$ does not depend on $k$. When we reason about the theory of the new derived model without concerning about any particular generic, we write $D^+(\M,\lambda)$. Recall that we set $\Lambda=\Sigma^h$.
\begin{proposition}\label{claim:inDM}
$\Lambda \cap \M(\mathbb{R}^*_k)\in D^+(\M, \l)$. Furthermore, in $D^+(\M,\l)$, $L(\Lambda,\mathbb{R})\vDash \sf{LSA}$.
\end{proposition}
\begin{proof}
First note that there is a term $\tau\in \M$ such that $(\M,\Psi,\tau)$ term captures $\Sigma^{h}$. More precisely, letting $i:\M \rightarrow \N$ be an iteration map according to $\Psi$, let $l$ be a $<i(\lambda)$-generic over $\N$, then $Code(\Sigma^{h})\cap \N[l] = i(\tau)_{l}$; this follows from results in Section 2 (cf. \rprop{generic interpretability in ehpm}). To see that in $\M(\mathbb{R}^*_k)$,  $L(\Lambda,\mathbb{R}) \models \sf{AD}$, suppose not. Let $x$ be a real and $A$ be the least $OD(\Lambda,x)$ counterexample to $\sf{AD}$ in  $L(\Lambda,\mathbb{R})$. Also, by minimizing the ordinal parameters, we may assume $A$ is definable from $x$ and $\Lambda$ in $L(\Lambda,\mathbb{R})$. Using the term $\tau$ for $\Lambda$, we can easily define a term $\sigma$ over $\M[x]$ such that $(\M[x],\Psi,\sigma)$ term captures $A$.\footnote{Note that there is a generic $k'\subseteq Coll(\omega,<\lambda)$ for $\M[x]$ such that $\mathbb{R}^*_k = \mathbb{R}^*_{k'}$.} Applying Neeman's theorem (cf. \cite{neeman1995optimal}), we get that $A$ is determined. Finally, let $i:\M[x]\rightarrow \N$ be a $\mathbb{R}^{\P[g*h]}$-genericity iteration according to $\Psi$.\footnote{In $\P[g*h]$, let $k\subseteq Coll(\omega,\mathbb{R})$ be generic and let $(x_i : i<\omega)$ be the generic enumeration of the reals. The iteration $i$ is the direct limit of the system $(\M_m, i_{m,m+1}: m<\omega)$ where $\M_0 = \M[x]$, for each $m$, $i_{m,m+1}:\M_m \rightarrow \M_{m+1}$ is the $x_m$-genericity iteration that makes $x_m$ generic at the image of $\delta_m$.} By the argument just given, in $\N(\mathbb{R}_{g*h})$, $A$ is determined. So $\M(\mathbb{R}^*_k)\models \sigma_k$ is determined. Contradiction. 

If $\sf{LSA}$ fails in $L(\Lambda,\mathbb{R})$, then $\Lambda$ is Suslin co-Suslin in $D^+(\M,\lambda)$, and the argument in \rprop{claim:not_uB} gives a contradiction. The point is that in $D^+(\M,\l)$, the Wadge ordinal of $\Gamma^b(\P_0,\Lambda)$ is a limit of Suslin cardinals, and the failure of $\sf{LSA}$ means that there is a larger Suslin cardinal above the Wadge ordinal\footnote{The supremum of Wadge ranks of the sets of reals in $\Gamma^b(\P_0,\Lambda)$.}  of $\Gamma^b(\P_0,\Lambda)$. So $\Lambda$ is Suslin co-Suslin in $D^+(\M,\lambda)$. Now we can run the argument in \rprop{claim:not_uB} to obtain a contradiction. Hence in $D^+(\M,\lambda)$, 
\begin{equation}\label{eqn:LSA}
L(\Lambda,\mathbb{R}) \vDash \sf{LSA}. 
\end{equation}
\end{proof}

Now perform a $\mathbb{R}^{\P[g*h]}$-genericity iteration according to $\Psi$ at $\lambda$, more precisely, there is an iteration $i:\M\rightarrow \N$ according to $\Psi$ such that letting $l \subseteq Coll(\omega, < i(\lambda))$ be $\N$-generic, letting $\mathbb{R}^*_{l} = \bigcup_{\xi<i(\lambda)}\mathbb{R}^{\N[l|\xi]}$, we get 
\begin{center}
$\mathbb{R}^{\P[g*h]} = \mathbb{R}^*_{l}$
\end{center}
and  
\begin{center}
$L(\Sigma^{h},\mathbb{R}^{\P[g*h]}) \subseteq D^+(\N,i(\lambda),l)$, 
\end{center}
hence by (\ref{eqn:LSA}),
\begin{center}
$L(\Sigma^{h},\mathbb{R}^{\P[g*h]}) \models \sf{LSA}$.
\end{center}
This completes the proof of clause 1 above and also the proof of $\sf{LSA-over-uB}$ in $\P[g]$.

\subsection{An upper bound for $\sf{Tower \ Sealing}$}

Let $(\P_0,\Sigma^-), \P, \Sigma, g$ be as before (see right after \rthm{thm:upper_bd}). We prove clause (2) of $\sf{Tower \ Sealing}$ holds in $\P[g]$. Clause (1) has already been established by the previous sections. Let $\mathbb{P}\in \P[g]$ be any poset, $h\subseteq \mathbb{P}$ be a $\P[g]$-generic and let $\delta$ be Woodin in $\P[g*h] =_{def} W$. Let $G\subseteq \mathbb{Q}_{<\delta}$ be $W$-generic (the argument for $\mathbb{P}_{<\delta}$ is the same) and $j: W\rightarrow M \subset W[G]$ be the generic elementary embedding induced by $G$\footnote{$\mathbb{Q}_{<\delta}$ and $\mathbb{P}_{<\delta}$ are the countable and full stationary tower forcings.}. 


Let $\Lambda$ be the canonical interpretation of $\Sigma^-$ in $W$ and $\Lambda^G$ be the canonical interpretation of $\Lambda$ in $W[G]$ (considered as the short-tree strategy of $\P_0$ acting on countable trees). Now, by the fact that $M$ is closed under countable sequences in $W[G]$ and the way $\Lambda^G$ is defined (using generic interpretability)
\begin{equation}\label{eqn:equal}
\Lambda^G = j(\Lambda).
\end{equation}
Here is the outline of the argument. Let $\T$ be countable and according to both $j(\Lambda)$ and $\Lambda^G$. Note that $\T\in M$. Suppose $\T$ is $\sf{nuvs}$ (the case $\T$ is $\sf{uvs}$ is similar). One gets that in $W[G]$, $\Q(\Lambda^G(\T),\T)$ exists and is authenticated by $\vec{C}$, a fully backgrounded authenticated construction in $W$ where extenders have critical point $> \delta$; note that we can take $\vec{C}\in W$. This implies that $\Q(\Lambda^G(\T),\T)$ is authenticated by $j(\vec{C})\in M$; $\Q(j(\Lambda)(\T),\T)$ is also authenticated by $j(\vec{C})$ in $M$. The details are very similar to the proof of Proposition \ref{generic interpretability in ehpm}. So $j(\Lambda)(\T) = \Lambda^G(\T)$. Hence $\Lambda^G\in M$.

By Lemma \ref{lem:Gamma}, 
\begin{center}
$(\Gamma^\infty)^{W[G]} = \Gamma^b(\P_0,\Lambda^G)$. 
\end{center}
By elementarity, the fact that $(\Gamma^\infty)^W = \Gamma^b(\P_0,\Lambda)$, \ref{eqn:equal}, and Lemma \ref{lem:Gamma}, 
\begin{center}
$j((\Gamma^\infty)^W) = \Gamma^b(\P_0,\Lambda^G)$. 
\end{center}
So indeed, $(\Gamma^\infty)^{W[G]}= j((\Gamma^\infty)^W)$ as desired.

\begin{remark}
Another proof of clause (2) of $\sf{Tower \ Sealing}$ is the following. We give a sketch: by results of \cite{hod_mice_LSA}, $\Lambda^G\in L(\pi^b_{\P_0,\infty}[\P_0], \mathcal{H},(\Gamma^\infty)^{W[G]})$. Here working in $W[G]$, let $\mathcal{H}^-$ be the direct limit of hod pairs $(\R,\Delta)\in L(\Gamma^\infty)$ such that in $ L(\Gamma^\infty)$, $\Delta$ is fullness preserving, has hull and branch condensation. $|\mathcal{H}^-| = V^{HOD}_\Theta$ in $L(\Gamma^\infty)$. Let $\mathcal{H} = \bigcup\{\M : \mathcal{H}^- \lhd \M, \M \textrm{ is a sound, hybrid } \textrm{countably iterable premouse such that } \rho_\omega(\M)\leq o(\mathcal{H}^-)\}$. For each $\M\lhd \mathcal{H}$ as above, for every $\N$ countably transitive such that $\N$ is embeddable into $\M$, $\N$ has an $\omega_1$-strategy in $L(\Gamma^\infty)$.

If $(\Gamma^\infty)^{W[G]}\neq j((\Gamma^\infty)^W)$, then suppose the former is a strict Wadge initial segment of the latter (the other case is handled similarly). So the model 
\begin{center}
$L(\pi^b_{\P_0,\infty}[\P_0], \mathcal{H},(\Gamma^\infty)^{W[G]})\in M$
\end{center}
as $M$ is closed under $\omega$-sequences in $W[G]$. In fact, we get that $\Lambda^G \in j((\Gamma^\infty)^W)$. By Generation of Mouse Full Pointclasses (applied in $L(j((\Gamma^\infty)^W))$ and a comparison argument as in Lemma \ref{claim:not_uB}, there is a (maximal) $\Lambda^G$-iterate $\S$ of $\P_0$ such that $\S$ has an iteration strategy $\Psi$ such that
\begin{itemize}
\item $\Psi^{stc} = (\Lambda^G)_\S$,
\item $\Gamma^b(\S,\Psi) = \Gamma(\P_0,\Lambda^G)$,
\item $\Psi \in j((\Gamma^\infty)^W)$.
\end{itemize}
By elementarity, the existence of $\S,\Psi$ holds in $L((\Gamma^\infty)^W)$. This contradicts the fact that $j(\Lambda)\notin (\Gamma^\infty)^W$.
\end{remark}

\section{Basic core model induction}\label{sec:uBstacks}

The notation introduced in the section will be used throughout this paper. It will be wise to refer back to this section for clarifications. From this point on the paper is devoted to proving that both $\sf{Sealing}$ and $\sf{LSA-over-uB}$ imply the existence of a (possibly class size) excellent hybrid premouse. As we have already shown that a forcing extension of an excellent hybrid premouse satisfies both $\sf{Sealing}$ and $\sf{LSA-over-uB}$, this will complete the proof of \rthm{thm:main_theorem}.

We will accomplish our goal by considering $\H$ of $L(\Gamma^\infty, \bR)$ and showing that, in some sense, it reaches an excellent hybrid premouse. Our first step towards this goal is to show that $\Theta$ is a limit point of the Solovay sequence of $L(\Gamma^\infty, \bR)$. 

\begin{proposition}\label{adr in gamma ub} Assume there are unboundedly many Woodin cardinals. Furthermore, assume either $\sf{Sealing}$, or $\sf{Tower \ Sealing}$, or $\sf{LSA-over-uB}$. Then for all set generic $g$, the following holds in $V[g]$. 
\begin{enumerate}
\item $\powerset(\bR)\cap L(\Gamma^\infty, \bR)=\Gamma^\infty$.
 \item $L(\Gamma^\infty, \bR)\models \sf{AD}_{\mathbb{R}}$.
 \end{enumerate}
\end{proposition}
\begin{proof} Towards a contradiction assume that $L(\Gamma^\infty, \bR)\models \neg \sf{AD_{\mathbb{R}}}$. By a result of Steel (\cite[Theorem 4.3]{DMT}), every set in $\Gamma^\infty$ has a scale in $\Gamma^\infty$. Notice then that clause 1 implies clause 2. This is because given clause 1, $L(\Gamma^\infty, \bR)$ satisfies that every set has a scale, and therefore, it satisfies $\sf{AD_{\mathbb{R}}}$\footnote{Recall that Martin and Woodin showed that under ${\sf{AD}}$, $\sf{AD_{\mathbb{R}}}$ is equivalent to the statement that every set of reals has a scale. See \cite{MartinWoodin}. Also, by results of Martin, Steel and Woodin, assuming class of Woodin cardinals, every uB set is determined. See \cite{DMT}.}.

It is then enough to show that clause 1 holds. It trivially follows from $\sf{Sealing}$ or $\sf{Tower \ Sealing}$. To see that it also follows from $\sf{LSA-over-uB}$, fix a set $A\subseteq \bR$ such that $\Gamma^\infty$ is the set of Suslin, co-Suslin sets of $L(A, \bR)$ and $L(A, \bR)\models \sf{LSA}$. It now follows that if $\k$ is the largest Suslin cardinal of $L(A, \bR)$ then, in $L(A, \bR)$, $\Gamma^\infty$ is the set of reals whose Wadge rank is $<\k$. Since $\kappa$ is on the Solovay sequence of $L(A,\mathbb{R})$, $\Gamma^\infty = L(\Gamma^\infty)\cap \powerset(\mathbb{R})$. Therefore clause 1 follows. 
\end{proof}

For the rest of this paper we write $\Gamma^\infty\models_{\Omega} \sf{AD}_{\bR}$ to mean that clause 1 and 2 above hold in all generic extensions. $\Omega$ here is a reference to Woodin's $\Omega$-logic. We develop the notations below under  $\Gamma^\infty\models_{\Omega} \sf{AD}_\bR$.

Suppose $\mu$ is a cardinal. Let $g\subset Col(\omega,<\mu)$ be $V$-generic. Working in $V$, we say that a pair $(\M,\Sigma)$ is a \textit{hod pair} at $\mu$ if 
\begin{enumerate}
\item $\M\in V_\mu$, 
\item $\Sigma$ is a  $(\mu,\mu)$-iteration strategy of $\M$ that is in $\Gamma^\infty$ in $V^{Coll(\omega,|\M|)}$ and is positional, commuting and has branch condensation, and
\item $\Sigma$ is fullness preserving with respect to mice with $\Gamma^\infty$-iteration strategy.
\end{enumerate}

Let $\mathcal{F}$ be the set of hod pairs at $\mu$. It is shown in \cite{hod_mice_LSA} that hod mice at $\mu$ can be compared (see \cite[Chapter 4.6 and 4.10]{hod_mice_LSA}.). More precisely, given any two hod pairs $(\M, \Sigma)$ and $(\N, \Lambda)$ in $\mathcal{F}$, there is a hod pair $(\S, \Psi)\in \mathcal{F}$ such that for some $\M^*\unlhd_{hod} \S$ and $\N^*\unlhd_{hod} \S$, 
\begin{enumerate}
\item $\M^*$ is a $\Sigma$-iterate of $\M$ such that the main branch of $\M$-to-$\M^*$ iteration doesn't drop,
\item $\N^*$ is a $\Lambda$-iterate of $\N$ such that the main branch of $\N$-to-$\N^*$ iteration doesn't drop,
\item $\Sigma_{\M^*}=\Psi_{\M^*}$ and $\Lambda_{\N^*}=\Psi_{\N^*}$ and
\item either $\S=\M^*$ or $\S=\N^*$.
\end{enumerate}

Working in $V[g]$, let $\mathcal{F}^+$ be the set of all hod pairs $(\M,\Sigma)$ such that $\M$ is countable  and $\Sigma$ is an $(\omega_1,\omega_1+1)$-strategy of $\M$ that is $\Gamma^\infty$-fullness preserving, positional, commuting, has branch condensation,\footnote{In \cite{hod_mice_LSA}, $\mathbb{P} = Col(\omega,\omega_2)$ but in our case, since $\mu$ is measurable, all results in \cite[Chapter 12]{hod_mice_LSA} hold in our context. The point is that we can work with stationaraily many hulls $X\prec H_\xi$ for some $\xi>>\Omega$ such that $X\cap \mu=\gamma$ is an inaccessible cardinal, $X^{<\gamma}\subseteq X$, and their corresponding uncollapse map $\pi_X: M_X\rightarrow H_\xi$. Or equivalently, we work with the ultrapower embedding $j_U: V \rightarrow Ult(V,U)$, noting that $j_U$ lifts to a generic elementary embedding on $V[G]$. By results in \cite{hod_mice_LSA}, $\Sigma$ has strong branch condensation and is strongly $\Delta$-fullness preserving.} and $\Sigma\rest \rm{HC}\in\Gamma^\infty$.

Because any two hod pairs in $\mathcal{F}^+$ can be compared, $\mathcal{F}$ \textit{covers} $\mathcal{F}^+$. More precisely, for each hod pair $(\M, \Sigma)\in \mathcal{F}^+$ there is $\Sigma$-iterate $\N$ of $\M$ such that the $\M$-to-$\N$ iteration doesn't drop on its main branch, $(\Sigma_\N\rest V)\in V$ and $\Sigma_\N$ is the unique extension of $(\Sigma_\N\rest V)$ to $V[g]$. 

Given any hod pair $(\M, \Sigma)$, let $I(\M, \Sigma)$ be the set of iterates $\N$ of $\M$ by $\Sigma$ such that the main branch of $\M$-to-$\N$ doesn't drop. Let $X\subseteq I(\M, \Sigma)$ be a \textit{directed} set, i.e., if $\N, \P\in  X$ then there is $\R\in X$ such that $\R$ is a $\Sigma_\N$-iterate of $\N$ and a $\Sigma_\P$-iterate of $\R$. We then let $\M_\infty(\M, \Sigma, X)$ be the direct limit of all iterates of $\M$ by $\Sigma$ that are in $X$. Usually $X$ will be clear from context and we will omit it.

Working in $V[g]$, let $\mathbb{R}_g = \mathbb{R}^{V[g]}$. Let $\mathcal{H}^{-}$ be the direct limit of hod pairs in $\mathcal{F}^+$. Because $\mathcal{F}$ covers $\mathcal{F}^+$, we also have that
$\mathcal{H}^-$ is the direct limit of hod pairs in $\mathcal{F}$.

Fix $(\M, \Sigma)\in \mathcal{F}^+$ such that $\M_\infty(\M, \Sigma )=_{def}\Q\unlhd_{hod} \mathcal{H}^-$. We let $\Psi_\Q=\Sigma^+_{\Q}$. $\Psi_\Q$ only depends on $\Q$ and does not depend on any particular choice of $(\M,\Sigma)\in \mathcal{F}^+$. Let $(\mH^-(\alpha): \alpha < \lambda)$ be the layers of $\mH^-$ (in the sense of \cite{hod_mice_LSA} and \cite{hod_mice}) and let $\Psi_\alpha$ be the strategy of $\mH^-(\alpha)$ for each $\alpha<\lambda$. $\Psi_\alpha$ is the tail strategy $\Sigma_\Q$ for $\Q = \M_\infty(\M,\Sigma)$ for any $(\M,\Sigma)\in \F^+$ such that $\M_\infty(\M,\Sigma)=\mH^-(\alpha)$. 
We now set
\begin{center}
$\Psi =_{def}\Psi_\mu=_{def} \oplus_{\alpha<\lambda^{\mathcal{H}^-}} \Psi_\alpha$.
\end{center}

\begin{definition}\label{cuB stack} Suppose $x$ is a set in $V(\bR_g)$ and $\Phi$ is an iteration strategy with hull condensation. Working in $V(\bR^*)$, let $Lp^{cuB, \Phi}(x)$ be the union of all sound $\Phi$-mice $\M$ over $x$ that project to $x$ and whenever $\pi: \N \rightarrow \M$ is elementary, $\N$ is countable, transitive then $\N$ has a universally Baire iteration strategy.
\end{definition}

Continuing, we set
\begin{enumerate}
\item $\mathcal{H} = Lp^{cuB, \Psi}(\mathcal{H}^-)$ (note that $\mathcal{H}\in V$),
\item $\Theta = o(\mathcal{H}^-)$,
\item $(\theta_\alpha : \alpha < \lambda)$ as the $\sf{Solovay\ Sequence}$ of $\Gamma^\infty$. Note that $\Theta = sup_\alpha \theta_\alpha$ and $\theta_\alpha = \delta^{\mH^-(\alpha)}$ for each $\alpha<\lambda$.
\end{enumerate}

We note that all objects defined in this section up to this point depend on $\mu$. To stress this, we will use $\mu$ as subscript. Thus, we will write, if needed, $\Psi_\mu$ or $\mH_\mu$ for $\Psi$ and $\mH$ respectively.  We will refer to the objects introduced above, e.g. $\mH_\mu$, $\Psi_\mu$ and etc, as the CMI objects at $\mu$. 

Given a hybrid strategy mouse $\Q$ and an iteration strategy $\Lambda$ for $\Q$, we say $\Lambda$ is \textit{potentially-universally Baire} if whenever $g\subseteq Coll(\omega, \Q)$ is generic there is a unique $\Phi\in V[g]$ such that
\begin{enumerate}
\item $\Phi\rest V=\Lambda$,
\item in $V[g]$, $\Phi$ is a uB iteration strategy for $\Q$.
\end{enumerate}
Similarly we can define potentially-$\eta$-uB iteration strategies. 
\begin{definition}\label{potentially ub stack} Suppose $\mu$ is a cardinal and $(\Q, \Lambda)$ is such that $\Q\in H_{\mu^+}$, $\Q$ is a hybrid strategy mouse and  $\Lambda$ is a potentially-uB strategy for $\Q$. Suppose $X\in H_{\mu^+}$. We then let $Lp^{puB, \Lambda}(X)$ be the union of all sound $\Lambda$-mice over $X$ that project to $X$ and have a potentially-uB iteration strategy. 
\end{definition}
Clearly $Lp^{puB, \Lambda}(X)\insegeq Lp^{cuB, \Lambda}(X)$. In many core model induction applications it is important to show that in fact $Lp^{puB, \Lambda}(X)= Lp^{cuB, \Lambda}(X)$. The reason this fact is important is that the first is the stack that we can prove is computed by the maximal model of determinacy containing $X$ after we collapse $X$ to be countable while if $\Q, X$ are already countable, the $OD(\Lambda, \Q, X)$ information inside the maximal model is captured by $Lp^{cuB, \Lambda}(X)$. This is because for countable $\Q, X$, $Lp^{cuB, \Lambda}(X)=Lp^{uB, \Lambda}(X)$ where the mice appearing in the latter stack have universally Baire strategies. The equality $Lp^{puB, \Lambda}(X)= Lp^{cuB, \Lambda}(X)$ is important for covering type arguments that appear in the proof of \rprop{lem:full_hull}.

\section{$Lp^{cuB}$ and $Lp^{puB}$ operators}

The following is the main result of this section, and it is the primary way we will translate strength from our hypothesis over to large cardinals. If $\mu$ is such that $Hom^*_g=\Gamma^\infty_g$ for any $g\subseteq Coll(\omega, <\mu)$, then we say that $\mu$ \textit{stabilizes uB}.

\begin{definition}\label{x captures} For each inaccessible cardinal $\mu$ let $A_\mu\subseteq \mu$ be a set that codes $V_\mu$.  We then say that $X\prec H_{\mu^+}$  \textbf{captures} $Lp^{cuB, \Psi_\mu}(A_\mu)$ if $Lp^{cuB, \Psi_\mu}(A_\mu)\in X$ and letting $\pi_X: M_X\rightarrow H_{\mu^+}$ be the uncollapse map and letting $\Lambda$ be the $\pi$-pullback of $\Psi_\mu$,
\begin{center} $\pi_X^{-1}(Lp^{cuB, \Psi_\mu}(A_\mu))=Lp^{cuB, \Lambda}(\pi_X^{-1}(A_\mu))$.\end{center}
\end{definition}

\begin{theorem}\label{main thm on stacks under sealing} Suppose there is a proper class of Woodin cardinals and a stationary class of measurable cardinals\footnote{And hence, a stationary class of measurable cardinals that are limit of Woodin cardinals.}. Suppose further that $\Gamma^\infty\models_{\Omega}\sf{AD}_{\mathbb{R}}$. There is then a stationary class $S$ of measurable cardinals that are limits of Woodin cardinals, a proper class $S_0\subseteq S$, and a regular cardinal $\nu\geq \omega_1$ such that the following holds:
\begin{enumerate}
\item for any $\mu\in S$, $\card{\mH_\mu}<\mu^+$, $\cf(Ord\cap \mH_\mu)<\mu$, and $\cf(Ord\cap Lp^{cuB,\Psi_\mu}(A_\mu))<\mu$;
\item for any $\mu\in S_0$,  $\mu$ stabilizes uB, $\cf(Ord\cap \mH_\mu)<\nu$, and $\cf(Ord\cap Lp^{cuB,\Psi_\mu}(A_\mu))<\nu$;
\item for any $\mu\in S_0$, there is $Y_\mu\in \powerset_\nu(H_{\mu^+})$ such that $A_\mu\in Y_\mu$ and whenever $X\prec H_{\mu^+}$ is of size $<\mu$, is $\nu$-closed and $Y_\mu\subseteq X$, $X$ captures $Lp^{cuB, \Psi_\mu}(A_\mu)$. 
\end{enumerate}
\end{theorem}

We emphasize that the arguments in this section (and in this paper) are carried out entirely in $\sf{ZFC}$, though it may appear that we are working with proper classes. See Remark \ref{rmk:regressive} for a more detailed discussion and summary. First, we prove a useful lemma, pointed out to us by Ralf Schindler. Below, by ``class", we of course mean ``definable class".

\begin{lemma}[$\sf{ZFC}$]\label{lem:pressing_down}
Suppose $S$ is a stationary class of ordinals. Suppose $f: S \rightarrow Ord$ is regressive, i.e. $f(\alpha)<\alpha$ for all $\alpha\in S$. There is an ordinal $\nu$ and a proper class $S_0\subseteq S$ such that $f[S_0] = \{\nu\}$.
\end{lemma}
\begin{proof}
Suppose not. For each $\nu$, let $\alpha_\nu = sup\{\alpha : f(\alpha) = \nu\}$ if $\nu \in rng(f)$ and $\nu+1$ otherwise. $\alpha_\nu$ exists because we are assuming that $f^{-1}(\{\nu\})$ is a set. Let $g: Ord \rightarrow Ord$ be the function: $\nu \mapsto \alpha_\nu$; hence $g(\nu)>\nu$ for all $\nu$. Let $C = \{\mu : g[\mu] \subseteq \mu \}$. So $C$ is a club class. Let $\alpha \in lim(C)\cap S$. We may assume for unboundedly many $\beta < \alpha$, $\beta \in rng(f)$. Then we easily get that $f(\alpha)$ is not $< \alpha$. Contradiction.
\end{proof}

Clause 1 of \rthm{main thm on stacks under sealing} follows easily from the above lemma.

\begin{proposition}\label{lem:dm_uB} Suppose there is a proper class of Woodin cardinals and $S$ is a stationary class of inaccessible cardinals that are limit of Woodin cardinals. Then there is a proper class $S^*\subseteq S$ such that whenever $\mu\in S^*$ and $g\subseteq Coll(\omega,<\mu)$ is $V$-generic, in $V[g]$, $Hom^*_g= \Gamma^\infty_g$.
\end{proposition}
\begin{proof}
Clearly $\Gamma^\infty_g\subseteq Hom^*_g$. Suppose then the claim is false. We then have a club $C$ such that whenever $\mu\in C\cap S$ and $g\subseteq Coll(\omega, <\mu)$, $\Gamma^\infty_g \not =Hom^*_g$. For each $\mu\in C\cap S$ let $\eta_\mu<\mu$ be least such that whenever $g\subseteq Coll(\omega, \eta_\mu)$, there are $\mu$-complementing trees $(T, U)\in V[g]$ with the property that $p[T]$ is not uB in $V[g][h]$ for any $V[g]$-generic $h\subseteq Coll(\omega, <\mu)$. By Lemma \ref{lem:pressing_down}, we then have a proper class $S_0\subseteq S$ such that for every $\mu_0<\mu_1\in S_0$, $\eta_{\mu_0}=\eta_{\mu_1}$. Let $\eta$ be this common value of $\eta_\mu$ for $\mu\in S_0$ and $g\subseteq Coll(\omega, \eta)$ be $V$-generic. For each $\mu\in S_0$ we have a pair $(T_\mu, U_\mu)\in V[g]$ that represents a $\mu$-uB set that is not uB. A simple counting argument then shows that for a proper class $S^{*}\subseteq S_0$ whenever $\mu_0, \mu_1\in S^{*}$, $V[g]\models p[T_{\mu_0}]=p[T_{\mu_1}]$. Letting $A=(p[T_\mu])^{V[g]}$ for some $\mu\in S^{*}$, we get a contradiction as $A$ is uB in $V[g]$.
\end{proof}

What follows is a sequence of propositions that collectively imply the remaining clauses of \rthm{main thm on stacks under sealing}. We start by establishing that the two stacks are almost the same. 

\begin{proposition}\label{cub=pub} Suppose $\mu, (\Q, \Lambda), X$ are as in \rdef{potentially ub stack} and suppose $\mu$ is in addition a measurable cardinal stabilizing uB. Let $j:V\rightarrow M$ be an embedding witnessing the measurability of $\mu$. Then $Lp^{cuB, \Lambda}(X) =(Lp^{puB, \Lambda}(X))^M$. 
\end{proposition}
\begin{proof} Let $j: V\rightarrow M$ be an embedding witnessing the measurability of $\mu$. Let $\M\insegeq Lp^{cuB, \Lambda}(X) $ be such that $\rho(\M)=X$. Let $h\subseteq Coll(\omega, <j(\mu))$ be generic. Consider $j(\M)$. In $M[h]$, $\M$, as it embeds into $j(\M)$, has a uB strategy. It follows that $\M$ has a potentially-uB strategy in $M$, and hence, $\M\insegeq (Lp^{puB, \Lambda}(X))^M$. Conversely, if $\M\insegeq (Lp^{puB, \Lambda}(X))^M$ is such that $\rho(\M)=X$ then in $M$, $\M$ has a potentially-uB strategy, and hence, in $V$, any countable $\pi:\M^*\rightarrow \M$ has a $\mu$-uB-strategy. As $\mu$-stabilizes uB, we have $\M\insegeq Lp^{cuB, \Lambda}(X)$.
\end{proof}


%
%

The next two propositions are rather important. Similar propositions hide behind any successful core model induction argument.

\begin{proposition}\label{pub stack small cof}  Suppose there are unboundedly many Woodin cardinals, $\mu$ is an inaccessible  cardinal and $\Gamma^\infty\models _\Omega \sf{AD}_{\bR}$. Suppose further that $\Lambda$ is a potentially-uB iteration strategy for some $\Q\in H_{\mu^+}$ and $X\in H_{\mu^+}$. Let $\M=Lp^{puB, \Lambda}(X)$. Then $\card{\M}<\mu^+$.
\end{proposition}  
\begin{proof}  Suppose $Ord\cap \M=\mu^+$. Let $g\subseteq Coll(\omega, \mu)$ be generic. Then 
\begin{center}
$(Lp^{puB, \Lambda}(X))^{V}=(Lp^{puB, \Lambda}(X))^{V[g]}$.
\end{center}
Moreover, $(Lp^{puB, \Lambda}(X))^{V[g]}\in L(\Gamma^\infty_g, \bR_g)$. Hence, $L(\Gamma^\infty_g, \bR_g)\models ``$there is an $\omega_1$-sequence of reals". This contradicts the fact that $L(\Gamma^\infty_g, \bR_g)\models \sf{AD^+}$. 
\end{proof}

\begin{corollary}\label{failure of lp covering} Suppose there are unboundedly many Woodin cardinals, $\mu$ is a measurable limit of Woodin cardinals that stabilizes uB and $\Gamma^\infty\models _\Omega \sf{AD}_{\bR}$. Suppose further that $\Lambda$ is a potentially-uB iteration strategy for some $\Q\in H_{\mu^+}$ and $X\in H_{\mu^+}$. Let $\M=Lp^{cuB, \Lambda}(X)$. Then $\card{\M}<\mu^+$ and $\cf(Ord\cap \M)<\mu$.
\end{corollary}
\begin{proof} Fix $j:V\rightarrow M$ witnessing the measurability of $\mu$. It follows from \rprop{cub=pub} that $\M=(Lp^{puB, \Lambda}(X))^M$. Applying \rprop{pub stack small cof} in $M$ we get that $\card{\M}<\mu^+$.

Assume next that $\cf(Ord\cap \M)=\mu$. Let $\eta=Ord\cap \M$ and let $\vec{C}$ be the $\square(\eta)$-sequence of $\M$. Because $\mu$ is measurable, we have that $\vec{C}$ is threadable. To see there is a thread $D$, note that $sup \ j[\eta] =_{def} \gamma < j(\eta)$. Let $E = j(\vec{C})_\gamma$ and $D = j^{-1}[E]$. Then $D$ is a thread through $\vec{C}$.

This implies that there is a $\Lambda$-mouse $\N$ extending $\M$ such that $\rho(\N)=\eta$ and every $<\mu$-submodel of $\N$ embeds into some $\N^*\insegeq \M$\footnote{This is a consequence of the proof of $\square$. $\N$ is a direct limit of $(\M_\a, j_{\a, \b}: \a<\b, \a, \b\in D)$ where $D\subseteq Ord\cap \M$ is cofinal in $Ord\cap \M$ and $\M_\a\insegeq \M$.}. It follows that $\N\insegeq Lp^{cuB, \Lambda}(X)$. 
\end{proof}

\begin{corollary}\label{prop:smallcof}
Assume there is a class of Woodin cardinals and let $\mu$ be a measurable limit of Woodin cardinals that stabilizes uB. Assume $\Gamma^\infty\models _\Omega \sf{AD}_{\bR}$. Let $\mathcal{H}^-, \mathcal{H}$ etc. be defined relative to $\mu$ as in \rsec{sec:uBstacks}. Set \begin{center}$\xi=\max(\cf^V(Ord\cap \mathcal{H}), \cf^V(Ord\cap Lp^{cuB, \Psi_\mu}(A_\mu))$.\end{center} Then $\xi < \mu$.
\end{corollary}
\begin{proof} We show that $\cf^V(Ord\cap \mathcal{H})<\mu$. The second inequality is very similar. Let $g\subseteq Coll(\omega, <\mu)$. Notice that $\card{\Gamma^\infty_g}^{V[g]}=\aleph_1=\mu$. It follows that $\card{\Theta}<\mu^+$ (recall that $\powerset(\bR_g)\cap L(\Gamma^\infty_g, \bR_g)=\Gamma^\infty_g$). The fact that $\card{\mH}^V<\mu^+$ follows from \rcor{failure of lp covering}. The fact that $\cf^V(Ord\cap \mathcal{H})<\mu$ follows from the fact that $\square(\mu)$ fails while letting $\zeta=Ord\cap \mathcal{H}$, $\mH$ has a $\square(\zeta)$-sequence. Let $\vec{C}$ be the $\square(\zeta)$-sequence constructed via the proof of $\square$ in $\mH$\footnote{Notice that this is the easy version of the proof of square, the construction of \cite{Jensen} is all we need.}. If $\cf(\zeta)=\mu$ then $\vec{C}$ has a thread $D$ by measurability of $\mu$; the existence of $D$ follows by an argument similar to that of Corollary \ref{failure of lp covering}. Because of the way $\vec{C}$ is defined, $D$ indexes a sequence of models $(\M_\a: \a\in D)$ such that 
\begin{enumerate}
\item for every $\a\in D$, $\M_\a\insegeq \mH$ and $\rho(\M_\a)=\Theta$, and
\item for $\a<\b, \a, \b\in D$, there is an embedding $\pi_{\a, \b}:\M_\a\rightarrow \M_\b$. 
\end{enumerate}
Let $\M$ be the direct limit along $(\M_\a, \pi_{\a, \b}: \a\in D)$. Then every countable submodel of $\M$ embeds into some $\M_\a$, implying that $\M\insegeq \mH$. However, as $D$ is a thread, $\mH\inseg \M$, contradiction.

\end{proof}

The next proposition shows that sufficiently closed Skolem hulls of $Lp^{cuB}$-operator condense. The proof of it is very much like the proof of \cite[Theorem 10.3]{sargsyan2013covering} and the proof of \cite[Theorem 9.2.6]{hod_mice_LSA}. The proof of \cite[Theorem 9.2.6]{hod_mice_LSA} is done for $\mH$ not $A_\mu$. The proof of \rprop{lem:full_hull}  can be obtained from the proof of \cite[Theorem 9.2.6]{hod_mice_LSA} by simply changing $\P$ to $A_\mu$ everywhere.

\begin{proposition} \label{lem:full_hull} Suppose there is a proper class of Woodin cardinals, $\mu$ is a measurable limit of  Woodin cardinals stabilizing uB and $\Gamma^\infty\models _\Omega \sf{AD}_{\bR}$.  There is then $\nu<\mu$ and $Y_0\in \powerset_{\nu}( H_{\mu^+})$ such that $A_\mu\in Y_0$ and for any $X\prec H_{\mu^+}$ of size $<\mu$ that is closed under $\nu$-sequences and $Y_0\subseteq X$, letting $\pi_X: M_X\rightarrow H_{\mu^+}$ be the uncollapse map and letting $\Lambda$ be the $\pi$-pullback of $\Psi$, $\pi_X^{-1}(Lp^{cuB, \Psi}(A_\mu))=Lp^{cuB, \Lambda}(\pi_X^{-1}(A_\mu))$.
\end{proposition}
\begin{proof} Let $j: V\rightarrow M$ be an embedding witnessing that $\mu$ is a measurable cardinal. Thus, $\cp(j)=\mu$.
It follows from \rcor{prop:smallcof} that $j[Ord\cap Lp^{cuB, \Psi}(A_\mu)]$ is cofinal in $Ord\cap j(Lp^{cuB, \Psi}(A_\mu))$.\footnote{Below, we often confuse strategies with their interpretations in relevant generic extensions or in relevant inner models. However, in some cases, the distinction between the two strategies is important, and in those situations we will either separate the two strategies or point out that the distinction is important.} Below we use the following information: 
\begin{itemize}
\item $\P=Lp^{cuB, \Psi}(A_\mu)$, $N=j(H_{\mu^+})$, $A=j(A_\mu)$, 
\item $Y_0^*=j[\P]$ and $Y_0=Hull^N(Y_0^*)$, and
\item if $Y\in \powerset_{\omega_1}(N)\cap M$ then we let $M_Y$ be the transitive collapse of $Y$, $\pi_Y:M_Y\rightarrow j(H_{\mu^+})$ be the inverse of the transitive collapse, $\P_Y=\pi^{-1}_Y(j(\P))$ and $\Sigma_Y$ be the $\pi_Y$-pullback of $\Psi$,
\item if $Y\subseteq Y'$ then we let $\pi_{Y, Y'}: M_Y\rightarrow M_{Y'}$ be the canonical embedding,
\item  $\mH$ and $\mH^-$ are defined relative to $\mu$ as in \rsec{sec:uBstacks}.
\end{itemize}

We want to show that \\\\
(a) if $Y\in \powerset_{j(\mu)}(N)\cap M$ is such that $Y_0\subseteq Y$ then $\pi_Y^{-1}(j(\P))=Lp^{cuB, \Sigma_Y}(\pi_Y^{-1}(A))$.\\\\
Towards a contradiction assume that (a) is false. Fix one such $Y$ that is a counterexample to (a), and let $\M\insegeq Lp^{cuB, \Sigma_Y}(\pi_Y^{-1}(A))$ be a sound $\Sigma_Y$-mouse over $\pi_Y^{-1}(A)$ such that $\M\not \insegeq \pi_Y^{-1}(j(\P))$ and $\rho(\M)=Ord\cap \pi_Y^{-1}(A)$. We can then find some $\Sigma_{Y_0}$-hod pair $(\P^+, \Pi)\in M$\footnote{Notice that $\Sigma_{Y_0}=\Psi$.} and a hod pair $(\S, \Phi)\in M$ such that
\begin{enumerate}
\item $\P^+\in H_{j(\mu)}^M$ and $\P^+$ is a hod premouse over $A_\mu$ extending $\P$,
\item $\Pi$ has strong branch condensation, 
\item $\P^+$ is meek and of limit type (see \cite[Definition 2.7.1]{hod_mice_LSA}), 
\item $\cf^{\P^+}(\d^{\P^+})=\omega$,
\item $(Y\cap j(\mH^-))\subseteq \rge(\pi^{\Phi}_{\S, \infty})$ and no proper complete layer of $\S$ has this property\footnote{I.e., if $\S'\inseg^c_{hod}\S$ then $(Y\cap j(\P|\d^\P))\not \subseteq \rge(\pi^{\Phi}_{\S', \infty})$. Complete layers are those layers $\S'$ of $\S$ for which $\d^{\S'}$ is a Woodin cardinal of $\S$ or is a limit of Woodin cardinals of $\S$.},
\item $\Pi\in M$ is a $(j(\mu), j(\mu))$-strategy for $\P^+$ such that if $h\subseteq Coll(\omega, <j(\mu))$ is $M$-generic then $\Pi$ can be uniquely extended to a strategy $\Pi^{h}\in (\Gamma^\infty)^{M[h]}$, and moreover, $\Pi$ witnesses that $\P^+$ is a $\Sigma_{Y_0}$-hod mouse.\footnote{For details, the reader may wish to check \cite[Theorem 9.2.6]{hod_mice_LSA}.}
\end{enumerate}
Let $\tau: M_{Y_0}\rightarrow M_Y$ be the canonical embedding, and let $E$ be the long extender of length $Ord\cap \pi_Y^{-1}(Lp^{cuB, \Psi}(A))$ derived from $\tau$. Because $\P^+$ might have cardinality $>\mu$, when we form $\P_Y^+=_{def}Ult(\P^+, E)$, we cannot conclude that $\P_Y^+$ is iterable in $M$. This is because we do not know that $j\rest \P^+\in M$. To resolve this issue we take a hull of size $\mu$. Let $\mu_1=(\mu^+)^V$.

We work in $M$. We can now find $m: W\rightarrow N$ such that
\begin{itemize}
\item $W\in M$ is transitive and $\mu+1\subseteq W$, 
\item $(j(\P), Y_0, Y, (\P^+, \Pi), (\S, \Phi))\in \rge(m)$.
\end{itemize}
Let $Z=m^{-1}(Y)$, $\N=m^{-1}(\M)$, $\R=m^{-1}(j(\P))$ and $k: \P\rightarrow \R$ be $m^{-1}(j\rest \P)$. Working in $M$, set 
\begin{itemize}
\item $\Q=(\P_{Z})^W$, 
\item $\sigma=(\pi_{m^{-1}(Y_0), Z}\rest \P)^{W}$ and $\tau=(\pi_Z\rest \Q)^{W}$,
\item $\overline{\P^+}=m^{-1}(\P^+)$ and $\overline{\Pi}=m^{-1}(\Pi)$,
\item $(\overline{\S}, \overline{\Phi})=m^{-1}(\S, \Phi)$.
\end{itemize}
Thus, we have that\\\\
(A) $k=\tau\circ \sigma$, $\sigma:\P\rightarrow \Q$ and $\tau:\Q\rightarrow \R$,\\\
(B) in $W$, 
\begin{enumerate}
\item $\N$ is a sound $\Sigma_{Z}$-mouse over $\Q$ that projects to $Ord\cap \Q$.
\item  in any derived model of $(\overline{\P^+}, \overline{\Pi})$ as computed by an $\bR$-genericity iteration, $\N$ has an $\omega_1$-iteration strategy witnessing that it is a $\Sigma_Z$-mouse,
\item $\N$ is not an initial segment of $\Q$. 
\item $\overline{\Phi}$ is in the derived model of $(\overline{\P^+}, \overline{\Pi})$ as computed by any $\bR$-genericity iteration,
\item letting $\xi:\sigma(\mH)\rightarrow \overline{\S}|\d^{\overline{\S}}$ be such that $\xi=(\pi^{\overline{\Phi}}_{\overline{\S}, \infty})^{-1} \circ \tau$, $\Sigma_Z=(\xi$-pullback of $\overline{\Phi}_{\overline{\S}|\d^{\overline{\S}}}$. 
\end{enumerate}
Let now $F$ be the long extender of length $\d^\Q$ derived from $\sigma$, and set $\Q^+=Ult(\overline{\P^+}, F)$. Let $\sigma^+=\pi_F^{\P^+}$. Notice that because $m\circ k=j\rest \P$, we have $\phi^+:\Q^+\rightarrow j(\overline{\P^+})$ such that\\\\ (C) $j\rest \overline{\P^+}=\phi^+\circ \sigma^+$.\\\\
Let $\overline{\Pi}^+$ be the $m\rest \overline{\P^+}$-pullback of $\Pi$\footnote{We confuse $\Pi$ with its extension to $N[g]$. Similarly, we think of $\overline{\Pi}^+$ as a strategy in $N$ as well as in $N[g]$. Same comment applies below to $\overline{\Pi}$ and $\overline{\Phi}$.} and let $\overline{\Phi}^+$ be the $m$-pullback of $\Phi$. Notice that\\\\
(D1) $\overline{\Pi}^+\rest W= \overline{\Pi}$\footnote{See proof of Claim 2 in the proof of \cite[Lemma 10.4]{sargsyan2013covering}.},\\
(D2) $\overline{\Pi}^+$ witnesses that $\overline{\P^+}$ is a $\Psi$-hod mouse\footnote{This follows from the fact that $\Pi$ witnesses that $\P^+$ is a $\Psi$-hod mouse and $m\rest \P=id$.},\\
(D3) $\overline{\Phi}^+\rest W= \overline{\Phi}$.\\\\
Notice now that we have \\\\
(F) in $M$, $j(\overline{\Pi}^+\rest H_{\mu_1}^{M})$ is a $(j(\mu), j(\mu))$-iteration strategy witnessing that $j(\overline{\P^+})$ is a $j(\Psi)$-hod mouse, and moreover, $j\rest \overline{\P^+}\in M$\footnote{Because $\card{\overline{\P^+}}=\mu$.}.\\\\
We let $\Gamma=(\Sigma_Z)^{W}$. Notice that in $W$, $\Gamma$ is the $\tau$-pullback of $m^{-1}(j(\Psi))$. Let $\Gamma^+$ be the $\phi^+\rest \Q=m\circ \tau\rest \Q$-pullback of $j(\Psi)$. It follows that\\\\
(G) $\Gamma^+$ is the $m\circ \xi$-pullback of $\Phi$, and it is also $\xi$-pullback of $\overline{\Phi}^+$.\\\\
We now claim that\\\\
(b) in $M$, in any derived model of $(\overline{\P^+}, \overline{\Pi}^+)$ as computed by an $\bR$-genericity iteration, $\N$ has an $\omega_1$-iteration strategy witnessesing that $\N$ is a $\Gamma^+$-mouse.\\\\
The proof of (b) is like the proof of Claim 1 of \cite[Lemma 10.4]{sargsyan2013covering} and it is also very similar to the proof of (b) that appears in the proof of \cite[Theorem 9.2.6]{hod_mice_LSA}. Because of this we skip the proof of (b). 

To finish the proof of \rprop{lem:full_hull}, it remains to implement the last portion of the proof of \cite[Theorem 10.3]{sargsyan2013covering}. Let $\Delta_0$ be $\phi^+$-pullback of $j^+(\overline{\Pi}^+\rest N)$. Notice that it follows from (F) that $\Delta_0$ witnesses that $\Q^+$ is a $\Gamma^+$-hod mouse. It then follows from (b) that\\\\
(H) in $M$, in any derived model of $(\Q^+, \Delta_0)$ as computed by an $\bR$-genericity iteration, $\N$ has an $\omega_1$-iteration strategy $\Delta$ witnessing that $\N$ is a $\Gamma^+$-mouse.\\\\
(H) gives contradiction, as it implies that\\\\
(I) $\Q^+\models ``Ord\cap \Q$ is not a cardinal"\footnote{This is because (K) implies that $\N$ is ordinal definable in $\Q^+$ and therefore, $\N\in \Q$.},\\\\
while clearly $\overline{\P^+}\models ``Ord\cap \P$ is a cardinal", contradicting the elementarity of $\phi^+$.\\\\
\end{proof}

\textit{Proof of \rthm{main thm on stacks under sealing}}\\\\
 We now prove \rthm{main thm on stacks under sealing}. First, take $S$ to be the stationary class of measurable cardinals which are limits of Woodin cardinals; for any $\mu\in S$, $\mu$ satisfies clause (1) of \rthm{main thm on stacks under sealing} by Corollary \ref{prop:smallcof}.  To get clauses (2) and (3), we apply Lemma \ref{lem:pressing_down} to the function $f$ on $S$ that maps each $\mu\in S$ to the maximum of the ordinals $\{\nu, \eta_\mu, \xi\}$, where $\nu$ appears in Proposition \ref{lem:full_hull}, $\eta_\nu$ appears in the proof of Proposition \ref{lem:dm_uB}, and $\xi$ appears in the statement of Corollary \ref{prop:smallcof}. Using Lemma \ref{lem:pressing_down}, we obtain proper class $S_0\subseteq S$ such that for each $\mu\in S_0$, $\nu$  witnesses clauses (2) and (3) of  \rthm{main thm on stacks under sealing} . This finishes the proof of \rthm{main thm on stacks under sealing}. $\qedsymbol{}$\\
\begin{remark}\label{rmk:regressive}
\begin{enumerate}
\item By Lemma \ref{lem:pressing_down}, the existence of $S,S_0,\nu$ above can be proved within $\sf{ZFC}$.
\item It may appear that we use second order set theory to ``pick" for each measurable limit of Woodin cardinals $\mu$ a set $A_\mu$ that codes $V_\mu$, but the theory $\sf{ZFC} + ``$there is (global) well-order of $V$" is conservative over $\sf{ZFC}$. Over any $V\models \sf{ZFC}$, we can find a (class) generic extension $V[g]$ of $V$ such that $V[g]\models ``\sf{ZFC} + $there is a global well order". 
\item The above two remarks simply say that  we may assume as part of the hypothesis that $V$ has a global well-order. This then allows us to get $S, S_0, \nu$ and the sequences $(Y_\mu: \mu \in S_0)$, $(A_\mu: \mu \in S)$ in \rthm{main thm on stacks under sealing}. 
\end{enumerate}
\end{remark}

The rest of the argument does not need the hypothesis that $\Gamma^\infty\models _{\Omega} \sf{AD}_{\bR}$. It only needs the conclusion of \rthm{main thm on stacks under sealing}. To stress this point we make the following definitions.

\begin{definition}\label{theory t} We let $T$ stand for the following theory.
\begin{enumerate}
\item $T_0$
\item There is a stationary class $S$, a proper class $S_0\subseteq S$, an infinite regular cardinal $\nu$ and two sequences $\vec{Y}=(Y_\mu: \mu\in S_0)$ and $\vec{A}=(A_\mu: \mu\in S)$ such that the following conditions hold for any $\mu\in S$
\begin{enumerate}
\item $\mu$ is a measurable limit of Woodin cardinals,
\item $\mu$ stabilizes uB,
\item $\card{\mH_\mu}<\mu^+$, 
\item $A_\mu\subseteq \mu$ codes $V_\mu$ and $\max(\cf(Ord\cap \mH_\mu), \cf(Ord\cap Lp^{cuB, \Psi_\mu}(A_\mu))<\mu$;
\setcounter{nameOfYourChoice}{\value{enumi}}
\end{enumerate}
furthermore, if $\mu\in S_0$, then the following hold:
\begin{enumerate}
\setcounter{enumi}{\value{nameOfYourChoice}}
\item $\max(\cf(Ord\cap \mH_\mu), \cf(Ord\cap Lp^{cuB, \Psi_\mu}(A_\mu))<\nu$,
\item $Y_\mu\in \powerset_\nu(H_{\mu^+})$,
\item $A_\mu\in Y_\mu$, and
\item whenever $X\prec H_{\mu^+}$ is of size $<\mu$, is $\nu$-closed and $Y_\mu\subseteq X$, $X$ captures\footnote{See \rdef{x captures}.} $Lp^{cuB, \Psi_\mu}(A_\mu)$. 
\end{enumerate}
\end{enumerate}
\end{definition}

\section{Condensing sets}

Here we review some facts about \textit{condensing sets} that were introduced in \cite{sargsyan2013covering} and developed further in \cite[Chapter 9.1]{hod_mice_LSA}. We develop this notion assuming the theory $T$ introduced in \rdef{theory t}. Let $(S, S_0, \nu_0,\vec{Y}, \vec{A})$   witness that $T$ is true. 

Fix $\mu\in S_0$ and let $g\subseteq Coll(\omega, <\mu)$ be generic.  We let $ \mathcal{H}, \Psi$ etc. stand for the CMI objects associated with $\mu$. We summarize some basic notions and results concerning condensing sets which will play a key role in our $K^c$-constructions. \cite[Chapter 9]{hod_mice_LSA} gives more details and proofs of basic facts about these objects.

The notion of fullness that we will use is full in $L(\Gamma^\infty_g, \bR_g)$. Notice that if $\Phi\in L(\Gamma^\infty_g, \bR_g)$ is an $\omega_1$-strategy with hull condensation then in $L(\Gamma^\infty_g, \bR_g)$,  for any $x\in \bR_g$, $OD(\Phi)$ is the stack of $\omega_1$-iterable $\Phi$-mice over $x$\footnote{This is an instance of the Mouse Set Conjecture, which is not known in full generality. However, we are working towards establishing the equiconsistency in \rthm{thm:main_theorem}. But the target large cardinal is weak, and so Mouse Capturing holds in $L(\Gamma^\infty_g, \bR_g)$. See \cite[Chapter 10.2]{hod_mice_LSA}.}. Because any such $\Phi$-mouse has an iteration strategy in $\Gamma^\infty_g$, it follows that ``full in $L(\Gamma^\infty_g, \bR_g)$" is equivalent to ``full with respect to $Lp^{cuB}$ in $V[g]$".
Thus, given $\M\in HC^{V[g]}$ we say $\M$ is $\Phi$-full if for any $\M$-cutpoint $\eta$, $Lp^{cuB, \Phi}(\M|\eta)\in \M$. If $\M$ is a $\Phi$-mouse over $\M|\eta$ then by ``$\M$ is $\Phi$-full" we in fact mean that $\M|(\eta^+)^\M=Lp^{cuB, \Phi}(\M|\eta)$. Here we note again that ``$Lp^{cuB,\Phi}$" is computed in $V[g]$.

We start working  in $V[g]$. Following  \cite[Chapter 9]{hod_mice_LSA}, for each $Z\subseteq \mH$, we let: 
\begin{itemize}
\item $\Q_Z$ be the transitive collapse of $Hull_1^\mathcal{H}(Z)$, 
\item $\tau_Z:\Q_Z\rightarrow \mathcal{H}$ be the uncollapse map, and
\item $\delta_Z = \delta^{\Q_Z}$, where $\tau_Z(\delta_Z) = \Theta = \delta^\mathcal{H}$. 
\end{itemize}

For $X\subseteq Y \in \powerset_{\omega_1}(\mathcal{H})$, let 
\begin{center}
$\tau_{X, Y}= \tau^{-1}_Y \circ \tau_X$
\end{center}

\begin{definition}\label{def_simext}
Let $Z\in \powerset_{\omega_1}(\mathcal{H})$. $Y \in \powerset_{\omega_1}(\mathcal{H}^{-})$ is a \textbf{simple extension} of $Z$ if
\begin{center}
$Hull_1^{\mathcal{H}} (Z \cup Y )\cap \mathcal{H}^- \subseteq Y.$
\end{center}
\end{definition}

Let $Z,Y$ be as in Definition \ref{def_simext}. Let 
\begin{equation}
Y \oplus Z  = Hull_1^{\mathcal{H}}(Z\cup Y)
\end{equation}

Let 
\begin{equation}
\tau^Z_Y = \tau_{Z\oplus Y},
\end{equation}
and
\begin{center}
$\pi^Z_Y: \Q_Z \rightarrow \Q_{Z\oplus Y}$ be  $\tau_{Z,Y\oplus Z}$; 
\end{center}
we also write $\Q^Z_Y$ for $\Q_{Z\oplus Y}$ and $\delta^Z_Y$ for $\delta^{\Q_{Z\oplus Y}}$.
We have that 
\begin{equation}
\tau_Z = \tau^Z_{Y}\circ \pi^Z_Y
\end{equation}

Given two simple extensions of $Z$, $Y_0 \subseteq Y_1$, we let $\pi^Z_{Y_0,Y_1}:\Q^Z_{Y_0} \rightarrow \Q^Z_{Y_1}$ be the natural map.  
We also let
\begin{center}
$\Psi^Z_Y = \tau_Y^Z$-pullback of $\Psi$.
\end{center}

\begin{definition}\label{def_ext}
$Y$ is an \textbf{extension} of $Z$ if $Y$ is a simple extension of $Z$ and $\pi^Z_Y\rest ( \Q_Z|\delta_Z)$ is the iteration embedding according to $\Psi^Z_Y$. Here we allow $Z$ to be an extension of itself.
\end{definition}

Suppose $Y$ is an extension of $Z$. Let $\sigma^Z_{Y}:\Q^Z_Y\rightarrow \mathcal{H}$ be given by
\begin{equation}
\sigma^Z_{Y}(q) = \tau_Z(f)(\pi^{\Psi^Z_Y}_{\Q^Z_Y,\infty}(a))
\end{equation}
where $a\in (\Q^Z_Y|\delta^Z_Y)^{<\omega}$ and $q = \pi^Z_{Y}(f)(a)$.

\begin{definition}\label{def_honext}
$Y$ is an \textbf{honest extension} of $Z$ if 
\begin{enumerate}
\item $Y$ is an extension of $Z$,
\item $\dom(\sigma_Y^Z)=\Q^Z_Y$ and $\sigma_Y^Z$ is elementary,
\item $\tau_Z =\sigma_Y^Z \circ \pi_Y^Z$\footnote{This condition follows from other conditions.}.
\end{enumerate}
\end{definition}

We say $Y$ is an \textit{iteration extension} of $Z$ if $Y$ is an honest extension of $Z$ and $Y=\sigma^Z_Y[\Q^Z_Y|\delta^Z_Y]$.

\begin{definition}\label{def_simcond}
We say $Z$ is a \textbf{simply condensing set} if
\begin{enumerate}
\item for any extension $Y$ of $Z$, $\Q_Y^Z$ is $\Psi^Z_Y$ -full, 
\item all extensions $Y$ of $Z$ are honest.
\end{enumerate}

We say Z is \textbf{condensing} if for every extension $Y$ of $Z$, $Z \oplus Y$ is a simply condensing set. 

\end{definition}

In $V[g]$, let 
\begin{center}
$Cnd(\mathcal{H}) = \{Z \in \powerset_{\omega_1}(\mathcal{H}) : Z \textrm{ is condensing}\}$.
\end{center}

\noindent Results in \cite[Chapter 9]{hod_mice_LSA} give
\begin{theorem}\label{thm:club}
In $V[g]$, $Cnd(\mathcal{H})$ is a club in $\powerset_{\omega_1}(\mathcal{H})$ (i.e. it is unbounded and is closed under countable unions). 

Furthermore, for any cardinal $\kappa\geq \nu_0$ and $\kappa<\mu$, $\{X\in V: X \in Cnd(\mathcal{H})\wedge |X|^V \leq \kappa \} $ is a club in $\powerset^V_{\kappa^+}(\mathcal{H})$. The same holds if $V$ is replaced by $V[g\cap Coll(\omega, <\kappa)]$.

Furthermore, for each $Z\in Cnd(\mH)$, if $Y$ is an honest extension of $Z$, then $Y$ is an iteration extension of $Z$.	
\end{theorem}

Also the following uniqueness fact is very important for this paper. It follows from \rprop{simple structure} and can be proved exactly the same way as \cite[Lemma 9.1.14]{hod_mice_LSA}.

\begin{proposition}\label{independence of psi} Suppose $Z$ is a condensing set. Suppose $Y$ and $W$ are extensions of $Z$ such that $\Q^Z_Y=\Q^Z_W$. Then $\Psi^Z_Y=\Psi^Z_W$.  
\end{proposition}

The following are easy corollaries of \rprop{independence of psi}.

\begin{corollary}\label{unique honest extensions} Suppose $Z$ is a condensing set and $\Q$ is such that for some extension $Y$ of $Z$, $\Q=\Q^Z_Y$. There is then a unique honest extension $W$ of $Z$ such that $\Q=\Q^Z_W$. 
\end{corollary}

\begin{corollary}\label{invariance under pullbacks} Suppose $Z$ is a condensing set. Suppose further that $Y$ and $W$ are two extensions of $Z$ such that there is an embedding $i:\Q^Z_Y\rightarrow_{\Sigma_1} \Q^Z_W$ such that $\tau^Z_W\circ i[\Q^Z_Y]$ is an extension of $Z$. Then the $i$-pullback of $\Psi^Z_W$ is $\Psi^Z_Y$. 
\end{corollary}
\begin{proof} Let $Y^*=\tau^Z_W\circ i[\Q^Z_Y]$. We have that $\Q^Z_{Y}=\Q^Z_{Y^*}$. Moreover, $\Psi^Z_Y=\Psi^Z_{Y^*}$ and $\Psi^Z_{Y^*}$ is the $i$-pullback of $\Psi^Z_{W}$. 
\end{proof}
  
\section{$Z$-realizable iterations}

In this section we fix a condensing set $Z$.

\begin{definition}\label{def_niceext}
Let $Z$ be a condensing set. $\Q$ \textbf{nicely extends} $\Q_Z$ if $\Q$ is non-meek\footnote{See 
\cite[Definition 2.7.1]{hod_mice_LSA}. $\Q$ is meek if either it has successor type or $\Q=\Q^b$. Otherwise, we say $\Q$ is non-meek.} and $\Q^b = \Q_Z$. We also say that $\Q$ is a \textbf{nice extension} of $\Q_Z$.
\end{definition} 
Suppose  $Y$ is an extension of $Z$ and $\Q$ nicely extends $\Q^Z_Y$. We would like to analyze the stacks on $\Q$, following the terminology and conventions used in \cite{hod_mice_LSA}. A stack\footnote{A stack of normal iteration trees.} $\T$ on $\Q$ has the form
\begin{center}
$\VT=((\M_\a)_{\a<\eta}, (E_\a)_{\a<\eta-1}, D, R, (\beta_\a, m_\a)_{\a\in R}, T)$,
\end{center}
where the displayed objects are introduced in \cite[Definition 2.4.1]{hod_mice_LSA}. The above notation is quite standard. $D$ is the set of drops, $R$ is the set of stages where player $I$ starts a new round of the iteration game, $(\beta_\a, m_\a)$ is the place player $I$ drops at the beginning of the $\a$th round, and $T$ is the tree order. We adopt an important convention introduced in \cite{hod_mice_LSA}. Namely, we assume that all our stacks are proper (see \cite[Remark 2.7.27]{hod_mice_LSA}). One of the key aspects of being proper is that if $\b<lh(\VT)$ is such that $\VT_{\geq \b}$ is a stack on $\M^\VT_\b$ then $\b\in R$\footnote{Thus, no normal component of $\VT$ can be split into two normal components.}. We will also use the notation introduced in \cite[Notation 2.4.4]{hod_mice_LSA}. In particular, for $\a\in R^\VT$, ${\sf{next}}^{\VT}(\a)=min(R^{\VT}-(\a+1))$ if this minimum exists and otherwise ${\sf{next}}^{\VT}(\a)=lh(\VT)$. For $\a\in R^\VT$, we also set ${\sf{nc}}^{\VT}_\a=\VT_{[\a, \a']}$ where $\a'={\sf{next}}^{\VT}(\a)$.


\begin{definition} \label{dfn:realizable_ext}Suppose $Z\in Cnd(\mH)$ and $Y$ is an extension of $Z$. Suppose further that $\Q$ nicely extends $\Q^Z_Y$. Given $E\in \vec{E}^\Q$ such that $\cp(E)=\d^{\Q^Z_Y}$, we say $E$ is \textbf{$(Z, Y)$-realizable} if there is $W$, an extension of  $Z\oplus Y$ such that $E=E^Z_{Y, W}$, where $E^Z_{Y,W}$ is the extender defined by\:
\begin{equation}\label{dfn:extender}
(a,A)\in E^Z_{Y,W} \Leftrightarrow \tau^Z_{W}(a) = \pi^{\Psi^Z_W}_{\Q^Z_W,\infty}(a)\in \tau^Z_Y(A),
\end{equation}
for any $a\in [lh(E)]^{<\omega}$ and $A\in \powerset(\cp(E))^{|a|}\cap \Q$. 

We are continuing with the notation of \rdef{dfn:realizable_ext}. Suppose $\VT$ is a stack on $\Q$. We say $\VT$ is a \textbf{$(Z, Y)$-realizable} iteration if there is a sequence $(W_\a: \a \in R^{\VT})$ such that
\begin{enumerate}
\item $W_0=Y$,
\item if $\a, \b\in R^{\VT}$ and $\a<\b$ then $W_\b$ is an extension of $Z\oplus W_\a$, 
\item if $\a, \b\in R^{\VT}$, $\a<\b$ and $\pi^{\VT, b}_{\a, \b}$ is defined then $\pi^{\VT, b}_{\a, \b}=\pi^Z_{W_\a, W_\b}$\footnote{The embedding $\pi^{\VT}_{\a, \b}$ is defined similarly to $\pi^{\VT, b}$, it is essentially the embedding $\pi^{\VT}_{\a, \b}\rest \M_\a^b$. See \cite[Chapter 2.8]{hod_mice_LSA}.}, and
\item if $\a\in R^\VT$ and $\VU$ is the largest fragment of $\VT_{\geq \a}$ that is based on $\M_\a^b$ then $\VU$ is according to $\Psi^Z_{W_\a}$.
\end{enumerate}
We say $\VT$ is \textbf{$Z$-realizable} if $Y$ is an honest extension of $Z$ and $\VT$ is  $(Z, Y)$-realizable. 
\end{definition}

The following lemma is a consequence os \rprop{independence of psi}, \rcor{unique honest extensions} and \rcor{invariance under pullbacks}.

\begin{lemma}\label{unique z-realizability} Suppose $Y$ is an extension of $Z$ and $\Q$ nicely extends $\Q^Z_Y$. Suppose $\VT$ is a $(Z, Y)$-realizable iteration as witnessed by $(W'_\a: \a \in R^{\VT})$. For $\a\in R^{\VT}$ let $W_\a$ be the unique honest extension of $Z$ with the property that $(\M_\a^{\VT})^b=\Q^Z_{W_\a}$. Then $(W_\a: \a\in R^{\VT})$ witnesses that $\Q$ is $Z$-realizable. 
\end{lemma}
\begin{proof} It is enough to show that if $\a, \b\in R^\VT$ and $\b=\min(R^\VT-(\a+1))$ then $\pi^\VT_{\a, \b}=\pi^Z_{W_\a, W_\b}$. First we show that $W_\a\subseteq W_\b$. We have that $x\in W_\a$ if for some $a\in \d_{Z\oplus W_\a}$ and some $f\in Z$, $x=\tau_Z(f)(\tau^Z_{W_\a}(a))$. Since $\tau^Z_{W_\a}\rest \d_{Z\oplus W_\a}$, $\pi^{\VT, b}_{\a, \b}\rest \d_{Z\oplus W_\a}$ and $\tau^Z_{W_\b}\circ \pi^{\VT, b}_{\a, \b}\rest \d_{Z\oplus W_\a}$ are all iteration embeddings according to $\Psi^Z_{W_\a}$, we have that $x=\tau_Z(f)(\tau^Z_{W_\b}(\pi^{\VT, b}_{\a, \b}(a))$. Thus, $W_\a\subseteq W_\b$. A similar argument shows that $\pi^{\VT, b}_{\a, \b}=\tau^Z_{W_\a, W_\b}$.
\end{proof}

\begin{remark} It follows from \rlem{unique z-realizability} that $(Z, Y)$-realizability is equivalent to $Z$-realizability. Because of this, in this paper, we will mostly use $Z$-realizability.
\end{remark}

Suppose $\Q$ nicely extends $\Q^Z_Y$ and $\VT$ is a $Z$-realizable iteration of $\Q$. We cannot in general prove that $\VT$ picks unique branches mainly because we say nothing about $\Q$-structures that appear in $\VT$ when we iterate above $\d^{\M_\a^b}$ for some $\a\in R^\VT$. The next definition introduces a notion of a premouse that resolves this issue. 

\begin{definition}\label{weakly suitable} We say $\R$ is \textbf{weakly $Z$-suitable} if $\R$ is a hod premouse of lsa type such that $\R=(\R|\d^\R)^{\#}$, $\R$ has no Woodin cardinals in the interval $(\d^{\R^b}, \d^\R)$ and for some extension $Y$ of $Z$, $\R$ nicely extends $\R^b=\Q^Z_Y$. 
\end{definition} 

The following lemma says that hulls of $Z$-realizable iterations are $Z$-realizable, and easily follows from \rcor{invariance under pullbacks}.

\begin{proposition}\label{hull of realizable iterations} Suppose $\R$ and $\S$ are weakly $Z$-suitable hod premice. Suppose further that $\VT$ is a $Z$-realizable iteration of $\S$ and $\VU$ is an iteration of $\R$ such that $(\R, \VU)$ is a hull\footnote{In the sense of \cite[Definition 1.30]{hod_mice}.} of $(\S, \VT)$. Then $\VU$ is $Z$-realizable. 
\end{proposition}

We now define the notion of \textit{$Z$-approved sts premouse of depth $n$} by induction on $n$. The induction ranges over all weakly $Z$-suitable hod premice. 
\begin{definition}\label{z approved of depth n} Suppose $\R$ is weakly $Z$-suitable hod premouse and for some extension $Y$ of $Z$, $\R$ nicely extends $\R^b=\Q^Z_Y$.
\begin{enumerate}
\item  We say that $\M$ is a \textbf{$Z$-approved sts premouse} over $\R$ of depth $0$ if $\M$ is an sts premouse over $\R$\footnote{This in particular means that the strategy indexed on the sequence of $\M$ is a strategy for $\R$.} such that if $\T\in \M$ is according to $S^\M$ then $\T$ is $(Z, Y)$-realizable.
\item Suppose $\R$ is weakly $Z$-suitable hod premouse. We say that $\M$ is a \textbf{$Z$-approved sts premouse} over $\R$  of depth $n+1$ if $\M$ is a $Z$-approved sts premouse over $\R$ of depth $n$ such that if $\T\in \M$ is $\sf{nuvs}$ and $S^\M(\T)$ is defined then letting $b=S^\M(\T)$, $\Q(b, \T)$ is a $Z$-approved sts premouse over $\m^+(\T)$ of depth $n$.
\end{enumerate}
\end{definition}

\begin{definition}\label{z-validated sts mice} We say $\M$ is a \textbf{$Z$-approved sts premouse} over $\R$ if for each $n<\omega$, $\M$ is a $Z$-approved sts premouse over $\R$ of depth $n$. We say $\M$ as above is a \textbf{$Z$-approved sts mouse} (over $\R$) if $\M$ has a $\mu$-strategy $\Sigma$ such that whenever $\N$ is a $\Sigma$-iterate of $\M$, $\N$ is a $Z$-approved sts premouse over $\R$. 
\end{definition}

The following proposition is an immediate consequence of our definitions, but perhaps is a bit tedious to prove.

\begin{proposition}\label{preservation of zv under embeddings} Suppose $\R$ and $\S$ are weakly $Z$-suitable, $\N$ is an sts premouse over $\R$ and $\M$ is a $Z$-approved premouse (mouse) over $\S$. Suppose $\pi:\N\rightarrow_{\Sigma_1} \M$. Then $\N$ is also a $Z$-approved premouse (mouse). 
\end{proposition}
\begin{proof} We only show that if $\T^*\in \N$ is according to $S^\N$ then $\T^*$ is $Z$-realizable. Even less, we show that if $\T^*=\T^\frown \U$ is such that $\pi^{\T, b}$ exists and $\U$ is based on $\pi^{\T, b}(\R^b)$ then there is an extension $Y$ of $Z$ such that $\pi^{\T, b}(\R^b)=\Q^Z_{Y}$ and $\U$ is according to $\Psi^Z_{Y}$. The rest of the proof is very similar.

Notice that by elementarity of $\pi$, $\pi(\T^*)$ is according to $S^\M$. Therefore, there is some extension $W$ of $Z$ such that $\pi(\pi^{\T, b}(\R^b))=\Q^Z_{W}$ and $\pi(\U)$ is according to $\Psi^Z_W$. Set $Y=\tau^Z_W\circ \pi[\pi^{\T, b}(\R^b)]$. Then $\pi^{\T, b}(\R^b)=\Q^Z_{Y}$ and $\U$ is according to the $\pi$-pullback of $\Psi^Z_{W}$. As the $\pi$-pullback of $\Psi^Z_{W}$ is just $\Psi^Z_{Y}$, we are done.
\end{proof}

\begin{definition}\label{z approved}
Suppose $\R$ is weakly $Z$-suitable. We let $Lp^{Za, sts}(\R)$ be the union of all $Z$-approved sound sts mice $\M$ over $\R$ such that $\rho(\M) \leq Ord\cap \R$. 
\end{definition}

Finally, we can define the \textit{correctly guided $Z$-realizable iterations}.

\begin{definition}\label{correctly guided iterations} Suppose $\R$ is a weakly $Z$-suitable hod premouse and $\VT$ is a $Z$-realizable iteration of $\R$. We say $\VT$ is \textbf{correctly guided} if whenever $\a\in R^\VT$, $\U=_{def}{\sf{nc}}^{\VT}_\a$ is above $\d^{\M_\a^b}$, $\a<lh(\U)$ is a limit ordinal such that $\m^+(\U\rest \a)\models ``\d(\U\rest \a)$ is a Woodin cardinal", then letting $b=[0, \a]_\U$, $\Q(b, \U\rest\alpha)$ is a $Z$-approved sts mouse over $\m^+(\U\rest \a)$.
\end{definition} 

Combining \rprop{hull of realizable iterations} and \rprop{preservation of zv under embeddings} we get the following.

\begin{corollary}\label{hulls of corectly guided iterations} Suppose $\R$ and $\S$ are weakly $Z$-suitable hod premice. Suppose further that $\VT$ is a correctly guided $Z$-realizable iteration of $\S$ and $\VU$ is an iteration of $\R$ such that $(\R, \VU)$ is a hull of $(\S, \VT)$ (in the sense of \cite[Definition 1.30]{hod_mice}). Then $\VU$ is also correctly guided $Z$-realizable iteration. 
\end{corollary}

Our uniqueness theorem applies to $\R$ that are not \textit{infinitely descending}.

\begin{definition}\label{inf desc} We say that a weakly $Z$-suitable hod premouse $\R$ is \textbf{infinitely descending} if there is a sequence $(p_i, \R_i, Y_i: i<\omega)$ such that 
\begin{enumerate}
\item $\R_0=\R$,
\item for every $i<\omega$, $\R_i$ is weakly $Z$-suitable,
\item for every $i<\omega$, $p_i$ is a correctly guided $Z$-realizable iteration of $\R_i$,
\item for every $i<\omega$, $p_i$ has a last normal component $\T_i$ of successor length\footnote{I.e., $R^{p_i}$ has a maximum element $\a$ and $\T_i=(p_i)_{\geq \a}$.} such that $\a_i=_{def}lh(\T_i)-1$ is a limit ordinal and 
$\R_{i+1}=\m^+(\T_i\rest \a_i)$,
\item for every $i<\omega$, setting $b_i=_{def}[0, \a_i)_{\T_i}$, $b_i$ is a cofinal branch of $\T_i$ such that $\Q(b_i, \T_i)$ exists and is $Z$-approved.
\end{enumerate}
\end{definition} 

Note that in the above definition, for each $i$, $\R_{i+1}$ is a strict initial segment of $\Q(b_i, \T_i)$. The following is our uniqueness result. 
\begin{proposition}\label{branch uniqueness} Suppose $\R$ is a weakly $Z$-suitable hod premouse that is not infinitely descending and $\VT$ is a correctly guided $Z$-realizable iteration of limit length on $\R$. Then there is at most one branch $b$ of $\VT$ such that $\VT^\frown \{b\}$ is correctly guided and $Z$-realizable.
\end{proposition}

Here we outline the proof. First notice that\\\\
(a) if $\VT$ doesn't have a last component or\\
(b) if there is $\a\in R^\VT$ such that $\VT_{\geq\a}$ is based on $\M_\a^b$\\\\
 then there is nothing to prove as letting $W_\S$ be as in \rdef{dfn:realizable_ext}, $\Psi^Z_{W_S}$ only depends on $\S^b$ (e.g. see \cite[Lemma 9.1.9]{hod_mice_LSA})\footnote{Notice that in this case there is a branch $b$ such that $\VT^\frown \{b\}$ is correctly guided and $Z$-realizable.}. Let now $\T={\sf{nc}}^\VT_\a$ be the last normal component of $\VT$. If $b, c$ are two different branches of $\T$ such that $\VT^\frown\{b\}$ and $\VT^\frown \{c\}$ are correctly guided $Z$-realizable iterations then $\Q(b, \T)\not =\Q(c, \T)$ and both are $Z$-approved sts mice over $\m^+(\T)$. It now follows from \cite[Lemma 4.7.2]{hod_mice_LSA} and the fact that $\R$ is not infinitely descending that we can reduce the disagreement of $\Q(b, \T)$ and $\Q(c, \T)$ to a disagreement between $\Psi^Z_X$ and $\Psi^Z_U$ for some extensions $X,U$ of $Z$ with $\Q^Z_X=\Q^Z_U$. However, this cannot happen by \rprop{independence of psi} (the proof is given by \cite[Lemma 9.1.9]{hod_mice_LSA}).

\section{$Z$-validated iterations}\label{sec:z-val}

We continue by assuming $T$. Let $(S, S_0, \nu_0,\vec{Y}, \vec{A})$ again witness that $T$ is true and let $\mu, g, \mathcal{H}$ etc. be defined as in \rsec{sec:uBstacks}. The goal of this section is to introduce some concepts to be used in the $K^c$ construction of the next section. The main new concept here is the concept of $Z$-validated iterations which are the kind of iterations that will appear in the $K^c$ construction of the next section.

The following definition is important for this paper. It introduces the hulls that we will use to $Z$-validate mice, iterations, etc. It goes back to Steel's \cite{PFA}.

\begin{definition} \label{def_goodhull} Suppose $\l\in S_0-\mu$ and $U \prec_1 H_{\l^+}$. We say $U$ is $(\mu, \l, Z)$-\textbf{good} if $\mu\in U$, $(Y_\mu\cup Y_\l \cup Z)\subseteq U$, $\card{U} <\mu$, $U\cap \mH^-$ is an honest extension of $Z$ and $U^{\nu_0}\subseteq U$. When $\mu$ and $\l$ are clear from the context or are not important, we simply say $U$ is a good hull. We say a good hull $U$ is transitive below $\mu$ if $U\cap \mu\in \mu$
\end{definition}

If $U$ is a good hull then we let $\pi_U : M_U = M \rightarrow H_{\l^+}$ be the inverse of the transitive collapse of $U$. If in addition $U$ is transitive below $\mu$, we let $\pi_U^+:M_U[g_\nu]\rightarrow H_\Omega[g]$ where $\nu=\cp(\pi_U)$ and $g_\nu=g\cap Coll(\omega, <\nu)$. 

\begin{definition}\label{z-validated sts mice above mu} Suppose
\begin{itemize}
\item $\R_0$ nicely extends $\mH$, 
\item $p$ is an iteration of $\R_0$,
\item if $p$ is ${\sf{nuvs}}$ then setting $\R=\m^+(p)$, $\M$ is an sts mouse over $\R$,  and 
\item $\l\in S$ is the least such that $(\R, \M, p)\in H_\l$. 
\end{itemize}
\begin{enumerate}
\item We say $\R$ is not infinitely descending if whenever  $U$ is a $(\mu, \l, Z)$-good hull such that $\R\in U$, $\pi_U^{-1}(\R)$ is not infinitely descending.
\item We say \textbf{$p$ is $Z$-validated} if whenever $U$ is a $(\mu, \l, Z)$-good hull such that $\{\R, p\}\subseteq U$, $\pi_U^{-1}(p)$ is a correctly guided\footnote{See \rdef{correctly guided iterations}.} $Z$-realizable iteration of $\pi_U^{-1}(\R)$.
\item Suppose $\R$ is a weakly $Z$-suitable hod premouse above $\mu$. We say \textbf{$\M$ is a $Z$-validated sts premouse} over $\R$ if  for every $(\mu, \l, Z)$-good hull $U$  such that $\{\R, \M\}\subseteq U$, letting $\N=\pi_U^{-1}(\M)$, $\N$ is a $Z$-approved\footnote{See \rdef{z-validated sts mice}.} sts premouse over $\pi^{-1}_U(\R)$. 
\item Suppose $\M$ is a $Z$-validated sts mouse over $\R$ and $\xi$ is an ordinal. We say $\M$ has a $Z$-validated $\xi$-iteration strategy if there is $\Sigma$ such that $\Sigma$ is a $\xi$-iteration strategy and whenever $\N$ is an iterate of $\M$ via $\Sigma$, $\N$ is a $Z$-validated sts mouse over $\R$. 
\item Suppose $q$ is an iteration of $\R$. We say $q$ is \textbf{$Z$-validated} if $p^\frown q$ is $Z$-validated.
\end{enumerate}
\end{definition}

The following proposition is very useful and is an immediate consequence of \rprop{preservation of zv under embeddings}. When $X$ is a good hull we will use it as a subscript to denote the $\pi_X$-preimages of objects that are in $X$.

\begin{proposition}\label{one hull witness for premice} Suppose $(\R_0, p, \R, \M, \l)$ are as in \rdef{z-validated sts mice above mu}. Suppose $U$ is a $(Z, \mu, \l)$-good hull such that $\{\R, \M\}\subseteq U$ and  $\M_U$ is not $Z$-approved. Then whenever $U^*$ is a $(Z, \mu, \l)$-good hull such that $U\cup\{U\}\subseteq U^*$, $\M_{U^*}$ is not $Z$-approved\footnote{Hence, $\M$ is not $Z$-validated.}.
\end{proposition}

Similarly for iterations.

\begin{proposition}\label{one hull witness for iterations} Suppose $(\R_0, p, \l)$ are as in \rdef{z-validated sts mice above mu}. Suppose $U$ is a $(\mu, \l, Z)$-good hull such that $\{\R_0, p\}\subseteq U$ and  $p_U$ is not $Z$-realizable. Then whenever $U^*$ is a $(\mu, \l, Z)$-good hull such that $U\cup\{U\}\subseteq U^*$, $p_{U^*}$ is not $Z$-realizable\footnote{Hence, $p$ is not $Z$-validated.}.
\end{proposition}

\begin{definition}\label{z stack}
 Suppose $(\R_0, p, \R, \mu)$ are as in \rdef{z-validated sts mice above mu}. We let $Lp^{Zv, sts}(\R)$ be the union of all $Z$-validated sound sts mice over $\R$ that project to $Ord\cap \R$.
 \end{definition}
 
The following proposition is a consequence of \rprop{branch uniqueness}.

\begin{proposition}\label{unique branches above mu}  Suppose $(\R_0, p, \R, \mu)$ are as in \rdef{z-validated sts mice above mu} and $\R$ is not infinitely descending. Suppose $\VT$ is a $Z$-validated iteration of $\R$ of limit length. Then there is at most one branch $b$ of $\VT$ such that $\VT^\frown\{b\}$ is $Z$-validated.
\end{proposition}

 \section{Realizability array}\label{sec: real array}

We continue with $(S, S_0, \nu_0,\vec{Y}, \vec{A})$, $\mu\in S_0$ etc. We define the notion of an \textit{array} at $\mu$ by induction. We say $\vec{\V}=\V_0$ is an array of length $0$ if $\V_0=\mH$. Suppose we have already defined the meaning of \textit{array} of length $<\eta$. We want to define the meaning of array of length $\eta$. 

\begin{definition}\label{array} We say $\vec{\V}=(\V_\a:\a\leq \eta)$ is an \textbf{array} of length $\eta$ at $\mu$ if the following conditions hold. 
\begin{enumerate}
\item For every $\a<\eta$, $(\V_\b:\b\leq \a)$ is an array of length $\a$ at $\mu$.
\item $\V_\eta$ nicely extends $\mH$ and is a hod premouse.
\item For all $\a<\eta$, if $\V_\a$ is weakly $Z$-suitable then there is $\b\leq \eta$ such that $\V_\b$ is a $Z$-validated sts mouse over $\V_\a$ and $rud(\V_\b)\models ``$there are no Woodin cardinals $>\d^{\mH}$".
\item For all $\a<\eta$, if $rud(\V_\a)\models ``$there are no Woodin cardinals $>\d^{\mH}$" then $\V_\a$ has a $Z$-validated iteration strategy. 
\end{enumerate}
We say $\vec{\V}$ is \textbf{small} if $rud(\V_\eta)\models ``$there are no Woodin cardinals $>\d^\mH$". We let $\eta=lh(\vec{\V})$ and for $\a\leq \eta$, we let $\vec{\V}\rest \a=(\V_\b:\b\leq \a)$.
\end{definition}

Recall the notions of $k$-maximal iteration trees in \cite[Definition 3.4]{steel2010outline}, weak $k$-embeddings \cite[Definition 4.1]{steel2010outline}. For an iteration tree $\T$ on $\M$, letting $\M^\T_\alpha$ be the $\alpha$-th model in the tree; for $\alpha+1<lh(\T)$, recall the notion of degree $deg^\T(\alpha+1)$ \cite[Definition 3.7]{steel2010outline}. Recall the definition of $D^\T$:  if $\alpha+1\in D$, then the extender $E^\T_{\alpha+1}$ is applied to a strict initial segment of $\M^\T_\beta$ where $\beta = T-pred(\alpha+1)$. For $\lambda$ limit, $deg^\T(\lambda)$ is the eventual values of $deg^\T(\alpha+1)$ for $\alpha+1\in [0,\lambda]_\T$. For a cofinal branch $b$ of $\T$, $deg^\T(b)$ is defined to be the eventual value of $deg^\T(\alpha+1)$ for $\alpha+1\in b$. We write $\mathcal{C}_k(\M)$ for the $k$-th core of $\M$. Sometimes, we confuse $\mathcal{C}_0(\M)$ with $\M$ itself.

\begin{definition}\label{realizability property} Suppose $\vec{\V}$ is an array at $\mu$. We say $\vec{\V}$ has the $Z$-\textbf{realizability property} if for all $\a<lh(\V)$, $\vec{\V}\rest \a$ has the $Z$-realizability property and whenever $g\subseteq Coll(\omega, <\mu)$ is generic, in $V[g]$, whenever $\pi:\W\rightarrow \mathcal{C}_k(\V_\eta)$ and $\T$ are such that 
\begin{enumerate}
\item $\pi$ is a weak $k$-embedding 
\item $Z\subseteq rng(\pi)$
\item $\W, \T\in HC$,
\item $\T$ is a correctly guided $Z$-realizable\footnote{See \rdef{correctly guided iterations} and \rdef{dfn:realizable_ext}.} $k$-maximal iteration of $\W$ that is above $\d^{\W^b}$
\end{enumerate}
one of the following holds (in $V[g]$).
\begin{enumerate}
\item $\T$ is of limit length and there is a cofinal well-founded branch $c$ such that $c$ has no drops in model (i.e. $D^\T\cap b = \emptyset$); letting $l=deg^\T(b)$, there is a weak $l$-embedding $\tau:\M^\T_c\rightarrow \mathcal{C}_l(\V_\eta)$ such that $\pi\rest \W=\tau\circ \pi^\T_c$.
\item $\T$ is of limit length and there is a cofinal well-founded branch $c$ such that $c$ has a drop in model, and there is $\b<\eta$ and a weak $l$-embedding $\tau:\M^\T_c\rightarrow \mathcal{C}_l(\V_\b)$ such that $\tau\rest (\M^\T_c)^b=\pi\rest (\M^\T_c)^b$, where $l=deg^\T(c)$.
\item $\T$ has a last model and letting $\gg=lh(\T)-1$, $[0, \gg]_\T\cap \mathcal{D}^\T=\emptyset$ and there is a weak $l$-embedding $\tau:\M^\T_\gg\rightarrow \mathcal{C}_l(\V_\eta)$ such that $\pi\rest \W=\tau\circ \pi^\T$, where $l = deg^\T(\gamma)$.
\item  $\T$ has a last model and letting $\gg=lh(\T)-1$, $[0, \gg]_\T\cap \mathcal{D}^\T\not =\emptyset$ and for some $\b<\eta$ there is a weak $l$-embedding $\tau:\M^\T_\gg\rightarrow \mathcal{C}_l(\V_\b)$ such that $\tau\rest (\M^\T_\gg)^b=\pi\rest (\M^\T_\gg)^b$, where $l = deg^\T(\gamma)$.
\end{enumerate}
When the above 4 clauses hold we say that $\T$ is $(\pi, \vec{\V})$-realizable.
\end{definition}

In the following, we will follow the convention in \cite[Section 1.3]{normalization_comparison}, a (hod, hybrid, or pure extender) premouse has the form $(\M,k)$, where $\M$ is a $k$-sound, acceptable $J$-structure. $k(\M) = k$ is the degree of soundness of $\M$. We write the core $\mathcal{C}(\M)$ for the ($k(\M)+1$-)core of $\M$ (if this makes sense, i.e. when $\M$ is $k(\M)+1$-solid). Similarly, we write $\rho(\M)$ for the $k(\M)+1$-projectum and $p(\M)$ for the $k(\M)+1$-standard parameters of $\M$. When $\mathcal{C}(\M)$ exists, $k(\mathcal{C}(\M)) = k(\M)+1$. $\M$ is sound iff $\M = \mathcal{C}(\M)$. We allow our iterations (e.g. $Z$-validated iterations) to consist of stacks of normal trees, where we may drop gratuitously at the start of a tree. 

\begin{proposition}\label{honest realizable} Suppose $\vec{\V}$ is an array with $Z$-realizability property. Assume further that $p$ is a $Z$-validated iteration of $\mathcal{C}_n(\V_\eta)$ (for some $n$) with last model $\R$ such that $\pi^p$ exists and all the generators of $p$ are contained in $\d^{\R^b}$. Suppose $U$ is a good hull such that $(\vec{\V}, \R, p)\in U$. Let $\R_U = \pi_U^{-1}(\R), p_U = \pi^{-1}(p), \W = \pi^{-1}(\mathcal{C}(\V_\eta))$. There is then a weak $n$-embedding $k:\R_U\rightarrow \mathcal{C}_n(\V_\eta)$ such that $\pi_U\rest \W=k\circ \pi^{p_U}$.
\end{proposition}
\begin{proof} As $p$ is $Z$-realizable, letting $X=U\cap \mH^-$, we can find a $Y$ extending $Z\oplus X$ such that $\Q^Z_Y=\R^b_U$ and $\tau^Z_X=\tau^Z_Y\circ \pi^{p_U, b}$. Let $E$ be the $(\d^{\R_U^b}, \d^{\mH})$-extender derived from $\tau^Z_Y$ and $i:Ult(\R_U, E)\rightarrow \mathcal{C}_n(\V_\eta)$ be the factor map given by $i(\pi^p(f)(a))=\pi_U(f)(\tau^Z_Y(a))$. It then follows that $i\circ \pi_E$ is as desired. 
\end{proof}

Next we introduce a weak notion of realizability. 
\begin{definition}\label{def embeddable}
Suppose $\vec{\V}$ is an array of length $\eta$ that has the $Z$-realizability property and $p$ is a $Z$-validated iteration of $\mathcal{C}_n(\V_\eta)$ with last model $\R$ such that $\pi^p$ exists and the generators of $p$ are contained in $\d^{\R^b}$. Suppose $\T$ is an $\sf{nuvs}$ iteration of $\R$ that is above $\d^{\R^b}$. We say $b$ is  $(Z, \vec{\V})$-embeddable branch of $\T$ if whenever $\l\in S$ is such that $(\R, \T, \vec{\V})\in V_\l$ and $U$ is a $(\mu, \l,Z)$-good hull with $(\vec{\V}, \R, \T, b)\in U$, there is $\a\leq lh(\vec{\V})$, some $l$, and a weak $l$-embedding $k:\M^{\T_U}_{b_U}\rightarrow \mathcal{C}_l(\V_\a)$.
\end{definition}

\begin{proposition}\label{z-validated are embeddable}
Suppose $\vec{\V}$ is a small array with the $Z$-realizability property. Set $\eta=lh(\vec{\V})$. Suppose further that $p$ is a $Z$-validated iteration of $\mathcal{C}_n(\V_\eta)$ with last model $\R$ such that $\pi^{p}$ exists and the generators of $p$ are contained in $\d^{\R^b}$. Additionally, suppose  that $\T$ is an iteration of $\R$ above $\d^{\R^b}$ such that $p^\frown \T$ is $Z$-validated  iteration of $\V_\eta$. Then for all limit $\a<lh(\T)$, if $\T\rest \a$ is $\sf{nuvs}$ then $[0, \a]_\T$ is the unique branch $c$ of $\T\rest \a$ such that $\Q(c, \T\rest \a)$ exists and is $(Z, \vec{\V})$-embeddable. 
\end{proposition}
\begin{proof} We first show that $c=_{def}[0, \a]_\T$ is $(Z, \vec{\V})$-embeddable. Towards contradiction assume not, and suppose $\a$ is least such that $\T\rest \a$ is $\sf{nuvs}$ but $[0, \a]_\T$ is not $(Z, \vec{\V})$-embeddable. Let $U$ be a $(\mu, \l, Z)$-good hull such that $(\R, \vec{\V},p, \T, \a)\in U$. Let $\V'=\pi_U^{-1}(\mathcal{C}_n(\V_\eta))$ and $k:\R_U\rightarrow \mathcal{C}_n(\V_\eta)$ be such that $\pi_U\rest \V'=k\circ \pi^{p_U}$.

We now have a cofinal branch $d$ of $\T_U\rest \a_U$ such that for some $\b\leq\eta$ there is $m:\M^{\T_U\rest \a_U}_d\rightarrow \V_\b$ and $\Q(d, \T_U\rest \a_U)$ exists\footnote{This is a consequence of the fact that $\vec{\V}$ is small.}. Let $\M=\Q(d, \T_U\rest \a_U)$ and $\N=\Q(c_U, \T_U\rest \a_U)$. Both $\M$ and $\N$ are $Z$-approved. Let $\S_0=\m^+(\T_U\rest \a_U)$. If we could conclude that $\M=\N$ then we would get that $c_U=d$, and that would finish the proof. To conclude that $\M=\N$, we need to argue that $\S_0$ is not infinitely descending\footnote{See \rdef{inf desc}.}. The reader may wish to review \rdef{inf desc} and the discussion after \rprop{branch uniqueness}\footnote{This discussion shows that $Z$-approved $\Q$-structures are the same provided they are  based on a $\#$-type lsa hod premouse which is not infinitely descending.}.\\

\textit{Claim.} Suppose $\S_0$ is infinitely descending. Then there is a sequence $(p_i, \S_i: i<\omega)$ witnessing that $\S_0$ is infinitely descending such that for some $\b'<\eta$ and for some $i_0<\omega$ for every $i<j\in (i_0, \omega)$ there are weak $n_i$-embeddings $m_i:\S_i\rightarrow \mathcal{C}_{n_i}(\V_{\b'})$ such that $m_i=m_j\circ \pi^{p_i}$. 
\begin{proof} Set $m_0=m$, $\S_0=\S$ and $\b_0=\beta$. We build the sequence by induction. As the successive steps of the induction are the same as the first step, we only do the first step. Let $(p_i', \S_i': i<\omega)$ be any sequence witnessing that $\S_0$ is infinitely descending. We now have two cases. Suppose first that there is $\b_1\leq \b_0$ and a weak $k$-embedding $m_1:\S_1'\rightarrow \mathcal{C}_k(\V_{\b_1})$ such that if $\b_1=\b_0$ then $m_0=m_1\circ \pi^{p_0'}$. In this case, set $p_0=p_0'$ and $\S_1=\S_1'$. Notice that $\S_1$ is infinitely descending. Suppose next that there is no such pair $(\b_1, m_1)$. In this case we have $d_1, \b_1, m_1, n_1$ such that
\begin{enumerate}
\item $\b_1\leq \b_0$,
\item $d_1$ is a maximal branch of $p_1'\rest \epsilon$ for some $\epsilon<lh(p_1')$,
\item $m_1:\M^{p_1'\rest \epsilon}_{d_1}\rightarrow \mathcal{C}_{n_1}(\V_{\b_1})$,
\item if $\b_1=\b_0$ then $m_0=m_1\circ \pi^{p_1'\rest \epsilon}_{d_1}$.
\end{enumerate}
In this case, set $p_1=p_1'\rest \epsilon^\frown \{d_1\}$ and $\S_1=\M^{p_1'\rest \epsilon}_{d_1}$, with $\b_1$ and $m_1$ as above. We now claim that $\S_1$ is still infinitely descending. To see this, let $c=[0, \epsilon)_{p'_1}$. Notice that we must have that $\Q(c, p'_1\rest \epsilon)\not =\Q(d_1, p'_1\rest \epsilon)$. As both are $Z$-approved, we must have that $\S_1$ is infinitely descending. Continuing in this manner, we get the sequence we desire.
\end{proof}

The existence of a sequence as in the claim above gives us a contradiction, as the sequence must have a well-founded branch. The uniqueness proof is similar to the proof of the claim above, and we leave it to our reader.
\end{proof}

\begin{remark}
If the iteration $p$ in Propositions \ref{honest realizable} and \ref{z-validated are embeddable} drops, we can still embed $\R_U$ by some map $k: \R_U \rightarrow \mathcal{C}_l(\V_\alpha)$ for some $\alpha < \epsilon$ and some $l<\omega$. In this case, there is some model $\M\in p$ such that $\pi^{p_U}_{\M_U,\R_U}$ exists and there is a weak $l$-embedding $\sigma: \M_U\rightarrow \mathcal{C}_l(\V_\alpha)$ such that $\sigma = k\circ \pi^{p_U}_{\M_U,\R_U}$.
\end{remark}

\begin{remark} The fact that $\V$ is small is key to the proof of \rprop{z-validated are embeddable}. See \rprop{break4} which partially deals with the situation when $\V$ is not small.
\end{remark}

Motivated by \rprop{z-validated are embeddable}, we make the following definitions.

\begin{definition}\label{weakly suitable above mu} We say $\R$ is \textbf{weakly $Z$-suitable above $\mu$} if $\R$ is a hod premouse of lsa type such that $\R=(\R|\d^\R)^{\#}$ and whenever $\l\in S$ is such that $\R\in V_\l$ and $U$ is a $(\mu, \l, Z)$-good hull, $\R_U$ is weakly $Z$-suitable\footnote{See \rdef{weakly suitable}.}. 
\end{definition} 

\begin{definition}\label{honest wzc} Suppose $\R$ is weakly $Z$-suitable above $\mu$. We say $\R$ is \textbf{honest} if there is an array $\vec{\V}=(\V_\a:\a\leq \eta)$ at $\mu$ with the $Z$-realizability property such that letting $\l\in S$ be the least such that $\R, \vec{\V}\in V_\l$, the following conditions hold. 
\begin{enumerate}
\item Either $\V_\eta=\R$ or there is a $Z$-validated iteration $p$ of $\V_\eta$ of limit length such that $\pi^{p, b}$ exists and $\R=\m^+(p)$.
\item $\vec{\V}$ is small if and only if  $\V_\eta \not=\R$.
\end{enumerate}
If $\R$ is honest and $\vec{\V}$ is as above then we say that $\vec{\V}$ is an honesty certificate for $\R$. 
\end{definition}

Suppose $\R$ is honest as witnessed by $(\vec{\V}, p)$. Then we say $\T$ is a $Z$-validated iteration of $\R$ if $p^\frown \T$ is a $Z$-validated iteration of $\V_\eta$ where $\eta+1=lh(\vec{\V})$. 

\begin{proposition}\label{bottom part realizability}
Suppose $\R$ is weakly $Z$-suitable above $\mu$, $(\vec{\V}, p)$ is an honesty witness, and suppose $\T$ is a $Z$-validated $\sf{nuvs}$ iteration of $\R$ with last model $\S$ such that $\pi^{\T}$ exists and the generators of $\T$ are contained in $\d^{\S^b}$. Suppose $U$ is a good hull such that $(\R, \vec{\V}, p, \T, \S)\in U$. There is then $\a\leq lh(\vec{\V})$, a $Z$-approved sts mouse $\M$ over $\R_U$, an embedding $k:\M\rightarrow \mathcal{C}(\V_\a)$ and an embedding $\sigma:\S_U\rightarrow \mathcal{C}(\V_\a)$ such that 
\begin{enumerate}
\item $\M\models ``\d^\R$ is a Woodin cardinal",
\item $k\rest \R_U=\sigma\circ \pi^{\T_U}$, 
\item $\M\not =\R_U$ if and only if $\R\not =\mathcal{C}(\V_\eta)$ (so $\vec{\V}$ is small), and
\item if $\M\not=\R_U$ then $rud(\M)\models ``\d^\R$ is not a Woodin cardinal".
\end{enumerate}
\end{proposition}
\begin{proof} First we claim that there is $\a\leq lh(\vec{\V})$, an $l<\omega$, and a weak $l$-embedding $k:\R_U\rightarrow \mathcal{C}_l(\V_\a)$ such that $Z\subseteq rng(k)$. If $\R=\V_\eta$ then set $\a=\eta$ and $k=\pi_U\rest \R_U$.

Suppose then $\R\not=\V_\eta$. In this case, $\vec{\V}$ is small. Let $\W$ be the largest node on $p$ such that $\pi^{p_{\leq \W}}$ exists and the generators of $p_{\leq \W}$ are contained in $\d^{\W^b}$. Then $p_{\geq \W}$ is above $\d^{\W^b}$. It follows from \rprop{z-validated are embeddable} and the remark after that $p_{\geq \W}$ is according to $(Z, \vec{\V})$-embeddable branches, and therefore we must have $\a\leq lh(\vec{\V})$ and a cofinal branch $c$ of $p_U$ such that there is an appropriate weak $l$-embedding $k:\M^{p_U}_c\rightarrow \mathcal{C}_l(\V_\a)$\footnote{$l$ is specified as in Definition \ref{realizability property}.} such that $Z\subseteq rng(k)$. Set then $\M=\M^{p_U}_c$; note that $\R_U \lhd \M$ and $rud(\M) \vDash ``\delta^\R$ is not Woodin" by smallness of $\vec{\V}$.

We continue with one such pair $(\a, k)$. Next, as $\T$ is $Z$-validated, we must have $Y$ an extension of $Z$ such that $X=_{def}k[\R_U^b]\subseteq Z\oplus Y$, $\S_U^b=\Q^Z_Y$ and $\tau^Z_X=\tau^Z_Y\circ \pi^{\T_U, b}$.  We can then lift $\tau^Z_Y$ to $\S$ and obtain some weak $l$-embedding $\sigma:\S_U\rightarrow \mathcal{C}_l(\V_\a)$ such that $k=\sigma\circ \pi^{\T_U}$. $\sigma$ is essentially the ultrapower map by the $(\d^{\S_U^b}, \d^{\mH})$-extender derived from $\tau^Z_Y$ (see the proof of Proposition \ref{honest realizable}).
\end{proof}

Finally we discuss iterations that are above $\S^b$ where $\S$ is as in \rprop{bottom part realizability}. The proof is just like the proof of \rprop{z-validated are embeddable}.

\begin{proposition}\label{rz-validated iterations} Suppose $\R$ is honest weakly $Z$-suitable above $\mu$ hod premouse and $(\vec{\V}, p)$ is an honesty witness for $\R$.  Suppose $\T$ is a normal $Z$-validated iteration of $\R$. Let $\a'\in R^\T$ be the largest such that setting $\S=\M_\a^\T$, $\pi^{\T_{\leq \S}}$ exists and the generators of $\T_{\leq \S}$ are contained in $\S^b$. Suppose $\T_{\geq \S}$ is above $Ord\cap \S^b$. Let $U$ be a good hull such that $\{\R,\vec{\V}, \T\}\in U$, and let $(\a, \M, k, \sigma)$ be as in \rprop{bottom part realizability}. 
Then $\T_{\geq \S}$ is $(\sigma, \vec{\V}\rest \a)$-realizable. Moreover, for each limit ordinal $\b<lh(\T_{\geq \S})$, if $\T_{\geq \S}\rest \b$ is $\sf{nuvs}$ then $d=_{def}[0, \b]_{\T_{\geq \S}}$ is the unqiue cofinal branch $d'$ of $\T_{\geq \S}\rest \b$ which is $(Z, \vec{\V})$-embeddable and $\Q(d', \T_{\geq \S}\rest \b)$ exists. 
\end{proposition}  

\begin{definition}\label{z-full} We say $\R$ is \textbf{$Z$-suitable above $\mu$} if it is weakly $Z$-suitable above $\mu$ and whenever $\M$ is a $Z$-validated sts mouse over $\R$, $\M\models ``\d^\R$ is a Woodin cardinal". 
\end{definition}

Our goal is to construct an $\R$ which is \textbf{$Z$-suitable above $\mu$}. We do this in \rprop{getting suitable}.

\section{$Z$-validated sts constructions}

We assume that the theory $T$ holds (see \rdef{theory t}). We then fix $(S, S_0, \nu_0,\vec{Y}, \vec{A})$ that witness $T$ and let $\mu\in S_0$. We will omit $\mu$ when discussing CMI objects at $\mu$.

The construction that we will perform in the next section will hand us a hod premouse $\R$ that is weakly $Z$-suitable above $\mu$. The rest of the construction that we will perform will be a fully backgrounded construction over $\R$ whose aim is to either find a $Z$-validated sts mouse destroying the Woodiness of $\d^\R$ or proving that no such structures exist. In the latter case, we will show that we must produce an excellent hybrid premouse. 

The construction that we describe in this section is a construction that is searching for the $Z$-validated sts mouse over $\R$ destroying the Woodiness of $\d^\R$. In this construction, we add two kinds of objects. The first type of objects are extenders, and they are handled exactly the same way that they are handled in all fully backgrounded constructions. The second kind of objects are iterations. Here the difference with the ordinary is that there is no strategy that we follow as we index branches of iterations that appear in the construction. Instead, when our sts scheme demands that a branch of some iteration $p$ must be indexed, we find an appropriate branch and index it. We will make sure that the iterations that we need to consider in the construction are all $Z$-validated. It must then be proved that given a $Z$-validated iteration there is always a branch that is $Z$-validated.

The solution here has a somewhat magical component to it. As we said above, the fully backgrounded $Z$-validated sts construction is not a construction relative to a strategy. This is an important point that will be useful to keep in mind. Instead, the construction follows the sts scheme, and the $Z$-validation method is used to find branches of iterations that come up in the construction. To see that we do not run into trouble we need to show that any such iteration $\T$ that needs to be indexed according to our sts scheme has a branch $b$ such that $\T^\frown \{b\}$ is $Z$-validated. Let $\M$ be the stage of the construction where $\T$ is produced. Recall now that we have two types of such iterations. If $\T$ is $\sf{uvs}$ then $Z$-validation will produce a branch in a more or less straightforward fashion (see \rprop{lem:break2a}). If $\T$ is $\sf{nuvs}$ then the fact that we need to index a branch of it suggests that we have also reached an authenticated $\Q$-structure for $\T$. We will then show that there must be a branch with this $\Q$-structure. This is the magical component we speak of above. In general, given an iteration $\T$ of a weakly $Z$-suitable $\R$ that is produced by $\sf{HFBC}(\mu)$ of the next section, there is no reason to believe that there is a $\Q$-structure for it of any kind. Even if there is a $\Q$-structure $\Q$ of some kind, there is no reason to believe that sufficiently  closed hulls of $\T$ will have branches determined by the pre-image of $\Q$. In our case, what helps is that $\Q\in \M$, and this condition, in the authors' opinion, is somewhat magical. \footnote{The proof of this fact is in Section \ref{sec:all_together}; it shows that if a level $\S$ of the construction has no Woodin cardinals, then if $\T$ is a tree on $\S$, then sufficently closed hulls of $\T$ will have branches determined by the pre-image of $\Q(\T)$.}

One particularly unpleasant problem is that we cannot in general prove that the non weakly $Z$-suitable levels of the fully backgrounded construction produced in the next sections are iterable. This unpleasantness causes us to work with weakly $Z$-suitable $\R$ that are iterates of a level of the fully backgrounded construction of the next section. In order to have an abstract exposition of the $Z$-validated sts construction, we introduce the concept of \textit{honest weakly $Z$-suitable} $\R$ over which we will perform our $Z$-validated sts constructions. The honest weakly $Z$-suitable hod premice will have honesty witnesses, and that is the concept we introduce first. The honesty witnesses are essentially models of a $K^c$-construction.

In this section and subsequent sections, we work with the fine structure in \cite{steel2010outline}.


%
%

\subsection{The $Z$-validated sts construction}

Suppose $\R$ is honest and $\vec{\V}$ is an honesty certificate for $\R$. We assume that $\R$ is a $\#$-lsa type hod premouse. Let $X$ be any transitive set such that $\R\in X$. Let $\l \in S_0$ be such that $\R, X\in H_\l$. In what follows we introduce the fully backgrounded $(Z, \l)$-validated sts construction over $X$.  
\begin{definition}\label{dfn:Kc_constr} We say $(\M_\xi, \N_\xi: \xi\leq \Omega^*)$ are the models of the \textit{fully backgrounded $(Z, \l)$-validated sts construction over $X$} if the following conditions hold:
\begin{enumerate}
\item $\Omega^*\leq \l$, for all $\xi<\l$ if $\M_\xi, \N_\xi$ are defined then $\M_\xi$ and $\N_\xi\in H_\l$.
\item For every $\xi\leq \Omega^*$, $\M_\xi$ and $\N_\xi$ are $Z$-validated sts hod premice over $X$.
\item Suppose the sequence $(\M_\xi, \N_\xi: \xi<\eta)$ and $\M_\eta$ have been constructed. Suppose further that
there is a total $(\k, \nu)$-extender $F$ such that letting $G=\M_\eta\cap F$, $(\M_\eta, G)$ is a $Z$-validated sts hod premouse over $X$. Let then $\N_\eta=(\M_\eta, G)$ and 
$\M_{\eta+1}=\mathcal{C}(\N_\eta)$\footnote{$\mathcal{C}(\N_\eta)$ is the core of $\N_\eta$.}.
\item Suppose the sequence $(\M_\xi, \N_\xi: \xi<\eta)$ and $\M_\eta$ have been constructed, and $\T\in \M_\eta$ is the $<_{\M_\eta}$-least $\sf{uvs}$ tree\footnote{$<_{\M_\eta}$ is the canonical well-ordering of $\M_\eta$.} without an indexed branch. Suppose further that there is a branch $b$ of $\T$ such that $(\M_\eta, b)$ is a $Z$-validated sts hod premouse\footnote{This in particular implies that $b\in \M_\eta$.} over $\R$. Let then $\N_\eta=(\M_\eta, b)$ and $\M_{\eta+1}=\mathcal{C}(\N_\eta)$. 
\item Suppose the sequence $(\M_\xi, \N_\xi: \xi<\eta)$ and $\M_\eta$ has been constructed, and for some $\sf{nuvs}$ tree $\T\in \M_\eta$ there is a branch $b\in \M_\eta$ such that $(\M_\eta, b)$ is a $Z$-validated sts hod premouse over $\R$. Let $\T$ be the $<_{\M_\eta}$-least such tree and $b$ be such a branch for $\T$. Then $\N_\eta=(\M_\eta, b)$ and $\M_{\eta+1}=\mathcal{C}(\N_\eta)$\footnote{Recall the Internal Definability of Authentication. In this case of the construction, the branch $b$ is chosen by a procedure internal to $\M_\eta$ and does not depend on any external factors. Because of this, proving $Z$-validity is not obvious at all. Also see the Anomaly in 3.b of \cite[Definition 4.2.1]{hod_mice_LSA}.}.
\item Suppose the sequence $(\M_\xi, \N_\xi: \xi<\eta)$ and $\M_\eta$ has been constructed and all of the above cases fail. In this case we let $\N_\eta=\mathcal{J}_1(\M_\eta)$ and provided $\N_\eta$ is a $Z$-validated sts hod premouse over $\R$, $\M_{\eta+1}=\mathcal{C}(\N_\eta)$.
\item Suppose the sequence $(\M_\xi, \N_\xi: \xi<\eta)$ has been constructed and $\eta$ is a limit ordinal. Then $\M_\eta=liminf_{\xi\rightarrow \eta}\M_\xi$.
\end{enumerate}
If $\lambda$ is clear from context then we will omit it from our notation.
\end{definition}

The fully backgrounded (f.b.) $Z$-validated sts construction can break down for several reasons. Let $p$ be the $\V_\eta$-to-$\R$ iteration witnessing that $\R$ is honest. Below we list all of these reasons. We say that f.b. $Z$-validated sts construction breaks down at $\eta$ if one of the following conditions holds.\\\\
\textbf{Break1.} $\M_\eta$ is not solid or universal.\\
\textbf{Break2.} $\M_\eta$ is not $Z$-validated\footnote{It is not hard to see that if $\M_\eta$ is $Z$-validated and $\N_\eta=\mathcal{J}_1(\M_\eta)$ then $\N_\eta$ is $Z$-validated. It is possible that $\M_\eta$ is $Z$-validated but $\N_\eta$ is not but these possibilities are covered by Break3 and Break4. Also notice that if $\N_\eta=(\M_\eta, G)$ where $G$ is an extender and $\M_\eta$ is $Z$-validated then $\N_\eta$ is $Z$-validated.}.\\
\textbf{Break3.} There is an $\sf{uvs}$ tree $\T\in \M_\eta$ such that the indexing scheme demands that a branch of $\T$ must be indexed yet $\T$ has no (cofinal well-founded) branch $b$ such that $(\M_\eta, b)$ is a $Z$-validated sts premouse over $\R$. \\
\textbf{Break4.} There is an $\sf{nuvs}$ tree $\T\in \M_\eta$ such that the indexing scheme demands that a branch of $\T$ must be indexed yet letting $b\in \M_\eta$ be the branch given by the authentication process,  $b$ is not $(Z, \vec{\V})$-embeddable branch of $\T^\frown p$.\footnote{See \rdef{def embeddable}. Break4 is similar to the Anomaly in 3.b of \cite[Definition 4.2.1]{hod_mice_LSA}.}\\
\textbf{Break5.} $\rho(\M_\eta)\leq \d^{\mH}$.\\

The argument that the construction doesn't break down because of Break1 is standard, cf. \cite{FSIT}. It is essentially enough to show that the countable substructures of $\M_\eta$ are iterable. We will show that much more complicated forms of iterability hold, and so to save ink and to not repeat ourselves, we will leave this portion to our kind reader.  To see that the construction doesn't break down because of Break2 is not too involved, and we will present that argument below. At this point, we cannot do much about Break5. We will deal with it when $X$ becomes a more meaningful object. The remaining cases will be handled in the next subsections.

\begin{proposition}\label{break2} Suppose $\R$ is an honest weakly $Z$-suitable hod premouse as witnessed by $\vec{\V}$, $X$ is a transitive set such that $\R\in X$ and $\l\in S_0$ is such that $X, \vec{\V}\in V_\l$.  Then the f.b. $(Z, \l)$-validated construction over $X$ does not break down because of Break2. 
\end{proposition}
\begin{proof} Towards contradiction assume that there is some $\a$ such that for all $\b<\a$, both $\M_\b$ and $\N_\b$ are $Z$-validated but $\M_\a$ is not $Z$-validated. Set $\W=\M_\a$ be the least such model. 

Suppose first $\a$ is a limit ordinal. Let $U$ be a $(\mu, \l, Z)$-good hull such that $ \{\R,  \W\}\subseteq U$ and $(\M_\b:\b\leq \a)\in U$. Let $(\K_\xi: \xi\leq \a_U)=\pi^{-1}_U(\M_\b:\b\leq \a)$. Fix $\T\in \K_{\a_U}$ according to $S^{\K_{\a_U}}$.  We need to see that $\T$ is $Z$-approved. Fix $\xi<\a_U$ such that $\T\in \K_\xi$ and is according to $S^{\K_\xi}$. Then $\pi_U(\T)\in \M_{\pi_U(\xi)}$ and is according to $S^{\M_{\pi_U(\xi)}}$. Therefore, $\T$ is $Z$-approved. 

Suppose next that $\a=\b+1$. Because we are assuming the least model that is not $Z$-validated is $\M_\a$ we must have that $\N_\b$ is $Z$-validated. Let now $U$ be a $(\mu, \l, Z)$-good hull such that $\{\R, \W\}\subseteq U$. But then $\pi_U^{-1}(\M_\a)=\mathcal{C}(\pi^{-1}_U(\N_\b))$. It then follows that $\pi_U^{-1}(\M_\a)$ is $Z$-approved (see \rprop{preservation of zv under embeddings}).
\end{proof}

\subsection{Break3 never happens} 

In this subsection, $\R$ is an honest weakly $Z$-suitable, $X$ is a transitive set containing $\R$ and $\l\in S_0$ is such that $(\R, X)\in V_\l$. Our main goal here is to prove that the $(Z, \l)$-validated sts construction over $X$ doesn't break down because of Break3. First, we prove the following general lemma. 

\begin{lemma}\label{general lemma on cof} Suppose $\M$ is a hod premouse for which $\M^b$ is defined, and $p$ is an iteration of $\M$ such that $\pi^{p, b}$ exists. Let $\d$ be a Woodin cardinal of $\pi^{p, b}(\M^b)$ and let $\xi$ be least such that $\pi^{p, b}(\xi)\geq \d$. Then $\cf(\d)=\cf((\xi^+)^\M)$. 
\end{lemma}
\begin{proof} Let $\Q$ be the least model of $p$ such that $\Q^b=\pi^{p, b}(\M^b)$ and set $q=p_{\leq \Q}$. Let $\N$ be the least model on $q$ such that $\d\in rng(\pi^q_{\N, \Q})$. Without losing generality we can assume $\Q=\N$ as $rng(\pi^q_{\N, \Q})\cap \d$ is cofinal in $\d$. As the iteration embeddings are cofinal at Woodin cardinals if $\pi^q(\xi)=\d$ then again there is nothing to prove. Assume then $\pi^q(\xi)>\d$. Without loss of generality we can assume that $\xi=\d^{\M^b}$. If $\xi<\d^{\M^b}$ then we need to redefine $\M$ as $\M|\zeta$ where $\zeta$ is the $\M$-successor of $o^\M(\xi)$. 

Because $\N$ is the least model that has $\d$ in it, it must be case that $\N=Ult(\W, E)$ where $\W$ is a node in $q$ and $E$ is an extender used in $q$ to obtain $\N$. Moreover, $\cp(E)=\d^{\W^b}$. Below $\pi_E$ is used for $\pi_E^\W$.

Suppose $\pi_E(f)(a)=\delta$ and $\pi_E(g)(a)= \nu_E$, where 
\begin{enumerate}
\item $\nu_E$ is the supremum of the generators of $E$
\item $f,g: \d^{\W^b} \rightarrow \d^{\W^b}$ are functions in $\W$, 
\item $a \in [\nu_E+1]^{<\omega}$.
\end{enumerate}
 Note that $\nu_E < \delta$.

We first show that \\\\
(1) $\sup(\{ \pi_E(k)(a): k:\d^{\W^b}\rightarrow \d^{\W^b}, k\in \W\}\cap \d)=\d$.\\\\
To see (1) fix $h:\d^{\W^b} \rightarrow \d^{\W^b}$ in $\W$ and let $s$ in $[\nu_E+1]^{<\omega}$ be such that $\pi_E(h)(s) < \delta$. We want to find $k$ such $\pi_E(k)(a)$ is in $[\pi_E(h)(s), \delta].$ Set $k(u)=$ the supremum of points of the form $h(t)$ such that $h(t)< f(u)$ and $t$ is a finite sequence from $g(u)$.  $f(u)$ is a Woodin cardinal (in $\R$), so $k(u)< f(u)$ for $E_a$-almost all $u$, so 
\begin{center}
$\pi_E(k)(a) < \delta = \pi_E(f)(a)$. 
\end{center}
Also 
\begin{center}
$\pi_E(h)(s)\leq \pi_E(k)(a)$
\end{center}
by the definition of $k$.

Let $\l=Ord\cap \W^b$. We have that $\cf(\l)=\cf(Ord\cap \M^b)$. Thus, it is enough to show that $\cf(\d)=\cf(\l)$. Let $\eta=\cf(\d)$ and let $(k_\a: \a<\eta)\subseteq \W$ be such that 
\begin{enumerate}
\item for $\a<\eta$, $k_\a:\d^{\W^b} \rightarrow \d^{\W^b}$,
\item for $\a<\eta$, $k_\a\in \W$,
\item for $\a<\b<\eta$, $\pi_E(k_\a)(a)<\pi_E(k_\b)(a)<\d$.
\end{enumerate}
Let $\vec{\gg}=(\gg_\a:\a<\eta)$ be increasing and such that
\begin{enumerate}
\item $k_\a\in \W|\gg_\a$,
\item $\rho(\W|\gg_\a)=\d^{\W^b}$.
\end{enumerate}
We claim that $\vec{\gg}$ is cofinal in $\l$. Suppose it is not. In that case, we can fix $\zeta>\sup\vec{\gg}$ such that $\rho_1(\W|\zeta)=\d^{\W^b}$. Let $p$ be the first standard parameter of $\W|\zeta$. For each $\a<\eta$, let $a_\a\in [\d^{\W^b}]^{<\omega}$ be such that $k_\a$ is definable from $p$ and $a_\a$ in $\W|\zeta$. It then follows that 
\begin{center}
$\sup(Hull_1^{\pi_E(\W|\zeta)}(\pi_E(p), \d^{\W^b})\cap \d)=\d$,
\end{center}
as witnessed by $(a_\a: \a<\eta)$. As $Hull_1^{\pi_E(\W|\zeta)}(\pi_E(p), \d^{\W^b})\in Ult(\W, E)=\N$, the above equality implies that $\d$ is singular in $\N$, contradiction. Thus, $\vec{\gg}$ must be cofinal in $\l$. Therefore, $\cf(\l)=\eta$.
\end{proof}

 Recall that we are working under theory $T$, see \rdef{theory t}.

\begin{corollary}\label{lem:small_cofa} 
Suppose $\T$ is a normal tree on $\R$ such that $\pi^{\T, b}$ exists and $\delta$ is a Woodin cardinal of $\pi^{\T, b}(\mH)$. Then $\cf(\delta) < \mu$ and if $\d>\sup(\pi^{\T, b}[\d^{\mH}])$ then $\cf(\d)<\nu_0$.
\end{corollary}
\begin{proof} First note that if $\delta$ is a Woodin cardinal of $\mathcal{H}$, then $\cf(\delta) < \mu$. This is because there is a hod pair $(\P,\Sigma)\in \mathcal{F}$, a $\delta^*$ such that $\P\models \delta^*$ is Woodin and $\delta = \pi_{\P,\infty}(\delta^*)$. 
Now, if $\d>\sup(\pi^{\T, b}[\d^{\mH}])$ then by \rlem{general lemma on cof}, $\cf(\d)=\cf(Ord\cap \mH)<\nu_0$.
\end{proof}

We now state and prove our main proposition of this subsection.

\begin{proposition}\label{lem:break2a} The $(Z, \l)$-validated sts construction over $X$ doesn't break down because of Break3. 
\end{proposition}
\begin{proof}  \cite[Section 12]{hod_mice_LSA} handles a similar situation, and the proof here is very much like the proofs in  \cite[Section 12]{hod_mice_LSA}. Because of this we give an outline of the proof.

Suppose $\M$ is a model appearing in the $(Z, \l)$-validated sts construction over $X$ and $\T^*\in \M$ is an $\sf{uvs}$ iteration of  $\R$ such that the indexing scheme requires that we index a branch of $\T^*$ at $Ord\cap \M$. We need to show that there is a branch $b$ of $\T^*$ such that $(\M, b)$ is $Z$-validated. Because of \rprop{unique branches above mu} and \rprop{rz-validated iterations}, there can be at most one such branch. 

Because $\T^*$ is $\sf{uvs}$, we have a normal iteration $\T\in \M$ with last model $\S$ such that $\pi^{\T}$ is defined and a normal iteration $\U$ based on $\S^b$ such that $\T^\frown \U=\T^*$. Because the construction doesn't break because of Break2 (see \rprop{break2}), we have that $\M$ is $Z$-validated and therefore, $\T$ is $Z$-validated. Also, we can assume that $\U$ is not based on $\S|\xi$ where $\xi=\sup(\pi^\T[\d^\mH])$, as otherwise the desired branch of $\U$ is given by $\Psi$. 

We now show that $\U$ has a branch $b$ such that $(\M, b)$ is $Z$-validated. Let $\l\in S$ be least such that $\R, \M\in H_\l$. Given a $(\mu, \l, Z)$-good hull $U$ such that $\{\M, \T, \S, \U\}\subseteq U$, let $b_U=\Psi^Z_W(\pi^{-1}_U(\U))$ where $W$ is any extension of $Z$ such that $\pi^{-1}_U(\S^b)=\Q^Z_W$. First we claim that for all $U$ as above, \\

\textit{Claim 1.} $b_U\in M_U$. 
\begin{proof} Given a $U$ as above, we will use it as a subscript to denote the $\pi_U$-preimages of the relevant objects. Fix then a $U$ as above. Suppose first that $\Q(b_U, \U_U)$ doesn't exist. As we are assuming $\U$ is not based on $\S|\xi$, \rcor{lem:small_cofa} implies that $\cf(\d(\U))\leq \nu_0$. Because $M_U$ is $\nu_0$-closed it follows that $b_U\in M_U$.

Suppose next that $\Q(b_U, \U_U)$ exists. Let $A_U$ be the preimage of $A_\l$. Notice now that letting $\Phi$ be the $\pi_U$-pullback of $\Psi_\l$, we have that $Lp^{cuB, \Phi}(A_U)\in M_U$. 

Let $Y=U\cap \mH$. Clearly $Y$ is an extension of $Z$ and because $\M$ is $Z$-validated, we must have $W^*$ an extension of $Z\cup Y$ such that $\S_U^b=\Q^Z_{W^*}$. Notice that because $\Psi^Z_{W^*}$ is computable from $\Phi$ and because $Lp^{cuB, \Phi}(A_U)\in M_U$, we must have that $\Q(b_U, \U_U)\in M_U$. Hence, $b_U\in M_U$.
\end{proof}

Suppose first that $\cf(lh(\U))>\omega$. In this case, let $\U$ be as above and set $c=\pi_U(b_U)$. Then $c$ is the unique well-founded branch of $\U$ and hence, for any $(\mu, \l, Z)$-good hull $X'$ such that  $U\cup \{(\M, c), U\}\in X'$, $c_{X'}=b_{X'}$. Hence, $(\M, c)$ is $Z$-validated (see  \rprop{one hull witness for premice}). 

Suppose then $lh(\U)=\omega$. We now claim that\\

\textit{Claim.} there is a $(\mu, \l, Z)$-good hull\footnote{See \rdef{def_goodhull}.} $X_0$ such that for all $(\mu, \l, Z)$-good hulls $Y$ such that $X_0\cup \{\M, X_0\}\in Y$, $\pi_{X_0, Y}(b_{X'})=b_Y$. 
\begin{proof}
Assuming not we get a  continuous chain $(X_\a: \a<\mu)$ such that
\begin{enumerate}
\item $\M, \U\in X_0$,
\item for all $\a<\mu$, $X_{\a+1}$ is a $(\mu, \l, Z)$-good hull,
\item for all $\a<\mu$, $X_{\a}\cup \{X_{\a}\}\in X_{\a+1}$,
\item for all $\a<\mu$, $\pi_{X_{\a+1}, X_{\a+2}}(b_{X_{\a+1}})\not =b_{X_{\a+1}}$.
\end{enumerate}
Let $\nu\in (\nu_0, \mu)$ be an inaccessible cardinal such that $X_\nu\cap \mu=\nu$. Fix now $\a<\nu$ such that 
\begin{center}
$\sup(b_{X_\nu}\cap  rng(\pi_{X_\a, X_\nu}))=lh(\U_{X_\nu})$.
\end{center}
 As $\cf(lh(\U_{X_\nu}))=\omega$ this is easy to achieve. For $\b\in [\a, \nu)$ let $c_\b$ be the $\pi_{X_\a, X_\nu}$-pullback of $b_{X_\nu}$. Let for $\b\in [\a, \nu]$, $W_\b$ be such that $\S^b_{X_\b}=\Q^Z_{W_\b}$. It follows that $c_\b$ is according to the $\pi_{X_\b, X_\nu}$-pullback of $\Psi^Z_{W_\nu}$. Because $\Psi^Z_{W_\b}$ depends only on $\S^b_{X_\b}$, we have that $c_\b=b_{X_\b}$ (this is because the  $\pi_{X_\b, X_\nu}$-pullback of $\Psi^Z_{W_\nu}$ is a strategy of the form $\Psi^Z_Y$ where $\Q^Z_Y=\S^b_{X_\b}$). It follows that for all $\b<\gg\in [\a, \nu)$, $\pi_{X_\b, X_\gg}(b_{X_\b})=b_{X_\gg}$. 
 \end{proof}
 
 Fix now an $X_0$ as in the Claim. Set $c=\pi_{X_0}(b_{X_0})$. The above property of $X_0$ guarantees that $(\M, c)$ is $Z$-validated. Indeed, fix a $(\mu, \l, Z)$-good hull $U$ such that $\M, c\in U$. Let $Y$ be a $(\mu, \l, Z)$-good hull such that $X_0\cup U \cup \{X_0, U\}\in Y$. Then $\pi_{U, Y}(c_U)=\pi_{X_0, Y}(b_{X_0})=b_Y$. It follows that $c_U$ is the $\pi_{U, Y}$-pullback of $\Psi^Z_W$ where $W$ is such that $\S^b_Y=\Q^Z_W$. Hence, $c_U=b_U$.  
\end{proof}

\subsection{Break4 never happens}

The following is the main proposition of this subsection. We continue with $(\R, X, \vec{\V}, \l)$ of the previous section. 

\begin{proposition}\label{break4} Suppose the $(Z, \l)$-validated sts construction over $X$ breaks because of Break4 and that $X$ is a transitive set such that $H_{\d^\R}^X$ is the universe of $\R|\d^\R$ and $\d^\R$ is a Woodin cardinal in $X$. Then 
\begin{enumerate}
\item $\vec{\V}$ is not small (implying $\V_\eta=\R$), and 
\item letting $\eta'$ be such that Break4 occurs at $\eta'$ and letting $(\T, b)\in \M_{\eta'}$ witness Break4 at $\eta'$, either $\d^\R$ is not a Woodin cardinal of $\M_{\eta'}$ or $H_{\d^\R}^{\M_{\eta'}}$ is not the universe of $\R|\d^\R$. 
\end{enumerate}
\end{proposition}
\begin{proof}  Let $p$ be the $\V_{\eta'}$-to-$\R$ iteration. Setting $\W=\M_{\eta'}$, we have that
\begin{enumerate}
\item $\W$ is $Z$-validated and 
\item $b$ is not $(Z, \vec{\V})$-embeddable (see \rdef{def embeddable}).
\end{enumerate}
Let $\b$ be such that $\W|\b$ authenticates $b$. Thus $\W|\b$ is a model of $\sf{ZFC}$ in which there is a limit of Woodin cardinals $\nu$ and the derived model of $\W|\b$ at $\nu$ has a strategy for $\Q(b, \T)$ that is $\W|\b$-authenticated. \\

\textit{Claim 1.} $p^\frown\T^\frown\{b\}$ is a $Z$-validated iteration (see \rdef{z-validated sts mice above mu}).
\begin{proof}
Towards a contradiction suppose $p^\frown \T^\frown\{b\}$ is not $Z$-validated. Fix now a $(\mu, \l, Z)$-good hull $U$ such that $(\R, \W, p, \T, b) \in U$ and $p_U^\frown \T_U^\frown \{b_U\}$ is not a correctly guided $Z$-realizable iteration of $\R_U$. Because $\W$ is $Z$-validated, we can assume that $p_U^\frown \T_U$ is a correctly guided $Z$-realizable iteration. It must then be that $\Q(b_U, \T_U)$ is not $Z$-approved. 

To save ink, let us prove that in fact $\N=_{def}\Q(b_U, \T_U)$ is $Z$-approved of depth 1. As the proof of depth $n$ is the same, we will leave the rest to the reader. To start with, notice that since $\T_U$ itself is correctly guided $Z$-realizable, we have that $\S=\m^+(\T_U)$ is weakly $Z$-suitable. To prove that $\N$ is $Z$-approved of depth 1 we need to show that if $\U\in \N$ is according to $S^{\N}$ then $\U$ is $Z$-realizable. 

Fix then $\a\in R^\U$ and set $\X=\M_\a^\U$. First let's show that there is $Z'$ an extension of $Z$ such that $\Q^Z_{Z'}=\X^b$. Because $\T_U^\frown \{b_U\}$ is authenticated inside $\W_U|\b_U$, we must have an iteration $\Y$ of $\R_U$ according to $S^{\W_U}$ with last model $\R_1$ such that there is an embedding $k:\X^b\rightarrow \R_1^b$ with the property that $\pi^{\Y, b}=k\circ\pi^{\U_{\leq \X}, b}$. Because $\Y$ is $Z$-realizable, we must have $Y$ an extension of $Z$ such that $\R_1^b=\Q^Z_Y$. Composing $k$ with $\tau^Z_Y$ we have that $\X^b=\Q^Z_{Z'}$ for some $Z'$.

The rest is similar. If $\U^*$ is the longest initial segment of $\U_{\geq \X}$ that is based on $\X^b$ then there are $\Y$ and $k$ as above such that $\U^*$ is according to the $k$-pullback of $S^{\M_U}_{\R_1^b}$. But because $\W_U$ is $Z$-approved, $S^{\W_U}_{\R_1^b}$ is a fragment of $\Psi^Z_Y$ where $Y$ is as above. Hence, $\U^*$ is according to $\Psi^Z_X$ for some $X$ (see \rcor{invariance under pullbacks}).
\end{proof}

Let now $U$ be a $(\mu, \l, Z)$-good hull such that $(\R, \W, \T,  b,\Q)\in U$. Because $\T$ is $Z$-validated, we have that the $\pi_U$-realizable branch $d$ of $\T_U$ is cofinal. Suppose then $\Q(d, \T_U)$ exists. Then because it is $Z$-approved, we must have that $\Q(d, \T_U)=\Q(b_U, \T_U)$ (for example see \rprop{z-validated are embeddable}). It follows that $d=b_U$, and so $b$ is $(Z, \vec{\V})$-embeddable.

Assume now that clause 1 fails. Because $\V$ is not small, we must have that $\Q(d, \T_U)$ exists (as $d$ is the realizable branch of $p^\frown \T_U$). Assume now that $\V$ is not small. This means that $\R=\V_\eta$. Assume now that clause 2 fails. Since $\d^\R$ is a regular cardinal of $\W$ , $H_{\d^\R}^{\W}$ is the universe of $\R|\d^\R$ and $\R=\V_\eta$, if $\Q(d, \T_\U)$ doesn't exist then the $d$ realizes back into $\W$. We now argue that $\Q(d, \T_U)$ exists.\\

\textit{Claim 2.} $\Q(d, \T_U)$ exists.
\begin{proof}  Towards a contradiction assume $\Q(d, \T_U)$ doesn't exist. Thus, $d\cap D^{\T_U}=\emptyset$ and $\pi^{\T_U}(\d^{\R_U})=\d(\T_U)$. Set $\N=\M^{\T_U}_d$ and $j=\pi^{\T_U}_d$\footnote{Here, $\T_U$ is a tree on $\W_U$ but it is based on $\R_U$.}. Because $\Q(d, \T_U)$ doesn't exist, we have that $\N\models ``\d(\T_U)$ is a Woodin cardinal".

 We have that $j(\Q(b_U, \T_U))\in \N$ and is authenticated in $\N$. Let $\gg=j(\b_U)$. Then $\N|\gg$ has Woodin cardinals that are bigger than $\d(j(\T_U))$. Let $\d$ be the least one that is bigger than $Ord\cap j(\Q(b_U, \T_U))$. We can now iterate $\N$ below $\d$ but above $Ord\cap j(\Q(b_U, \T_U))$ to make $\Q(b_U, \T_U)$ generic for the extender algebra at the image of $\d$. This iteration produces $i:\N\rightarrow \N_1$ such that $\cp(i)>\d(j(\T_U))$. Letting $h\subseteq Coll(\omega, i(\d))$ be $\N_1$-generic such that $\Q(b_U, \T_U)\in \N_1[h]$, we can find
 \begin{center}
  $l:\Q(b_U, \T_U)\rightarrow j(\Q(b_U, \T_U))$\footnote{As $\cp(i)>\d(j(\T_U))$, $i(j(\Q(b_U, \T_U)))=j(\Q(b_U, \T_U))$.}
  \end{center} such that 
  \begin{itemize}
  \item $l\in \N_1[h]$,
  \item $l\rest (\m^+(\T_U))^b=\pi^{j(\T_U), b}$. 
  \end{itemize}
  As $\N_1[h]\models ``j(\Q(b_U, \T_U))$ is authenticated and has an authenticated strategy", $\N_1[h]\models ``\Q(b_U, \T_U)$ has an  authenticated iteration strategy", and hence $\Q(b_U, \T_U)$ is definable in  $\N_1[h]$ from objects in $\N_1$. It follows that $\Q(b_U, \T_U)\in \N_1$, implying that $\N_1\models ``\d(\T_U)$ is not a Woodin cardinal". Hence, $\N\models ``\d(\T_U)$ is not a Woodin cardinal". Therefore, $\Q(d, \T_U)$ exists.
\end{proof}

\end{proof}

\subsection{A conclusion}

\begin{proposition}\label{getting suitable} Suppose $\vec{\V}$ is a small array with the $Z$-realizability property. Then either 
\begin{enumerate}
\item $\V_\eta$ has a $Z$-validated iteration strategy

or
\item there is a $Z$-validated $\sf{nuvs}$ iteration $p$ of $\V_\eta$ such that $\m^+(p)$ is $Z$-suitable\footnote{See \rdef{z-full}.}. 
\end{enumerate}
\end{proposition}
\begin{proof} The proof has already been given in the previous subsections. Suppose that $\V_\eta$ does not have a $Z$-validated iteration strategy. The proof of \rprop{lem:break2a} shows that if $p$ is a $Z$-validated $\sf{uvs}$ iteration of $\V_\eta$ of limit length then there is a unique branch $b$ of $p$ such that $p^\frown\{b\}$ is $Z$-validated. Therefore, since picking $Z$-validated branches is not defining an iteration strategy for $\V_\eta$, we must have a $\sf{nuvs}$ $Z$-validated iteration $p$ of $\V_\eta$ which does not have a $Z$-validated branch.\footnote{Note that there may not be any $\Q$-structure for $p$.} 

We now claim that $\m^+(p)$ is a $Z$-suitable hod premouse. Indeed, suppose there is some $Z$-validated sts premouse $\Q$ extending $\R=_{def}\m^+(p)$ such that $\Q$ is a $\Q$-structure for $p$. Let then $U$ be a good hull such that$\{\vec{\V}, p, \Q\}\in U$. Appealing to \rprop{bottom part realizability}, we now have $\b\leq lh(\vec{\V})$, a branch $b$ of $p_U$ such that $\Q(b, p_U)$ exists and a weak $l$-embedding $k:\M_b^{p_U}\rightarrow \mathcal{C}_l(\V_\b)$ for an appropriate $l$. It follows that $\Q(b, p_U)$ is $Z$-approved and hence, $\Q(b, p_U)=\Q_U$. Because $\Q_U\in M_U$, we have that $b\in M_U$. Then $c=_{def}\pi_U(b)$ is a (cofinal) branch of $p$ such that $p^\frown\{c\}$ is $Z$-validated. 
\end{proof}

\section{Hybrid fully backgrounded constructions}\label{sec:Kc}

The goal of this section is to adopt the $K^c$-construction used in \cite{hod_mice_LSA} to our current situation. As we have the large cardinals in $V$, it is easier to perform fully backgrounded constructions than using partial background certificates. For instance, proofs of iterability will be easier. 

The construction that we intend to perform will produce an almost excellent (see \rdef{dfn:hod_pm}) hod premouse $\P$ extending $\mH$. The construction will be done in $V$. 

The fully backgrounded construction that we have in mind has two different backgrounding conditions for extenders. The extenders with critical point $>\Theta =_{def} \delta^\mathcal{H}$ will have total extenders as their background certificates. The extenders with critical point $\Theta$ will be authenticated  by good hulls. We call this construction the \textit{the hybrid fully backgrounded construction} over $\mH$, and denote it by $\sf{HFBC(\mu)}$. 
 
We fix a condensing set $Z\in Cnd(\mH)\cap V$. While $Z$ will appear in our authentication definitions, it can be shown that $\sf{HFBC(\mu)}$ does not depend on $Z$. $\sf{HFBC(\mu)}$ proceeds more or less according to the usual procedure for building hod pairs until we reach a weakly $Z$-suitable stage $\R$. At this stage, we must continue with a fully backgrounded $Z$-validated sts construction over $\R$. If this construction produces a $\Q$-structure for $\R$ then we attempt to construct a $Z$-validated strategy for it. Failing to do so will produce our honest $Z$-suitable $\R$ as in \rprop{getting suitable}.


$\sf{HFBC}$  can fail in the usual ways, by producing a level whose countable substructures are not iterable\footnote{The real reason for a failure of such constructions is failure of universality or solidity both of which are consequences of iterability.}. However, our constructions are aimed at producing models with strictly weaker large cardinal structure than those for which we know how to prove iterability. In particular, the main theorem of \cite{Neeman} implies that $\sf{HFBC}$ does not fail because of issues having to do with iterability.

We should say that the construction that follows is an adaptation of a similar construction introduced in \cite[Section 10.2.9 and 12.2]{hod_mice_LSA}. Because of this we will not dwell too much on how extenders with critical point $\d^\mH$ are chosen. The reader may consult \cite[Lemma 12.3.15]{hod_mice_LSA}. The first of these constructions used fully backgrounded certificates like we will do in the next subsection. It  was used to prove the Mouse Set Conjecture in the minimal model of $\sf{LSA}$. The second was used to construct a model of $\sf{LSA}$ from $\sf{PFA}$.

\subsection{The levels of $\sf{HFBC}(\mu)$}\label{sec:nonsuitable}

We assume that $T$ holds and let $(S, S_0, \nu_0,\vec{Y}, \vec{A})$ witness it. Let $\mu\in S_0$. When discussing the CMI objects at $\mu$ we will omit $\mu$ from our notation. 

Below we will define the sequence $(\M_\xi,\Sigma_\xi: \xi \leq \Omega')$, where $\Omega' \leq Ord$, of levels of the hybrid fully backgrounded construction at $\mu$. Here, we develop the terminology that we will use to describe the passage from $\M_\xi$ to $\M_{\xi+1}$.

$\sf{HFBC}$ resembles the $K^c$-construction of \cite[Section 12.2]{hod_mice_LSA} except that we require that the extenders with critical point $>\delta^\mH$ used in the construction have total certificates in the sense of \cite[Chapter 12]{FSIT}.  Because our construction does not reach a Woodin cardinal that is a limit of Woodin cardinals, the results of \cite{Neeman} apply. For instance, \cite[Theorem 1.1]{Neeman} will be used to conclude that the countable submodels of each $\M_\xi$ are $\omega_1+1$-iterable. Other theorems that we will use from \cite{Neeman} are \cite[Theorem 2.1, 2.10, 3.11 and Corollary 3.14]{Neeman}. 

Say that $\M$ \textit{nicely extends} $\mH$ if there is a $(\mu, \mu, Z)$-good hull $U$ such that $\M_U$ nicely extends $\Q^Z_{U\cap \mH}$.  Suppose now that $\M$ is a $Z$-validated hod premouse nicely extending $\mH$. Set $mo(\M)=o^\M(\d^{\mH})$\footnote{``mo" stands for the ``Mitchell Order".}. 
\begin{definition}\label{appropriate}
Given a $Z$-validated hod premouse nicely extending $\mH$ we say $\M$ is \textbf{appropriate} if $Ord\cap \M=mo(\M)$ and $\M\models ``$there are no Woodin cardinals in the interval $[\d^{\mH}, mo(\M))"$.
\end{definition}
Given an appropriate $\M$, we would like to describe the next appropriate $Z$-validated hod premouse. We do this by preparing $\M$, which involves building over $\M$ some mild structures in order to reach the next stage that is either a stage where we can add an extender or is a weakly $Z$-suitable stage. In the latter case we will put $\sf{HFBC}(\mu)$ on hold and continue with the $Z$-validated sts construction. The preparation of $\M$ has two stages. We first add a sharp to $\M$ and then close the resulting hod premouse under its strategy. Each of these constructions can change $\M$ as they can reach levels that project across $\M$.
The functions that we referred to above are $next_{\#}$, $next_s$, $next_{bex}$ and $next_{\Theta-ex}$. \\\\
$next_{bex}(\M)$

This function simply adds a backgrounded extender to $\M$. Suppose that $\M$ is appropriate. We say $next_{bex}(\M)$ is \textit{almost successful} if there is a triple $(\k, \l, F)$ such that 
\begin{enumerate}
\item $\k<\l$ are inaccessible cardinals $>\d^\mH$,
\item $F$ is a $(\k, \l)$-extender such that $V_\l\subseteq Ult(V, F)$,
\item letting $G=F\cap \M$, $(\M, G)$ is hod premouse.
\end{enumerate}
We say $next_{bex}(\M)$ is \textit{successful} if it is almost successful and there is a unique triple $(\k, \l, F)$ as above such that if $G=F\cap \M$, $(\M, G)$ is a solid and universal $Z$-validated hod premouse with a $Z$-validated strategy and such that $\rho(\M, G)>\d^{\mH}$. 

Suppose now that $\M$ is appropriate. We say $\M$ has \textit{badness  type 0} if $next_{bex}(\M)$ is almost successful but it is not successful. We write $bad(\M)=0$. If $next_{bex}(\M)$ is successful or not almost successful then
\begin{enumerate}
\item if $next_{bex}(\M)$ is successful then letting $(\k, \l, F)$ be the unique triple witnessing the success of $next_{bex}(\M)$, we let $next_{bex}(\M)=(\M, G)$ where  $G=F\cap \M$.
\item If $next_{bex}(\M)$ is not almost successful then we let $next_{bex}(\M)=\M$.
\end{enumerate}
$next_\#(\M)$

 Suppose $\M$ is appropriate and $bad(\M)\not=0$. We let $next_\#(\M)$ be build as follows: Let $(\M_i: i\leq k)$ be a sequence of $Z$-validated hod premice  defined as follows:
\begin{enumerate}
\item $\M_0=next_{bex}(\M)$.
\item If $i+1\leq k$ then there is $\M^*$ that is an initial segment of $J[\M_i]$ such that $\rho(\M^*)<mo(\M_i)$ and letting $\M^*$ be the least such initial segment of $J[\M_i]$, $\M^*$ is solid and universal, $\rho(\M^*)>\d^{\mH}$ and  $\M_{i+1}=\mathcal{C}(\M^*)$.
\item $k$ is least such that either (i) no level of $J[\M_k]$ projects across $mo(\M_k)$ or (ii) some level of $J[\M_k]$ projects to or below $\d^{\mH}$.
\end{enumerate}  
We say $next_\#(\M)$ is \textit{successful} if 
\begin{enumerate}[(a)]
\item clause 3(ii) doesn't happen, 
\item $\M_k^{\#}$ is solid and universal,
\item $\rho(\M_k^{\#})>\d^{\mH}$.
\end{enumerate}
 If $next_\#(\M)$ is successful then let $next_{\#}(\M)=\mathcal{C}(\M_k^\#)$. 

Suppose now that $\M$ is appropriate. We say $\M$ has \textit{badness  type is 1} if $next_{\#}(\M)$ is not successful. We write $bad(\M)=1$. 

Suppose now that $\M$ is appropriate and $bad(\M)\not =0, 1$. We say $\M$ has \textit{badness type 2} if $next_{\#}(\M)$ is not weakly $Z$-suitable and  does not have a $Z$-validated  strategy. We write $bad(\M)=2$. \\\\
$next_s(\M)$

Suppose now that $\M$ is appropriate, $bad(\M)\not=0,1, 2$ and $next_{\#}(\M)$ is not weakly $Z$-suitable. Let $\Sigma$ be the unique $Z$-validated strategy of $next_\#(\M)$. We now define $next_s(\M)$ which, in a sense, adds $Lp^{\Sigma}(next_s(\M))$ to $\M$.  

We let $next_s(\M)$ be build as follows: Let $(\M_i, \Sigma_i: i\leq k)$ be a sequence of $Z$-validated hod premice along with their $Z$-validated strategies defined as follows:
\begin{enumerate}
\item $\M_0=\N$ and $\Sigma_0=\Sigma$.
\item If $i+1\leq k$ then there is $\M^*$ that is an initial segment of $J[\vec{E}, \Sigma_i](\M_i)$\footnote{Here and below by $J[\vec{E}, \Sigma_i]$ we mean the fully backgrounded construction relative to $\Sigma_i$. $J[\vec{E}, \Sigma_i](A)$ is the aforementioned construction done over $A$.}  such that $\rho(\M^*)<mo(\M_i)$ and letting $\M^*$ be the least such initial segment of $J[\vec{E}, \Sigma_i](\M_i)$, $\M^*$ is solid and universal, $\rho(\M^*)>\d^\mH$, $\M_{i+1}=\mathcal{C}(\M^*)$ and $\Sigma_{i+1}$ is the unique $Z$-validated strategy of $\M_{i+1}$. 
\item $k$ is least such that either (i) no level of $J[\vec{E}, \Sigma_k](\M_k)$ projects across $mo(\M_k)$ or (ii) $\M_k$ does not have a $Z$-validated strategy, or (iii) some level of $J[\vec{E}, \Sigma_k](\M_k)$ projects to or below $\d^{\mH}$.
\end{enumerate}
We say $next_s(\M)$ is \textit{successful} if clause 3(ii)-(iii) don't happen. If $next_s(\M)$ is successful then let $next_{s}(\M)=J[\vec{E}, \Sigma_k](\M_k)|\a$ where \begin{center} $\a=(mo(\M_k)^+)^{J[\vec{E}, \Sigma_k](\M_k)}$.\end{center}
Suppose now that $\M$ is appropriate. We say $\M$ has \textit{badness  type  3} if $bad(\M)\not=0,1, 2$ and $next_{s}(\M)$ is not successful. We write $bad(\M)=3$. We say $\M$ has  \textit{badness  type is 4} if $bad(\M)\not=0,1,2, 3$ and $next_s(\M)$ doesn't have a $Z$-validated strategy.\\\\
$next_{\Theta-ex}$

Suppose now that $\M$ is appropriate and $bad(\M)\not \in 5$. Let $\Sigma$ be the unique $Z$-validated strategy of $next_s(\M)$. 
We say that $next_{\Theta-ex}(\M)$ is successful if there is a unique $\M$-extender $F$ such that
\begin{enumerate}
\item $\cp(F)=\d^{\mH}$,
\item $(\M, F)$ is a $Z$-validated hod mouse,
\item $\rho((\M, F))>\d^{\mH}$.
\end{enumerate}
If $next_{\Theta-ex}(\M)$ is successful then we let $next_{\Theta-ex}(\M)=\mathcal{C}((\M, F))$ where $F$ is as above. 
 We say $\M$ has \textit{badness type 5} if $next_{\Theta-ex}(\M)$ is not successful. In this case, we write $bad(\M)=5$.
 
 \begin{definition}\label{bad m} Suppose $\M$ is appropriate. We say $\M$ is \textbf{bad} if $bad(\M)$ is defined.
 \end{definition}
 \begin{definition}\label{next m} Suppose $\M$ is appropriate. 
 If $\M$ is not bad then we let $next(\M)=next_{\Theta-ex}(\M)$. 
 \end{definition}
 
 \begin{remark} In order for $\M\in dom(next)$ it is necessary that $next_{\#}(\M)$ is not a weakly $Z$-suitable level. 
 \end{remark}
 
The $next$ function defined above gives us the next model in $\sf{HFBC}(\mu)$, but it doesn't tell us how to start the construction. We will start $\sf{HFBC}(\mu)$ with $\mH$, which is an appropriate hod premouse. However, if we encounter $\M$ such that $next_{\#}(\M)$ is weakly $Z$-suitable then we have to continue with the $Z$-validated sts construction. We get back to $\sf{HFBC}(\mu)$ once we produce the canonical witness to non-Woodiness of $\d^{next_{\#}(\M)}$. What we do next is we define the $start$ function whose domain will consist of objects that the $Z$-validated sts construction produces on top of $next_{\#}(\M)$.\\\\
$start(\R)$

Suppose $\R$ is a $Z$-validated hod mouse such that 
\begin{enumerate}
\item $mo(\R)$ is a Woodin cardinal of $\R$,
\item $(\R|mo(\R))^\#\insegeq \R$,
\item $\R$ is sound,
\item $\R$ is an sts premouse over $(\R|mo(\R))^\#$ such that $rud(\R)\models ``mo(\R)$ is not a Woodin cardinal".
\end{enumerate}
If $\R$ is as above then we write $\R\in dom(stop)$. Let $\Sigma$ be the $Z$-validated strategy of $\R$ if it exists; in the case it does not exist, we declare $start_0(\R)$ is unsuccessful and letting $p$ on $\R$ as in Proposition \ref{getting suitable}, we then switch to the f.b. $(Z, \lambda)$-validated sts construction over $\m^+(p)$. In the case $\Sigma$ exists, we define $start_0(\R)$ just like we defined $next_{s}(\M)$ above. If $start_0(\R)$ is successful then it will output a $Z$-validated hod mouse $\W$ that nicely extends $\mH$ and has a $Z$-validated strategy $\Lambda$. Moreover, for some $\d\in \W$, 
\begin{enumerate}
\item $(\W|\d)^\#\insegeq \W$ is of lsa type,
\item $\W\models ``\d$ is not a Woodin cardinal",
\item letting $\W^*\inseg\W$ be largest such that $\W^*\models ``\d$ is a Woodin cardinal",  $\W^*$ is a $\Lambda^{stc}_{(\W|\d)^\#}$-sts mouse over  $(\W|\d)^\#$ and $\W=J[\vec{E}, \Lambda_{\W^*}]|(\d^+)^{J[\vec{E}, \Lambda_{\W^*}]}$. 
\end{enumerate}
Next let $start_1(\R)$ be defined just like $next_{\Theta-ex}(\M)$ starting with $start_0(\R)$. For $\R\in dom(start)$ we say $start(\R)$ is successful if both $start_0(\R)$ and $start_1(\R)$ are successful, and we let $start(\R)$ be the model that $start_1(\R)$ outputs. 

\begin{definition}\label{ready for hfbc} Suppose $\R\in dom(start)$. We say $\R$ is ready for $\sf{HFBC}(\mu)$ if $start(\R)$ is successful. 
\end{definition}
Notice that if $start(\R)$ is successful then $mo(start(\R))=Ord\cap start(\R)$. We end this subsection with the definition of $\sf{HFBC}(\mu)$. \\\\
\textbf{Levels of $\sf{HFBC}$}

Suppose $\M=\mH$ or $\M=start(\R)$ for some $\R\in dom(start)$ such that $start(\R)$ is successful. Let $\Sigma$ be the unique $Z$-validated strategy of $\M$. 
\begin{definition}\label{hfbc} We say $(\M_\xi,\Sigma_\xi, : \xi < \Omega')$, where $\Omega' \leq Ord$, are the levels of the hybrid fully backgrounded construction at $\mu$ ($\sf{HFBC}(\mu)$) done with respect to $(\M, \Sigma)$ if the following conditions hold.
\begin{enumerate}
\item $\M_0=\M$ and $\Sigma_0=\Sigma$.
\item For each $\xi<\Omega'$, $\M_\xi$ is appropriate and $\Sigma_\xi$ is the unique $Z$-validated strategy of $\M_\xi$.
\item For all $\xi<\Omega'$ if $\xi+1<\Omega'$ then $\M_\xi$ is not bad and $\M_{\xi+1}=next(\M_\xi)$.
\item For all $\xi< \Omega'$, if $\xi$ is a limit ordinal then letting $\M^*_\xi=_{def}liminf_{\a\rightarrow \xi}\M_\a$, $\M^*_\xi$ is appropriate and not bad and $\M_\xi=next(\M_\xi^*)$.
\item $\Omega'$ is the least ordinal $\a$ such that one of the following conditions hold:
\begin{enumerate}
\item $\a$ is a limit ordinal and $\M^*_\a$ is bad.
\item $\a$ is a limit ordinal and $next_{\#}(\M^*_\a)$ is weakly $Z$-suitable.
\item $\a=\b+1$ and $\M_\b$ is bad.
\end{enumerate}
\end{enumerate}
\end{definition}

We say $\sf{HFBC}(\mu)$ converges if $\Omega'$ is as in clause 5(b), i.e., $next_{\#}(\M^*_{\Omega'})$ is weakly $Z$-suitable. The following proposition is essentially \cite[Theorem 11.3]{FSIT}.
\begin{proposition}\label{reaching suitable levels} Suppose $\d>\mu$ is a Woodin cardinal. Then if for all $\xi<\d$, $\M_\xi$ is defined then $\Omega'=\d$ and $\sf{HFBC}$ converges. Moreover, letting $\P^-=liminf_{\xi\rightarrow \d}\M_\xi$ and $\P=(\P^-)^\#$ then $\P$ is weakly $Z$-suitable.
\end{proposition}

Letting $\d>\mu$ be the least Woodin cardinal $>\mu$, we need to show that $\sf{HFBC}(\mu)$ either lasts $\d$ steps or encounters a weakly $Z$-suitable stage. Recall that we defined $\sf{HFBC}(\mu)$ over some $(\M, \Sigma)$.

\section{Putting it all together}\label{sec:all_together}

We are assuming theory $T$ and let $(S, S_0, \nu_0, \vec{Y}, \vec{A})$ witness it; let $\mu\in S_0$. Combining $\sf{HFBC}(\mu)$ with the fully backgrounded $Z$-validated sts construction, as shown by \rprop{getting suitable}, we see that we reach an honest $Z$-suitable $\R$. In this section, we would like to continue the $Z$-validated sts construction over $\R$ and show that it must reach an excellent $\P$. To do this, we will stack fully backgrounded $Z$-validated sts constructions one on the top of another to reach an almost excellent hybrid premouse which we will show has external iterability. We will then need some arguments that translate iterable almost excellent hybrid premice into an excellent ones. 

This stacking idea might be a little bit unnatural but it seems the most straightforward way of dealing with the two main issues at hand. What we would really like to do is to perform the fully backgrounded $Z$-validated sts construction over $\R$ and hope that it will reach an excellent hybrid premouse. There are two key issues that arise. The final model of our construction has to inherit a stationary class of measurable cardinals. Perhaps the most straightforward way of dealing with this issue is to attempt to show that every measurable cardinal $\k$ such that no cardinal is $\k$-strong remains measurable in the output of the backgrounded construction. We do not know how to show this without working with more complex forms of backgrounded constructions. Our solution involves just adding the measure by ``brute force". Once the construction reaches one such $\k$, we will continue by adding the measure coarsely, much like one does in the construction of $L[\mu]$.

The next issue is to guarantee window based iterability. The most natural way of accomplishing this is by showing that the models of our backgrounded construction are iterable. However, this may not work and if it fails, it fails as follows. Suppose $\N$ is a model appearing in the fully backgrounded $Z$-validated sts construction over $\R$ and $\k$ is a cutpoint cardinal of $\N$. Suppose $\N$ has no Woodin cardinals above $\k$.  We now seek a $Z$-validated strategy for $\N$ that acts on iterations above $\k$. If such a strategy doesn't exist then we must have a tree $\T$ on $\N$ above $\kappa$ which does not have a $Z$-validated $\Q$-structure. Let then $\N_1^-=\m(\T)$. It follows that if we perform a fully backgrounded $Z$-validated sts construction over $\N_1^-$ we will not reach a $\Q$-structure for $\T$. Let then $\N_1$ be the one cardinal extension of $\N_1^-$ built by the fully backgrounded $Z$-validated sts construction over $\N_1^-$. We thus have that $\N_1\models ``\d(\T)$ is a Woodin cardinal".

We now want to see that $\N_1$ has a window based iterability. Let then $\eta_1\in (\k, \d(\T))$ be a regular cardinal of $\N_1$, and we want to argue that $\N_1|\eta_1$ is iterable. The strategy we seek is again a $Z$-validated strategy. If it doesn't exist then we get a tree $\T_1$ on $\N_1|\eta_1$ such that $\T_1$ does not have a $Z$-validated $\Q$-structure. The construction above produced $\N_2$ extending $\m(\T_1)$.  The goal now is to show that $\N_2$ has window based iterability. Failure of such a strategy produced $\eta_2\in (\k, \d(\T_1))$, $\T_2$ based on $\N_2$ that is above $\k$ and a model $\N_3$ extending $\m(\T_2)$. The process outlined above cannot last $\omega$ many steps, for if it did we will have a sequence $(\N_i, \T_i: i<\omega)$ and a reflected version of this sequence cannot have a well-founded direct limit along the realizable branches. 

There is yet another issue that we need to deal with which is not connected with the stacking construction, but has to do with other aspects of the construction. We will need arguments that will show window based iterability in $V$ can somehow be reflected inside the sts premice alluded above. To show this, we will need to break into cases and examine exactly how we ended up with the model we seek. For this reason, we isolate the following hypothesis. \\\\
$\sf{Hypo}:$ For some $X$ containing $\R$ there is a sound $Z$-validated almost excellent mouse $\M$ over $X$ that is based on $\R$.\\\\
The following essentially follows from the main results of \cite{Neeman}.
\begin{proposition}\label{hypo implies hypo1}
Assume $\neg \sf{Hypo}$. Then for any $X$ containing $\R$, letting $\d$ be the least Woodin cardinal such that $X\in V_\d$, no model of the $Z$-validated sts construction of $V_\d$ that is based on $\R$ and is done over $X$, reaches an almost excellent hybrid premouse.
\end{proposition}

\subsection{The prototypical branch existence argument}\label{prototypical argument}

Here we present an argument due to John Steel that we will use over and over again. The argument is general and can be used in many settings. We will refer to this argument as the \textit{prototypical branch existence argument}. In the sequel, when we need to prove something via the same argument we will just say that ``the prototypical branch existence argument shows...". \\

\textbf{The prototypical branch existence argument}\\

Suppose $\d$ is a Woodin cardinal, $X\in V_\d$ is a set such that $\R\in X$ and $(\M_\a, \N_\a: \a\leq \d)$ are the models of the fully backgrounded $Z$-validated sts construction done over $X$. Fix $\a\leq \d$ and suppose that $\N_\a$ has no Woodin cardinals (as an sts premouse over $X$). Let $\T$ be a normal $Z$-validated iteration of $\N_\a$ such that for every limit $\b<lh(\T)$ if $c_\b=[0, \b]_T$ then $\Q(c_\b, \T\rest \b)$ exists and is $Z$-validated. Suppose that $\T$ has limit length and there is a $Z$-validated sts mouse $\Q$ such that $\m(\T)\insegeq \Q$ and $rud(\Q)\models ``\d(\T)$ is not a Woodin cardinal". Then there is a branch $b$ of $\T$ such that $\Q(b, \T)$ exists and is equal to $\Q$. 

The argument proceeds as follows. Fix some $\lambda\in S_0-\d$ and let $\pi_U: M_U\rightarrow H_\zeta$ be a $(\mu, \l, Z)$-good hull such that $\T, \Q\in rng(\pi)$. Let $\nu=\card{M_U}$ and let $g\subseteq Coll(\omega, \nu)$ be generic. Then there is a maximal branch $c$ of $\T_U$, $\b\leq \a$ and a (weak) embedding $\sigma:\M^{\T_U}_c\rightarrow \N_\b$ such that if $c$ is non-dropping then $\b=\a$ and $\pi_U=\sigma\circ \pi^{\T_U}$. Arguing as in \rprop{z-validated are embeddable}, we get that $c$ must be a cofinal branch and that $\Q(c, \T_U)$ must exist and be equal to $\Q_U$. It follows then that $c\in M_U$. Hence, $\pi_U(c)$ is as desired. 

\begin{remark} It is important to keep in mind that the argument doesn't work when $\N_\a$ has Woodin cardinals as then $\Q(c, \T_U)$ may not exist. Thus, this argument cannot in general be used to show that levels of $K^c$ are short tree iterable.
\end{remark}

\subsection{One step construction}

Suppose $X$ is a set such that $\R\in X$. The main goal of this section is to produce a \textit{short-tree-iterable} $Z$-suitable sts hod premouse over $X$. Here short tree iterability is in the sense of the $\H$ analysis (cf. \cite{hod_as_core_model}).
\begin{definition}
Suppose $\P$ is a $Z$-validated sts premouse over $X$ based on $\R$. We say $\P$ is  \textbf{almost $Z$-good} if $\P$, as an sts premouse over $X$, has a unique Woodin cardinal $\d^\P$ such that 
\begin{enumerate}
\item $\P=(\P|\d^\P)^\#$,
\item if $\M$ is a sound $Z$-validated sts mouse over $\P$  then $\M\models ``\d^\P$ is a Woodin cardinal".
\end{enumerate}
We say $\P$ is \textbf{$Z$-good} if $\P$ has a unique Woodin cardinal $\d^\P$ such that 
\begin{enumerate}
\item $(\P|\d^\P)^\#$ is almost $Z$-good,
\item $\P=Lp^{Zv, sts}(\P|\d^\P)$,
\item for every regular cardinal $\eta<\d^\P$, $\P|\eta$ has a $Z$-validated strategy.
\end{enumerate}
\end{definition}

We say that the $Z$-good $\P$ is fully backgrounded if for some maximal window $w$ and for some $\xi\in w$, $\P|\d^\P$ is a model appearing in the fully backgrounded $Z$-validated sts construction of $V_\eta$ which uses extenders with critical point $>\xi$.  

\begin{proposition}\label{one step construction} Assume $\neg \sf{Hypo}$. There is a $Z$-good fully backgrounded sts premouse  over  $X$ based on $\R$. 
\end{proposition}

We spend this entire subsection proving \rprop{one step construction}. We will do it in two steps. In the first step we will produce a fully backgrounded almost $Z$-good $\N$. Then we will obtain a fully backgrounded $Z$-good $\P$. We start with the first step.

\begin{lemma} Assume $\neg \sf{Hypo}$. There is an almost  $Z$-good fully backgrounded sts premouse over  $X$ based on $\R$. 
\end{lemma}
\begin{proof} Let $\d$ be the least Woodin cardinal of $V$ such that $X\in V_\d$. Let $(\M_\xi, \N_\xi: \xi\leq \Omega^*)$ be the models of the fully backgrounded  $Z$-validated sts construction of $V_\d$ done over $X$ (based on $\R$). Because we are assuming $\neg \sf{Hypo}$, $\Omega^*=\d$. We claim that\\

\textit{Claim.} there is $\xi\leq \d$ such that $\N_\xi$ is an almost $Z$-good sts premouse.
\begin{proof} Suppose for every $\xi<\d$, $\N_\xi$ is not almost $Z$-good. We show that $\N=(\N_\d)^\#$ is an almost $Z$-good sts premouse. A standard reflection argument shows that $\rho_\omega(\N)<\d$. Suppose then $\N$ is not almost $Z$-good and fix $\M$ such that
\begin{enumerate}
\item $\N\insegeq \M$,
\item $\M$ is sound above $\d$, 
\item $\rho_\omega(\d)<\d$,
\item $\M$ is $Z$-validated and has a $Z$-validated $Ord$-strategy. 
\end{enumerate}
 As we are are assuming $\neg \sf{Hypo}$, $\d$ is not a limit of Woodin cardinals in $\N$. Let then $\pi: \M^*\rightarrow \M$ be such that letting $\cp(\pi)=\nu$, $\pi(\nu)=\d$ and $\N$ has no Woodin cardinals in the interval $[\nu, \d)$. 

Working inside $\N$, let $\N'$ be the output of the $(\R, \R^b, S^\N)$-authenticated construction done over $\N|\nu+1$ using extenders with critical points $>\nu$. Let $\M'$ be the result of translating $\M$ over to $\N'$ via the $S$-construction (see \cite[Chapter 6.4]{hod_mice_LSA}). Similarly, for each $\N$-cardinal $\xi>\nu$ such that $(\N|\xi)^\#\models ``\xi$ is a Woodin cardinal" let $\M_\xi$ witness that our proposition fails for $(\N|\xi)^\#$. For each such $\xi$ let $\M'_\xi$ be the result of translating $\M_\xi$ over to $\N'|\xi$.  

We now compare $\M^*$ with the construction producing $\N'$. In this comparison, only $\M^*$ is moving. We claim that this comparison lasts $\d+1$-steps producing a tree $\T$ on $\M^*$ with last model $\M'$. Indeed, given $\T\rest \a$ where $\a\leq \d$ is a limit ordinal, if $\m^+(\T\rest \a)\models ``\d(\T\rest \a)$ is a Woodin cardinal" then $\M_\a'$ is defined. Because $\M$ is a $Z$-validated sts mouse, we must have a unique cofinal well-founded branch $b_\a$ of $\T\rest \a$ such that $\Q(b, \T\rest \a)$ is defined and is equal to $\M_\a'$. We then pick this branch $b_\a$ at stage $\a$. 

It must now be clear that the existence of $\T$ violates universality; this implies by standard results that there must be a superstrong cardinal in $\N$. The reader may consult \cite[Lemma 5.4]{hod_mice}.
\end{proof}
The claim finishes the proof.
\end{proof}

We now start the proof of \rprop{one step construction}. Towards a contradiction, we assume that there is no fully backgrounded $Z$-good sts premouse over $X$ based on $\R$. Let $\N_0$ be a fully backgrounded almost $Z$-good sts premouse over $X$ based on $\R$. 

Below given $\S$, $\d^\S$ will always denote the largest Woodin of $\S$ and $w^\S$ will denote the maximal window $w$ of $\S$ such that $\d^w=\d^\S$.

We now by induction produce an infinite sequence $(\N_i, \nu_i, \T_i: i<\omega)$ such that 
\begin{enumerate}
\item for every $i<\omega$, $\N_i$ is a fully backgrounded almost $Z$-good sts premouse over $X$ based on $\R$,
\item for every $i<\omega$, $\nu_i$ is a successor cardinal of $\N_i$,
\item for every $i<\omega$, $\T_i$ is a normal $Z$-validated iteration of $\N_i|\nu_i$ such that $\T_i$ has no cofinal well-founded branch $b$ such that $\Q(b, \T_i)$ exists and is a $Z$-validated mouse,
\item for every $i<\omega$, $\N_{i+1}=\m^+(\T_i)$.
\end{enumerate}
A simple reflection argument shows that such a sequence cannot exist. The fact that for $i>0$, $\N_i$ is fully backgrounded is irrelevant for the reflected argument alluded in the previous sentence. It is enough that $\N_0$ is fully backgrounded. 

Assume then we have built $\N_i$ and we now describe the procedure for  getting $(\nu_i, \T_i, \N_{i+1})$. Because $\N_i$ is not $Z$-good, there is a $\nu_i\in w^{\N_i}$ which is  a successor cardinal of $\N_i$ and $\N_i|\nu_i$ does not have a $Z$-validated strategy.

 Let $\eta$ be some Woodin cardinal such that $\N_i\in V_\eta$ and let $\xi_\eta<\eta$ be such that $\N_i\in H_{\xi_\eta}$ and there are no Woodin cardinals in the interval $(\xi_\eta, \eta)$. Let $(\M^\eta_\a, \S^\eta_\a: \a\leq \eta)$ be the models of the fully backgrounded $Z$-validated sts construction of $V_\eta$ done over $X$ using extenders with critical points $>\xi_\eta$.

As $\N_i\in H_\xi$, in the comparison of $\W=_{def}\N_i|\nu_i$ with the construction $(\M_\a^\eta, \S^\eta_\a: \a\leq \eta)$ only $\W$ moves. We now analyze the tree on $\W$. Suppose $\T^\eta$ is the tree on $\W$ built via the above comparison and suppose $\T^\eta$ has a limit length. We now have two cases. \\\\
\textbf{Case1.}  Suppose there is $\a$ such that $\S^\eta_\a\models ``\d(\T)$ is not a Woodin cardinal". It follows from the prototypical argument that there must be a cofinal branch $b$ of $\T$ such that $\Q(b, \T)$ exists and is a $Z$-validated mouse. We then extend $\T^\eta$ by adding $b$. \\
\textbf{Case2.} Suppose there is no $\a\leq \eta$ such that $\S^\eta_\a\models ``\d(\T)$ is not a Woodin cardinal". In this case, we stop the construction and set $\N_{i+1}=\m^+(\T)$ and $\T_i=\T$. \\\\
We stop the construction if either Case2 holds or for some $\a$, $\M^\eta_\a$ is the last model of $\T$ and $\pi^\T$ exists. 

We now claim that for some $\eta$, the construction of $\T^\eta$ stops because of Case2. Assume otherwise. Then for each Woodin cardinal $\eta$ we have $\a_\eta$ and an embedding $\pi:\W\rightarrow \M^\eta_{\a_\eta}$. As $\W$ has no Woodin cardinals, it follows that for every $\eta$, $\W$ is $\xi_\eta$-iterable via a $Z$-validated strategy. As $Ord=\bigcup_\eta \xi_\eta$, we have that $\W$ is $Ord$-iterable via a $Z$-validated strategy. Hence, for some $\eta$, Case2 must be the cause for stopping the construction of $\T^\eta$. Below we drop $\eta$ from subscripts. 

To finish the proof of \rprop{one step construction} we need to show that $\N_{i+1}$ is almost $Z$-good. This easily follows from universality. Because $\d(\T )$ is Woodin in $\S_\eta$, we must have that $\N_{i+1}\insegeq \S_\eta$. If now $\M$ is a $\d(\T)$-sound $Z$-validated sts mouse then because $\T$ has no $\Q$-structure, we must have that $\rho(\M)\geq \d(\T)$ and $\M\models ``\d(\T)$ is a Woodin cardinal" (as otherwise the prototypical argument would yield a branch of $\T$).

\subsection{Stacking suitable sts mice}\label{stacking suitable sts mice}

In this section, assuming $\neg\sf{Hypo}$ we build an almost excellent hybrid premouse. We achieve this by stacking fully backgrounded $Z$-good sts premice. As we said in the introduction to this section, we will make sure that a stationary set of measurable cardinals will remain measurable in the model produced by our construction. This will be achieved by adding each such measures by brute force. 

By induction we define a sequence $\vec{\K}=(\K_\a: \a\in \Omega)$, called a $Z$-\textit{good stack}\footnote{There can be many such stacks.}, such that
\begin{enumerate}
\item $\K_0$ is a fully backgrounded $Z$-good sts premouse over $\R$,
\item  for every $\a$, $\K_{\a+1}$ is a fully backgrounded $Z$-good sts hod premouse over $\K_\a$,
\item if $\a$ is a limit ordinal and $Ord\cap \bigcup_{\b<\a}\K_\b\not \in S$ then $\K_\a=Lp^{Zv, sts}(\bigcup_{\b<\a}\K_\b)$,
\item if $\a$ is a limit ordinal and $\l=_{def}Ord\cap \bigcup_{\b<\a}\K_\b \in S$ then letting $U$ be a normal measure on $\l$ and setting $\K'_\a=\bigcup_{\b<\a}\K_\b$ and $\K''_\a=\pi_U(\K'_\a)|(\l^{++})^{\pi_U(\K'_\a)}$, $\K'''_\a=(\K''_\a, E)$ where $E$ is the $(\l, (\l^{+})^{\pi_U(\K'_\a)})$-extender derived from $\pi_U$ and $\K_\a$ is the core of $\K'''_\a$. 
\end{enumerate}
We call $\vec{\K}$ the $(f.b. Z)$-validated stack. The construction of $\vec{\K}$ is straightforward. However, we need to verify the following three statements. \\\\
(S1) For $\a<\b$, if $\d$ is a Woodin cardinal of $\K_\a$ then no level of $\K_\b$ projects across $\d$ and $\K_\b\models ``\d$ is a Woodin cardinal".\\
(S2) If $Ord\cap \bigcup_{\b<\l}\K_\b \in S$ then $\l$ is a measurable cardinal in $\K_{\l+1}$. \\
(S3) The class of $\l$ such that $Ord\cap \bigcup_{\b<\l}\K_\b \in S$ is stationary.\\\\
We now prove the above three clauses by proving a sequence of lemmas.

\begin{lemma}\label{iterability} For every $\a$, $\K_\a$ and  is a $Z$-validated mouse. If $\a\in S$ is such that $Ord\cap \bigcup_{\b<\l}\K_\b \in S$, then $\K'_\a$, $\K''_\a$, $\K'''_\a$ are also $Z$-validated mice. 
\end{lemma}

\rlem{iterability} is a consequence of \cite[Corollary 3.16]{Neeman}. This corollary shows that if $U$ is a good hull then the pre-images of the relevant objects have iteration strategies that pick realizable branches, which implies that they have $Z$-approved strategies.

\begin{lemma}\label{s1}(S1) holds.
\end{lemma}
\begin{proof} Fix $\a<\b$ and $\d$ as in the statement of (S1). Suppose $\M\insegeq \K_{\b}$ is such that $\rho(\M)< \d$. It follows from our construction that for some $\gg+1\leq \a$, $\d=\d^{\K_{\gg+1}}$. Let $p=p_{n+1}(\M)$ be the standard parameter of $\M$ and $n$ be least such that $\rho_{n+1}(\M)=\rho(\M)<\d$. Let $\W$ be the canonical decoding structure of $Hull_1^{\M^n}(\d\cup\{p\})$, where $\M^n$ is the $n$-th reduct of $\M$. As $\M$ is a $Z$-validated sts premouse and $\K_{\gg+1}$ is $Z$-suitable, we must have that $rud(\W)\models ``\d$ is a Woodin cardinal". Hence, $\rho(\W)=\d$, contradiction. A similar argument shows that $\K_\b\models ``\d$ is a Woodin cardinal". 
\end{proof}

\begin{lemma}\label{s2}(S2) holds.
\end{lemma}
\begin{proof} First we claim that $\rho(\K_\l)>\l$. \rlem{s1} shows that $\rho(\K_\l)\geq \l$. Assume then that $\rho(\K_\l)=\l$. Let $\W=Core(\K_\l)$ and $U$ be the normal measure on $\l$. Because of our definition of $\vec{\K}$, we have that 
\begin{center}
$Ult(V, U)\models \K''_\l=Lp^{Zv, sts}(\K'_\l)$. 
\end{center}
Let now $F$ be the last extender of $\W$. As $\card{\W}=\l$, we have $\sigma: Ult(\W, F)\rightarrow \pi_U(\W)$ such that $\sigma\in Ult(V, U)$. It follows that
\begin{center}
$Ult(V, U)\models  Ult(\W, F)$ is a $Z$-validated sts premouse over $\K'_\l$. 
\end{center}
Hence, $Ult(\W, F)\insegeq \K''_\l$ implying that $\W\in \K''_\l$. Thus, $\rho(\K_\l)>\l$. 

The same argument shows that if $\M\insegeq \K_{\l+1}$ then $\rho(\M)>\l$. Thus, $\l$ must be a measurable cardinal in $\K_{\l+1}$. 
\end{proof}

(S3) is trivial. It then follows that $\bigcup_{\a\in Ord}\K_\a$ is an almost excellent hybrid premouse. 

\subsection{The conclusion assuming $\neg\sf{Hypo}$}

We remind our reader that we have gotten to this point by assuming that $\neg\sf{Hypo}$ holds. The following summarizes the results of the previous subsection.  

\begin{corollary}\label{summary of chap 10} Assume $\neg\sf{Hypo}$. Then there is an honest $Z$-suitable $\R$ and a $Z$-validated almost excellent class size premouse $\K$ based on $\R$ satisfying the following conditions. 
\begin{enumerate}
\item  For each maximal window $w$ of $\K$ and for each $\eta\in (\nu^w, \d^w)$ that is a regular cardinal in $\K$, $\K$ has a $Z$-validated iteration strategy $\Sigma$ that acts on normal iterations that are based on $\K|\eta$ and are above $\nu^w$.
\item For each maximal window $w$ of $\K$, $\K|\d^w$ is a fully backgrounded $Z$-good sts premouse over $\K|\nu^w$.
\item For each Woodin cardinal $\d$ of $\K$ and for each $Z$-validated sound sts mouse $\M$ such that $\K|\d\insegeq \M$, $\M\models ``\d$ is a Woodin cardinal".
\end{enumerate}
\end{corollary}

The next proposition completes the proof  \rthm{thm:main_theorem} and \rthm{thm:tower_sealing} assuming $\neg\sf{Hypo}$. 

\begin{proposition}\label{conclusion from neg hypo} Assume $\neg\sf{Hypo}$. Then there is a class size excellent hybrid premouse.
\end{proposition}

We spend the rest of this subsection proving \rprop{conclusion from neg hypo}. Let $\R$ and $\K$ be as in \rcor{summary of chap 10}. We claim that in fact $\K$ is excellent. To see this let $w$ be a maximal window of $\K$ and let $\eta\in (\nu^w, \d^w)$ be a regular cardinal of $\K$. We want to see that in $\K$, $\K|\eta$ has an iteration strategy that acts on normal iterations that are above $\nu^w$. Let $\Sigma$ be the $Z$-validated strategy of $\K|\eta$ that acts on normal iterations that are above $\nu^w$. It is enough to show that $\Sigma\rest \K$ is definable over $\K$.

We work inside $\K$. Given a normal iteration $\T$ of $\K|\eta$ that is above $\nu^w$, we will say $\T$ has a \textit{correct} $\Q$-structure if letting $(u, \zeta, (\M_\xi, \N_\xi: \xi\leq\d^u))$ be such that
\begin{enumerate}
\item $u$ is the least maximal window of $\K$ with the property that $\T\in \K|\d^u$,  
\item $\zeta\in (\nu^u, \d^u)$ is such that $\T\in \K|\zeta$,
\item $(\M_\xi, \N_\xi: \xi<\d^u)$ are the models of the fully backgrounded $(\R, \R^b, S^\K)$-authenticated sts construction of $\K|\d^u$ done over $\m(\T)$ using extenders with critical points $>\zeta$,
\end{enumerate}
for some $\xi<\d^u$, $\M_\xi\models``\d(\T)$ is not a Woodin cardinal". We then say that $\M_\xi$ is the correct $\Q$-structure for $\T$. We have that $\M_\xi$ has a $Z$-validated iteration strategy, and hence if it exists it is unique (i.e. does not depend on $\zeta$).

Continuing our work in $\K$, given $\T$ as above we say $\T$ is \textit{correctly guided} if for every limit $\a<lh(\T)$, letting $b=[0, \a]_\T$, $\Q(b, \T\rest \a)$ exists and is the correct $\Q$-structure of $\T$. The following lemma finishes the proof of \rprop{conclusion from neg hypo}. 

\begin{lemma}\label{strategy capturing} Suppose $\T\in \K$ is a normal iteration of $\K|\eta$ of limit length that is according to $\Sigma$. Then $\T$ is correctly guided and if $\T$ is of limit length then $\T$ has a correct $\Q$-structure. 
\end{lemma}
\begin{proof} The second part of the conclusion of the lemma implies the first as we can apply it to the initial segments of $\T$.  Thus assume that $\T$ is correctly guided and is of limit length. Let $b=\Sigma(\T)$. Then $\Q(b, \T)$ exists and $Z$-validated. Set $\Q=\Q(b, \T)$.

Let $(u, \zeta, (\M_\xi, \N_\xi: \xi\leq\d^u))$ be as in the definition of the correct $\Q$-structure. Towards a contradiction assume that letting $\N=_{def}\M_{\d^u}$, $\N\models ``\d(\T)$ is a Woodin cardinal". Notice that $\K|\d^u$ is generic over $\N$, implying that we can translate $\K$  via $S$-constructions into an sts premouse over $\N$, call it $\K'$. We have that $\K'\models ``\d(\T)$ is a Woodin cardinal" and $\K'$ is almost excellent. 

Next we compare  $\Q$ with $\K'$. All of the extenders on the extender sequence of $\K'$ have fully backgrounded certificates,  which implies that in the aforementioned comparison only the $\Q$-side moves. Let $\Lambda$ be the unique $Z$-validated strategy of $\Q$ and let $\U$ be the tree on $\Q$ of limit length such that $\m(\U)=\K'|\d^u$. Set $c=\Lambda(\U)$ and $\M=\Q(c, \U)$\footnote{Notice that $\Q(c, \U)$ must exist as $\Q$ projects to $\d(\T)$.}.

Notice now that $\M$ is $\d^u$-sound, $\d^u$ is a cutpoint in $\M$ and $\M$ has no extenders with critical point $\d^u$. Moreover, $\M$ is a $Z$-validated sts mouse over $\K'|\d^u$, and therefore, it can be translated into a $Z$-validated sts mouse $\X$ over $\K|\d^u$. We must then have that $rud(\X)\models ``\d^u$ is a Woodin cardinal". But then $rud(\M)\models ``\d^u$ is a Woodin cardinal", contradiction. 

\end{proof}

Notice now that for $\T\in \K$ we have the following equivalences.
\begin{enumerate}
\item $\T\in dom(\Sigma)$ if and only if $\K\models ``\T$ is correctly guided".
\item $\Sigma(\T)=b$ if and only if $\Q(b, \T)$ exists and $\K\models ``\Q(b, \T)$ is the correct $\Q$-structure for $\T$".
\end{enumerate}

\subsection{Excellent hybrid premouse from $\sf{Hypo}$}

Finally we show how to get an excellent hybrid premouse from $\sf{Hypo}$. This will complete the proof of \rthm{thm:main_theorem} and \rthm{thm:tower_sealing}. Suppose then $\R$ is an honest $Z$-suitable hod premouse, $X$ is a set such that $\R\in X$ and $\M$ is an almost excellent $Z$-validated sts premouse over $X$. In particular, $\M$ is a model of $\sf{ZFC}$. It must be clear from our construction that we can assume that $X$ has a well-ordering in $L[X]$. Let then $\k=(\card{X}^+)^\M$. Let $\d$ be the least Woodin cardinal of $\M$ that is $>\k$, and let $\N$ be the output of the fully backgrounded $(\R, \R^b, S^\M)$-authenticated sts construction of $\M|\d$ done over $\R$. Once again using $S$-constructions, we can translate $\M$ over to an sts mouse $\P$ over $\N$ such that $\P[\M|\d]=\M$. Moreover, $\P$ is almost excellent and is a $Z$-validated sts mouse over $\R$. Because good hulls of $\P$ are iterable via a $Z$-approved strategy, we can assume, by minimizing if necessary, that $\P$ is minimal in the sense that for each $\eta\in (\d^\R, Ord\cap \P)$, $\P|\eta$ is not an almost excellent sts premouse over $\R$.  

The rest of the proof follows the same argument as the one given in the previous subsection. We show that $\P$ is in fact excellent. As before, this amounts to showing that for an window $w$ of $\P$ such that $\d^w>\d^\R$, and for any $\eta\in (\nu^w, \d^w)$, if $\eta$ is a regular cardinal of $\P$ then $\P\models ``\P|\eta$ has $Ord$-iteration strategy". Towards a contradiction assume not.

Let $U$ be a good hull such that $\P\in U$. Set $\S=\P_U$, $\W=\R_U$ and $\l=\eta_U$. Let $\Lambda$ be the $Z$-approved strategy of $\S$ and let $\Sigma$ be the fragment of $\Lambda$ that acts on normal iterations of $\S|\l$ that are above $\nu^{w_U}$. The following lemma can be proved via a proof almost identical to the proof of \rlem{strategy capturing}. We define correct $\Q$-structure and correctly guided exactly the same way as we defined them in the previous subsection, except the definition now takes place in $\S$. 

\begin{lemma}\label{strategy capturing 2} Suppose $\T\in \S$ is a normal iteration of $\S|\eta$ of limit length that is according to $\Sigma$. Then $\T$ is correctly guided and if $\T$ is of limit length then $\T$ has a correct $\Q$-structure. 
\end{lemma}
There is only one difference between the proofs of \rlem{strategy capturing} and \rlem{strategy capturing 2}. In the proof of \rlem{strategy capturing}, we concluded that $rud(\X)\models ``\d^u$ is a Woodin cardinal" using the fact that $(\K|\d^u)^\#$ is $Z$-suitable sts premouse. Here we no longer have such a fact, but here we can use minimality of $\P$ to derive the same conclusion. 

As in the previous subsection, \rlem{strategy capturing 2} easily implies that $\Sigma\rest \P$ is definable over $\P$. This completes the proof of \rthm{thm:main_theorem} and \rthm{thm:tower_sealing}.

\section{Open problems and questions}

The rather mild assumption that the class of measurable cardinals is stationary is used in various ``pressing down" arguments in the proof of Theorem \ref{thm:main_theorem} and \rthm{thm:tower_sealing}, and also in stabilization arguments like those of \rthm{main thm on stacks under sealing}. This assumption is probably not needed, though proving some sort of stabilization lemma like the aforementioned one is probably necessary. 

\begin{question} Are the following theories equiconsistent?
\begin{enumerate}
\item $\sf{Sealing}\ $ + ``There is a proper class of Woodin cardinals".
\item $\sf{LSA-over-uB}\ $ + ``There is a proper class of Woodin cardinals". 
\item $\sf{Tower \ Sealing}\ $ + ``There is a proper class of Woodin cardinals".
\end{enumerate}
\end{question}

As mentioned above, CMI becomes very difficult past $\sf{Sealing}$. A good test question for CMI practitioners is.
\begin{open}
Prove that Con$(\sf{PFA})$ implies Con($\sf{WLW}$).
\end{open}

We know from the results above that $\sf{WLW}$ is stronger than $\sf{Sealing}$ and is roughly the strongest natural theory at the limit of traditional methods for proving iterability. We believe it is plausible to develop CMI methods for obtaining canonical models of $\sf{WLW}$ from just $\sf{PFA}$.\footnote{The second author observes that assuming $\sf{PFA}$ and there is a Woodin cardinal, there is a canonical model of $\sf{WLW}$. The proof is not via CMI methods, but just an observation that the full-backgrounded construction as done in \cite{Neeman} reaches a model of $\sf{WLW}$. The Woodin cardinal assumption is important here. The argument would not work if one assumes just $\sf{PFA}$ and/or a large cardinal milder than a Woodin cardinal, e.g. a measurable cardinal or a strong cardinal.}

\begin{remark}\label{sealing+}
\begin{enumerate}
\item By the above discussion, we also get in $V=\P[g]$ that for every generic $h$, $(\Gamma^\infty_h)^\#$ exists and by Lemma \ref{lem:Gamma}, $L(\Gamma^\infty_h) \vDash \sf{AD}_\mathbb{R} + \Theta$ is regular. So we have the following strengthening of $\sf{Sealing}$: $\sf{Sealing}^+ $ for all $V$-generic $g$, in $V[g]$,  $L(\Gamma^\infty,\mathbb{R}) \vDash \sf{AD}_\mathbb{R} + \Theta$ is regular. We call this theory $\sf{Sealing}^+$.
\item The conclusion of $\sf{LSA-over-uB}$ can be weakened to: For all $V$-generic $g$ there is $A\subseteq \bR^{V[g]}$ such that $L(A, \bR^{V[g]})\models \sf{LSA}$ and $\Gamma^\infty_g$ is contained in $L(A, \bR^{V[g]})$. We call this theory $\sf{LSA-over-uB}^-$. 
\item The results of this paper show the following. Let $T=``$there exists a proper class of Woodin cardinals and the class of measurable cardinals is stationary". Then the following theories are equiconsistent:
\begin{enumerate}
\item $\sf{Sealing}+T$.
\item $\sf{Sealing}^++T$.
\item $\sf{Tower \ Sealing} + T$.
\item $\sf{LSA-over-uB}+T$.
\item $\sf{LSA-over-uB}^-+T$.
\item $\sf{Weak \ Sealing} + T$.
\item $\sf{Sealing^-} + T$
\end{enumerate}
\end{enumerate}
\end{remark}

We end  the paper with the following conjecture, if true, would be an ultimate analog of the main result of \cite{steel2002core}.

\begin{conjecture}\label{equivalence conjecture} Suppose there are unboundedly many Woodin cardinals and the class of measurable cardinals is stationary. Then the following are equivalent.
\begin{enumerate}
\item $\sf{Sealing}$.
\item $\sf{Sealing^+}$.
\item $\sf{Weak \ Sealing}$.
\item $\sf{Sealing^-}$.
\item $\sf{Tower \ Sealing}$.
\end{enumerate}
\end{conjecture}

\bibliographystyle{plain}
\bibliography{Gamma_uB2023.bib}
\end{document}